\documentclass[11pt]{article}

\usepackage[T1]{fontenc}
\usepackage[utf8]{inputenc}
\usepackage{textcomp}

\usepackage[margin=1in]{geometry}
\usepackage{mathptmx}

\usepackage{xcolor}
\usepackage{tcolorbox}
\usepackage{amsmath,amssymb}
\usepackage{stmaryrd}
\usepackage{amsthm}
\usepackage{indentfirst}
\let\LocNoIndent\noindent

\makeatletter
\renewcommand\section{\@startsection{section}{1}{\z@}{4.0ex plus 1.2ex minus .3ex}{1.4ex plus .3ex}{\normalfont\large\bfseries}}
\renewcommand\subsection{\@startsection{subsection}{2}{\z@}{3.25ex plus 1ex minus .25ex}{1.15ex plus .25ex}{\normalfont\normalsize\bfseries}}
\renewcommand\subsubsection{\@startsection{subsubsection}{3}{\z@}{2.5ex plus .8ex minus .2ex}{0.9ex plus .2ex}{\normalfont\normalsize\itshape}}
\renewcommand\paragraph{\@startsection{paragraph}{4}{\z@}{2.0ex plus .6ex minus .15ex}{0.75em}{\normalfont\normalsize\bfseries}}
\renewcommand\l@section{\@dottedtocline{1}{0em}{1.5em}}
\renewcommand\l@subsection{\@dottedtocline{2}{0em}{2.3em}}
\renewcommand\l@subsubsection{\@dottedtocline{3}{0em}{3.2em}}
\makeatother

\usepackage{longtable,booktabs,array}
\usepackage{multirow}
\usepackage{calc}
\usepackage{float}
\usepackage{etoolbox}
\makeatletter
\patchcmd\longtable{\par}{\if@noskipsec\mbox{}\fi\par}{}{}
\makeatother
\IfFileExists{footnotehyper.sty}{\usepackage{footnotehyper}}{\usepackage{footnote}}
\makesavenoteenv{longtable}

\usepackage{graphicx}
\IfFileExists{caption.sty}{\usepackage[labelsep=period]{caption}}{}
\IfFileExists{subcaption.sty}{\usepackage{subcaption}}{}
\usepackage{tikz}
\IfFileExists{xurl.sty}{\usepackage{xurl}}{}
\IfFileExists{natbib.sty}{%
  \usepackage[numbers,sort,compress]{natbib}%
}{%
  \IfFileExists{cite.sty}{\usepackage{cite}}{}%
}
\usepackage{hyperref}
\hypersetup{hidelinks}
\usepackage{bookmark}

\makeatletter
\providecommand{\BKM@entry}[2]{}
\makeatother

\makeatletter
\def\BKM@entry#1#2{}
\makeatother

\makeatletter
\newsavebox\pandoc@box
\newcommand*\pandocbounded[1]{%
  \sbox\pandoc@box{#1}%
  \Gscale@div\@tempa{\textheight}{\dimexpr\ht\pandoc@box+\dp\pandoc@box\relax}%
  \Gscale@div\@tempb{\linewidth}{\wd\pandoc@box}%
  \ifdim\@tempb\p@<\@tempa\p@\let\@tempa\@tempb\fi%
  \ifdim\@tempa\p@<\p@\scalebox{\@tempa}{\usebox\pandoc@box}%
  \else\usebox{\pandoc@box}%
  \fi%
}
\def\fps@figure{htbp}
\makeatother

\newtheorem{remark}{Remark}

\title{Mesh-Intrinsic GFEM:\\[6pt] Asymptotic Smoothness on $C^0$ Unstructured Meshes\\[6pt]}

\author{Rong Tian\\Institute of Applied Physics and Computational Mathematics, Beijing 100086, China}
\date{}

\newcommand{\captionline}[1]{\par\smallskip\LocNoIndent\hfill\begin{minipage}{0.95\linewidth}\small #1\end{minipage}\hfill\par\smallskip}
\newcommand{\keywords}[1]{\par\LocNoIndent\textbf{Keywords:} #1\par}

\begin{document}
\maketitle

\begin{abstract}
\noindent
The Mesh-Intrinsic Generalized Finite Element Method (MiGFEM) constructs enriched local approximations from the original nodal degrees of freedom of a finite element mesh, without introducing extra enrichment unknowns. It therefore avoids the linear dependence, ill-conditioning, energy inconsistency, and mass-lumping difficulties that persist in conventional enriched formulations.

This paper establishes an additional property of MiGFEM: \emph{asymptotic smoothness} on standard $C^0$ unstructured meshes. Although the global approximation is assembled from $C^0$ partition-of-unity functions, its inter-element derivative jumps satisfy
\[
\|\llbracket D^\alpha u^h\rrbracket\| = \mathcal{O}(h^{p+1-|\alpha|}),
\qquad |\alpha|\le p,
\]
for $u\in C^{p+1}(\Omega)$, and vanish identically when the exact solution belongs to the local reconstruction space on every patch. The mechanism rests on two ingredients: \emph{interface coherence} of the patchwise local approximations, and a derivative-space partition-of-unity identity (Partition of Zero) that cancels the coherent interfacial traces.

Asymptotic smoothness motivates the \emph{leading Leibniz derivative}, which retains only the principal term of the full Leibniz expansion of $D^\alpha u^h$. This operator is globally continuous by construction; the neglected remainder is of order $\mathcal{O}(h^{p+1-|\alpha|})$. Consequently, the same $C^0$ MiGFEM trial space supports both weak Galerkin discretization via the full Leibniz derivative and strong-form collocation via the leading Leibniz derivative, on identical nodal unknowns.

Numerical experiments verify the predicted jump decay and confirm exact cancellation under patchwise reproduction. The unified weak--strong framework is demonstrated on elasticity and biharmonic problems using unstructured meshes. A singular-enrichment test further shows that strong-form collocation reproduces enriched basis functions exactly in the absence of quadrature constraints.
\end{abstract}


\section{Introduction}

Standard $C^0$ finite elements remain the dominant discretization technology in CAE industry because they combine robustness, geometric flexibility on unstructured meshes, and compatibility with mature assembly and solver infrastructures \cite{ref1,ref2,ref3,ref40,ref48}. At the same time, many important classes of problems exhibit local solution features---high gradients, singular fields, boundary layers, oscillations, or known analytical patterns---that are inefficient to capture with low-order polynomial bases alone. This observation has motivated the development of enrichment-based methods, including the Numerical Manifold Method (NMM) \cite{ref20,ref21}, Partition-of-Unity methods (PUM/PUFEM) \cite{ref22,ref52}, the Generalized Finite Element Method (GFEM) \cite{ref23,ref24,ref53}, and the eXtended Finite Element Method (XFEM) \cite{ref25,ref26,FriesBelytschko2010XFEMoverview}. Their common principle is to use a partition of unity (PoU) to blend problem-tailored local approximation spaces into a globally conforming approximation without changing the underlying mesh topology.

This PoU paradigm is powerful because it separates mesh topology from approximation richness. Once suitable local enrichments are chosen, one can substantially improve approximation quality without uniform mesh refinement. In particular, local trial spaces may incorporate high-order polynomials, singular functions, oscillatory functions, crack-tip asymptotics, or other basis functions reflecting the local structure of the exact solution. In principle, this gives GFEM/XFEM a much broader approximation capability than standard FEM.

However, the classical PoU enrichment framework is accompanied by several persistent practical difficulties.
The first is the well-known linear dependence or near-linear dependence problem, which causes severe ill-conditioning of the global stiffness matrix \cite{ref22,ref52,ref23,ref24,ref53}. This issue is intrinsic to the coexistence of the PoU basis and enriched basis functions and becomes more pronounced when enrichments overlap strongly or reproduce low-order polynomial content already contained in the underlying finite element space. Various approaches have been proposed to address the linear dependence problem \cite{ref43,ref45,ref44,ref46,RajendranZhang2007,OhKimHong2008,ZhangRajendran2008,CaiZhuangAugarde2010,ref55,XuRajendran2011,ref56,ref57,ref47,TianWen2015XFEM,GuptaDuarteBabuskaBanerjee2015SGFEM,ref58,ref59,ref60,Cui2020SGFEMelasticity}. To the best of the author's knowledge, however, the difficulty is still more numerically alleviated than theoretically solved.

A second difficulty arises from the additional degrees of freedom introduced by classical enrichment. In static analysis, these extra unknowns enlarge the global system and complicate implementation. In dynamics, they create more fundamental issues. In XFEM for propagating cracks, for example, the set of enrichment unknowns may change from one time step to the next as the crack advances \cite{BelytschkoChen2003,RethoreGravouil2005,SongAreiasBelytschko2006,Elguedj2018}. This disrupts the continuity of information transfer among nodal unknowns and has long been associated with the energy inconsistency observed in dynamic crack growth simulations \cite{RethoreGravouil2005,GregoireMaigreRethore2007,GravouilElguedj2009}. Historically, this difficulty partly motivated reduced crack-tip enrichment strategies and the preference for cohesive-zone-type crack-tip treatements in commercial XFEM implementations for propagating cracks. The same extra-DoF structure also complicates mass lumping \cite{GregoireMaigreRethore2007,MenouillardRethore2008,GravouilElguedj2009} and partially helps explain why GFEM/XFEM is rarely adopted in explicit dynamics \cite{WenTian2016XFEM,Elguedj2018}.

A third limitation concerns numerical efficiency. Classical GFEM most naturally uses the standard finite element shape functions as the PoU basis, which is elegant because the method reduces to standard FEM when the local approximation is taken as constant. Yet this same $C^0$ PoU structure ties the formulation naturally to weak forms and numerical integration. When enrichments are non-polynomial, discontinuous across internal subcells, or singular, accurate quadrature often requires heavy element subdivision, special integration rules, or mesh cutting. In such cases, the theoretical approximation power of local enrichment is offset by substantial implementation and computational overhead. From an industrial viewpoint, this is a major barrier to broader deployment. 

In parallel with enrichment-based methods, alternative discretization paradigms that avoid element-based quadrature have been developed, including meshfree methods based on moving least-squares or reproducing kernel approximations \cite{ref31,ref32,ref29,ref30,ref51}, radial basis function (RBF) methods \cite{ref33,ref34,ref35,ref36,ref37,ref51}, generalized finite difference schemes \cite{PerroneKao1975,LiszkaOrkisz1980,BenitoUrenaGavete2001,FornbergLehtoPowell2013,FlyerFornbergBayonaBarnett2016,BayonaFlyerFornbergBarnett2017}, the meshless generalized finite difference method \cite{TinocoGuerreroDominguezMotaGuzmanTorresPedrazaJimenezTinocoRuiz2025}, and strong-form meshfree frameworks \cite{SlakKosec2021Medusa}. Relatedly, alternative discretization frameworks such as the Virtual Element Method (VEM) \cite{refVEM2013basic,refVEM2016Ciarlet} and the Hybrid High-Order (HHO) method \cite{refHHO2015} have been developed for general polygonal/polyhedral meshes.

The present work revisits the extra-DoF-free GFEM in \cite{ref47} and Improved XFEM in \cite{TianWen2015XFEM,WenTian2016XFEM}, referred to here as \emph{mesh-intrinsic enrichment} (MiGFEM). On each nodal patch, a local enriched approximation is constructed from the original finite element nodal values only and blended through the standard $C^0$ PoU. No additional global unknowns are introduced; the enrichment is intrinsic to the mesh data structure. Consequently, the method preserves the original nodal layout while avoiding the three difficulties described above---linear dependence, extra-DoF complications in dynamics, and quadrature overhead \cite{TianWen2015XFEM,WenTian2016XFEM,TianWenWang2019IXFEM,XiaoWenTian2021IXFEM,XiaoWenTianZhang2023IXFEM,WenTianWangFeng2023IXFEM,WangWenTianFeng2024IXFEM,ref50}.

This paper identifies a property not previously emphasized: \emph{asymptotic smoothness on standard $C^0$ unstructured meshes}. Although the global MiGFEM approximation is assembled from a $C^0$ PoU, its inter-element derivative jumps satisfy
\[
\|\llbracket D^\alpha u^h\rrbracket\| = \mathcal O(h^{p+1-|\alpha|}),
\qquad |\alpha|\le p,
\]
for $u\in C^{p+1}(\Omega)$, and vanish identically under patchwise exact reproduction. The mechanism rests on two ingredients. (i) A derivative-space PoU identity---the \emph{Partition of Zero} (PoZ): $\sum_i D^\alpha N_i\equiv 0$ for $|\alpha|\ge 1$. (ii) \emph{Interface coherence} of the local patch reconstructions: when neighboring $U_i$ share interfacial jets, PoZ cancels the corresponding derivative jumps of the blended global approximation. For general smooth fields, the local reconstructions approximate a common exact solution on overlapping patches, inducing approximate coherence and the decay rate above. This smoothness motivates the \emph{leading Leibniz derivative}, $\widetilde D^\alpha u^h = \sum_i N_i D^\alpha U_i$, which retains the principal term of the Leibniz expansion and is globally continuous by construction; the remainder is $\mathcal{O}(h^{p+1-|\alpha|})$ (Appendix~\ref{app:leibniz_remainder}). Consequently, the same $C^0$ MiGFEM trial space supports weak Galerkin discretization via the full Leibniz derivative and strong-form collocation via the leading derivative, on identical nodal unknowns.

This result has two important implications: First, asymptotic smoothness explains why mesh-intrinsic enrichment exhibits stronger inter-element regularity than standard $C^0$ approximations. Second, it enables enriched strong-form discretization without quadrature or mesh-cutting costs. Numerical results confirm high accuracy for both weak and strong formulations on the same approximation space (Section~\ref{sec:numerical_verification}).

The remainder of the paper is organized as follows. Section~\ref{sec:migfem_approximation} introduces the mesh-intrinsic enrichment formulation and the associated PoU-blended global approximation. Section~\ref{sec:smoothness_mechanism} develops the asymptotic smoothness mechanism through PoZ, interface coherence, and inter-element jump analysis. Section~\ref{sec:leibniz_leading_derivative} defines the leading Leibniz derivative and relates it to the full product-rule derivative. Section~\ref{sec:unified_discretizations} formulates the weak Galerkin and strong collocation discretizations on the same mesh-intrinsic trial space, including the NC, CC, and SD variants. Section~\ref{sec:numerical_verification} reports numerical tests on jump decay, polynomial reproduction, elasticity, biharmonic problems, and singular enrichment. Section~\ref{sec:conclusions} concludes the paper.

Supporting materials are provided in appendices: Appendix~\ref{app:cwls_kronecker} describes the CWLS patch reconstruction and nodal Kronecker property; Appendix~\ref{app:jump_decay_proof} gives proof details for the inter-element jump bounds in Section~\ref{sec:asymp_jump_taylor}; Appendix~\ref{app:leibniz_remainder} gives the Leibniz remainder estimate; Appendix~\ref{app:biharmonic_order_gap} reports the nodal-injection diagnostic for the biharmonic trial space; Appendix~\ref{subsec:mi_pfem_pu_mini_case} compares MiGFEM with standard $p$-FEM and PUFEM.

\section{Mesh-intrinsic Enrichment---Local Enrichment Without Extra DoFs}
\label{sec:migfem_approximation}

The key idea of MiGFEM is to enrich the approximation \emph{locally on each mesh patch}, while retaining the \emph{original nodal degrees of freedom} of the underlying finite element mesh. In this sense, the enrichment is intrinsic to the mesh itself: no additional global unknowns are introduced, and the enriched approximation remains assembled from the standard nodal data.

Let $\mathcal T_h$ be a conforming mesh of the computational domain $\Omega$, and let $\mathcal N$ denote the set of mesh nodes. For each node $i\in\mathcal N$, we associate a patch $\omega_i$, typically defined as the union of elements sharing node $i$. Let $N_i(x)$ be the standard $C^0$ finite element partition-of-unity (PoU) shape function attached to node $i$, satisfying
\begin{equation}
\sum_{i\in\mathcal N} N_i(x)=1,\qquad x\in\Omega.
\end{equation}

On each patch $\omega_i$, we construct a local enriched approximation $U_i(x)$ from the nodal values in that patch by a weighted least-squares reconstruction over a prescribed polynomial or enriched local trial space. Denote by
\begin{equation}
\mathcal S_i \subset \mathcal N
\end{equation}
the set of nodal indices used in the reconstruction on $\omega_i$. Then the local approximant can be written in the nodal-value form
\begin{equation}
\label{eq:Ui_local_expansion}
U_i(x)=\sum_{k\in\mathcal S_i}\Psi_k^{(i)}(x)\,u_k,
\end{equation}
where $\Psi_k^{(i)}(x)$ are the patch-dependent reconstruction shape functions and $\{u_k\}$ are the original nodal unknowns. The crucial point is that the coefficients of $U_i$ are determined entirely from the existing nodal values on the patch; hence no generalized enrichment DoFs are introduced.

The global MiGFEM approximation is then defined by blending the local patch reconstructions with the standard finite element PoU:
\begin{equation}
\label{eq:migfem_blending}
u^h(x)=\sum_{i\in\mathcal N} N_i(x)\,U_i(x).
\end{equation}
Substituting \eqref{eq:Ui_local_expansion} into \eqref{eq:migfem_blending} gives
\begin{equation}
u^h(x)
=
\sum_{i\in\mathcal N} N_i(x)\sum_{k\in\mathcal S_i}\Psi_k^{(i)}(x)\,u_k
=
\sum_{k\in\mathcal N}\phi_k(x)\,u_k,
\end{equation}
where the global MiGFEM shape function associated with node $k$ is
\begin{equation}
\label{eq:phi_k_def}
\phi_k(x)=\sum_{i\in\mathcal N_k}N_i(x)\Psi_k^{(i)}(x),
\qquad
\mathcal N_k:=\{\,i\in\mathcal N:\,k\in\mathcal S_i\,\}.
\end{equation}
Therefore, although the enrichment is performed locally and patchwise, the resulting approximation still admits a standard FEM-like global expansion in terms of the original nodal DoFs only. The global MiGFEM trial space is then
\begin{equation}
\mathcal V^h := \operatorname{span}\{\phi_k\}_{k\in\mathcal N}.
\label{eq:global_trial_space}
\end{equation}

This feature distinguishes MiGFEM from conventional GFEM and PUFEM formulations. In standard GFEM, enrichment functions are typically introduced together with additional nodal or elemental unknowns, which may enlarge the global algebraic system and can lead to conditioning or linear-dependence issues \cite{ref43,ref47}. In contrast, MiGFEM embeds the enrichment directly into the local reconstruction process and reuses the existing nodal DoFs. The method therefore preserves the compact unknown structure of the original FEM while still enhancing the local approximation space.

Another important consequence of \eqref{eq:migfem_blending} is that neighboring patches are not independent. Since adjacent local reconstructions are all built from overlapping sets of nodal data and are blended through the same PoU, the resulting global approximation inherits a built-in interface coordination mechanism. This mechanism does not generally enforce exact higher-order continuity at a fixed mesh size, but it promotes cancellation of inter-element derivative jumps and leads to the asymptotic smoothness behavior analyzed in Section~\ref{sec:smoothness_mechanism}.

\subsubsection{Kronecker delta property}
\label{sec:properties_approx}

Provided that each local reconstruction is constrained to interpolate the patch-star node---the patch's topological center, namely
\begin{equation}
U_i(x_i)=u_i,
\end{equation}
the global MiGFEM approximation preserves the nodal interpolation property of the underlying finite element mesh. Indeed, using the Kronecker delta property of the standard PoU basis, $N_i(x_j)=\delta_{ij}$, we obtain
\begin{equation}
\label{eq:kronecker_from_cwls}
u^h(x_j)
=
\sum_{i\in\mathcal N}N_i(x_j)U_i(x_j)
=
U_j(x_j)
=
u_j.
\end{equation}
Hence MiGFEM remains nodally interpolatory whenever the patchwise reconstruction is enforced to pass through the corresponding patch-star node.

\begin{figure}[htbp]
  \centering
  \includegraphics[width=0.72\linewidth]{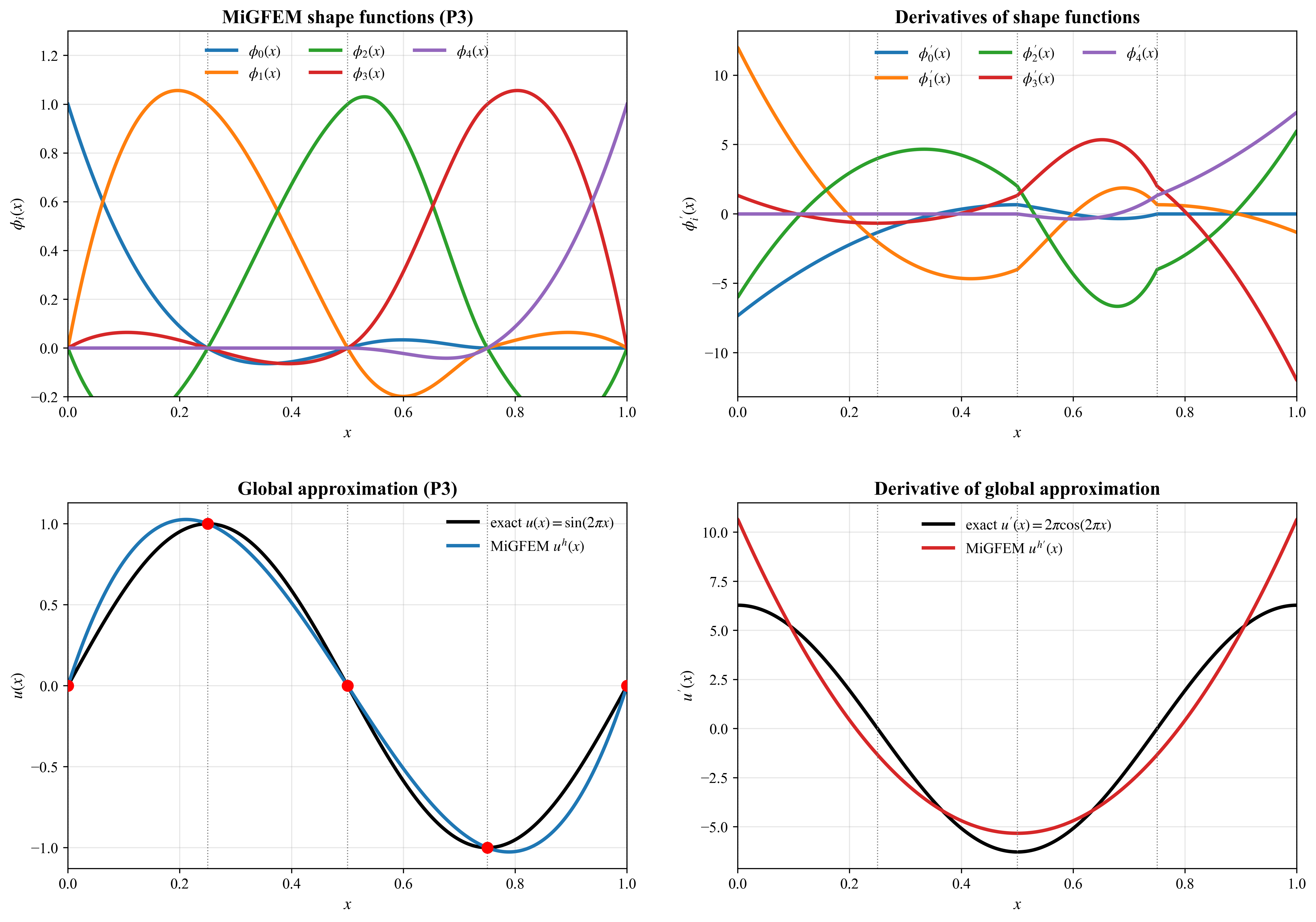}
  \caption{One-dimensional MiGFEM illustration. Upper row: MiGFEM shape functions $\phi_k$ and their derivatives. Lower row: global MiGFEM approximation $u^h(x)$ and its derivative compared with the exact solution $u(x)=\sin(2\pi x)$.}
  \label{fig:method_1d}
\end{figure}

Figure~\ref{fig:method_1d} illustrates the resulting enriched basis for the case of polynomial degree $p=3$. The upper row displays the MiGFEM global shape functions $\phi_k$ and their derivatives, while the lower row shows the corresponding approximation $u^h(x)=\sum_k\phi_k(x)u_k$ and its derivative. The approximation uses only the original nodal unknowns $\{u_k\}$, yet the reconstructed patch fields provide substantially richer local approximation behavior than standard linear FEM. This combination of local enrichment and unchanged global DoF structure is the defining feature of MiGFEM.

\section{Asymptotic Smoothness Mechanism}
\label{sec:smoothness_mechanism}

This section explains why MiGFEM exhibits improved inter-element regularity, despite being constructed from standard $C^0$ finite element partition-of-unity functions. The key mechanism is that derivative jumps across element interfaces are not governed solely by the continuity of the partition-of-unity functions, but also by the consistency and interface coherence of the patchwise local approximants. The resulting cancellation is encoded in a derivative-space partition-of-unity identity, which we call the Partition of Zero (PoZ).

\subsection{PoZ identity}

Let $\{N_i\}_{i\in\mathcal N}$ be a standard partition of unity, so that
\[
\sum_{i\in\mathcal N} N_i(x)=1 \qquad \text{for all } x\in\Omega.
\]
Differentiating this identity yields, for every multi-index $\alpha$ with $|\alpha|\ge 1$,
\begin{equation}
\label{eq:poz_def}
\sum_{i\in\mathcal N} D^\alpha N_i(x)=0,
\end{equation}
where the identity holds elementwise, and therefore also on each side of an internal interface wherever the traces are defined. We refer to \eqref{eq:poz_def} as the \emph{Partition of Zero} (PoZ), since the derivatives of the partition-of-unity functions sum to zero.

Let $\Gamma$ be an internal interface shared by two adjacent elements $E^+$ and $E^-$. Denote by
\[
\mathcal I(\Gamma):=\{\,i\in\mathcal N:\Gamma\subset \operatorname{supp}(N_i)\,\}
\]
the set of active partition-of-unity indices on $\Gamma$. Since $N_i$ vanishes in a neighborhood of $\Gamma$ whenever $i\notin\mathcal I(\Gamma)$, taking the traces of \eqref{eq:poz_def} from $E^+$ and $E^-$ and subtracting gives the interface form of PoZ:
\begin{equation}
\label{eq:poz_jump}
\sum_{i\in\mathcal I(\Gamma)} \llbracket D^\alpha N_i\rrbracket_\Gamma = 0,
\qquad |\alpha|\ge 1.
\end{equation}
This identity is the algebraic source of the cancellation mechanism developed below.

\subsection{Jump decomposition and cancellation}

From the MiGFEM representation $u^h=\sum_i N_i U_i$, the multi-index Leibniz rule gives, on each element $E$,
\begin{equation}
\label{eq:tag18}
D^{\alpha}u^{h}
=
\sum_{i\in\mathcal I(E)}\sum_{0\le \beta\le \alpha}
\binom{\alpha}{\beta}
(D^\beta N_i)(D^{\alpha-\beta}U_i),
\end{equation}
where $\mathcal I(E):=\{i\in\mathcal N:\operatorname{supp}(N_i)\cap E\neq\emptyset\}$.

Fix an internal interface $\Gamma$ shared by $E^+$ and $E^-$, and let $\mathcal I(\Gamma)$ be the active PoU index set on $\Gamma$. Since each $U_i$ is a single smooth field on its patch $\omega_i$, we have $\llbracket D^\mu U_i\rrbracket_\Gamma=0$ for $|\mu|\le p$. The jump of the product then factors as
\[
\llbracket (D^\beta N_i)(D^{\alpha-\beta}U_i)\rrbracket_\Gamma
=
\llbracket D^\beta N_i\rrbracket_\Gamma\,(D^{\alpha-\beta}U_i)|_\Gamma.
\]

Restricting the sum to $\mathcal I(\Gamma)$ and factoring out the $U_i$ traces yields the exact jump decomposition
\begin{equation}
\label{eq:tag19}
\llbracket D^\alpha u^h\rrbracket_\Gamma
=
\sum_{i\in\mathcal I(\Gamma)}\sum_{0\le \beta\le \alpha}
\binom{\alpha}{\beta}
\llbracket D^\beta N_i\rrbracket_\Gamma\,
(D^{\alpha-\beta}U_i)|_\Gamma.
\end{equation}
Isolating the $\beta=\alpha$ term gives the two-part form used later:
\begin{equation}
\label{eq:tag21}
\llbracket D^\alpha u^h\rrbracket_\Gamma
=
\underbrace{\sum_{i}\llbracket D^\alpha N_i\rrbracket_\Gamma\,U_i|_\Gamma}_{\text{Term I}}
\;+\;
\underbrace{\sum_{i}\sum_{0<\beta<\alpha}
\binom{\alpha}{\beta}
\llbracket D^\beta N_i\rrbracket_\Gamma\,(D^{\alpha-\beta}U_i)|_\Gamma}_{\text{Term II}}.
\end{equation}

\subsubsection{Exact cancellation under interface coherence}

Let the local approximants be \emph{$k$-th order derivative-coherent} on $\Gamma$:
\[
D^\mu U_i|_\Gamma = D^\mu U_j|_\Gamma =: D^\mu U_\Gamma,
\qquad i,j\in\mathcal I(\Gamma),\ |\mu|\le k.
\]
Then for $|\alpha|\le k+1$, the common traces factor out of \eqref{eq:tag21}:
\begin{equation}
\label{eq:jump_factorization}
\llbracket D^\alpha u^h\rrbracket_\Gamma
=
U_\Gamma\sum_{i}\llbracket D^\alpha N_i\rrbracket_\Gamma
\;+\;
\sum_{0<\beta<\alpha}
\binom{\alpha}{\beta}(D^{\alpha-\beta}U_\Gamma)
\sum_{i}\llbracket D^\beta N_i\rrbracket_\Gamma.
\end{equation}
Each sum $\sum_i\llbracket D^\beta N_i\rrbracket_\Gamma$ vanishes by the interface PoZ identity \eqref{eq:poz_jump}, hence
\begin{equation}
\label{eq:exact_jump_cancellation}
\llbracket D^\alpha u^h\rrbracket_\Gamma = 0,
\qquad |\alpha|\le k+1.
\end{equation}

Patchwise exact reproduction $U_i\equiv u$ ($k$ arbitrarily large) is the benchmark case; several settings that realize exact interface coherence at finite $h$ are catalogued in Section~\ref{sec:exact_regularity_mechanisms}. For a general $C^{p+1}$ field not contained in $\mathcal V_i$, the exact solution $u$ is single-valued; its contribution to the jump factors as in \eqref{eq:jump_factorization} and cancels by PoZ, leaving only the reconstruction-error contribution $J_e$ (explicit form in Appendix~\ref{app:jump_decay_proof}).

\subsection{\texorpdfstring{Asymptotic jump decay for smooth fields}{Asymptotic jump decay for smooth fields}}
\label{sec:asymp_jump_taylor}

Assume a shape-regular family of simplicial meshes $\mathcal T_h$ with characteristic size $h$; $\{N_i\}$ are the standard $C^0$ linear FE shape functions; on each nodal patch $\omega_i$, the local approximant $U_i\in\mathcal V_i\supset\mathcal P_p$ is obtained by the CWLS reconstruction of Appendix~\ref{app:cwls_kronecker}. Let $u\in C^{p+1}(\Omega)$. Decompose the local approximant as
\begin{equation}
\label{eq:local_approx_split}
U_i = u + e_i,
\qquad e_i := U_i - u.
\end{equation}
On each element $E\subset\operatorname{supp}(N_i)$, the inverse estimate $\|D^\beta N_i\|_{L^\infty(E)}\le C_N h^{-|\beta|}$ holds \cite{ref40,ref48}. Under standard MLS/WLS theory \cite{ref41,ref42,Wendland2004}, the reconstruction error satisfies
\begin{equation}
\label{eq:ei_mls_bound}
\|D^\mu e_i\|_{L^\infty(\omega_i)}
\le C_e h^{p+1-|\mu|},
\qquad |\mu|\le p.
\end{equation}

Substituting $U_i = u + e_i$ into the jump decomposition \eqref{eq:tag21} gives $\llbracket D^\alpha u^h\rrbracket_\Gamma = J_u + J_e$, where $J_u$ collects all $u$-terms and $J_e$ collects all $e_i$-terms. Because $u$ has a single set of interfacial jets, $J_u$ factors as in \eqref{eq:jump_factorization} and vanishes by PoZ. Hence only the reconstruction error contributes:
\begin{equation}
\label{eq:jump_error_only}
\llbracket D^\alpha u^h\rrbracket_\Gamma = J_e(\alpha;\Gamma),
\qquad |\alpha|\le p.
\end{equation}
Each summand in $J_e$ inherits the form $\llbracket D^\beta N_i\rrbracket_\Gamma\,(D^{\alpha-\beta}e_i)|_\Gamma$ from \eqref{eq:tag21} (explicit expression in Appendix~\ref{app:jump_decay_proof}). Combining the inverse estimate with \eqref{eq:ei_mls_bound} on each such term, and summing over the $\mathcal O(1)$ active indices $\mathcal I(\Gamma)$, yields
\begin{equation}
\label{eq:asymp_cp_rate}
\|\llbracket D^\alpha u^h\rrbracket_\Gamma\|_{L^\infty(\Gamma)}
\le
\widehat C\, h^{p+1-|\alpha|},
\qquad |\alpha|\le p,
\end{equation}
with $\widehat C$ independent of $h$. A detailed derivation is given in Appendix~\ref{app:jump_decay_proof}.

\noindent\textbf{Asymptotic smoothness.}
Estimate \eqref{eq:asymp_cp_rate} states that, although MiGFEM is \emph{not} globally $C^p$ at finite $h$, its inter-element derivative jumps decay at the same rate as the patchwise reconstruction error for $u\in C^{p+1}(\Omega)$.

\begin{remark}
Since jumps in \eqref{eq:jump_error_only} are driven entirely by the local errors $e_i$, enrichment that strengthens the local approximant simultaneously suppresses inter-element derivative jumps. Thus mesh-intrinsic enrichment improves both patchwise accuracy and inter-element regularity at the same rate.
\end{remark}

\subsection{Routes to exact interface coherence}
\label{sec:exact_regularity_mechanisms}

The interface coherence required for exact jump cancellation \eqref{eq:exact_jump_cancellation} can be realized at finite mesh size in two distinct ways, summarized in Table~\ref{tab:exact_interface_coherence}.

\paragraph{(1) Via mesh/patch topology.}
The traces $D^\mu U_i|_\Gamma$ are matched because of how patches and local approximants are defined, for \emph{any} nodal data after a standard solve---without requiring $u\in\mathcal{V}_i$. On structured meshes, if the patch reconstruction is taken as a tensor-product Lagrange interpolant (a square system, not the overdetermined CWLS of Appendix~\ref{app:cwls_kronecker}), adjacent patches share the same face values and are function-coherent ($k=0$); \eqref{eq:exact_jump_cancellation} then yields strict $C^1$ regularity \cite{ref49}. If all nodes in a subdomain use one grand nodal patch and a single blended local field, the $U_i$ coincide throughout that subdomain and coherence is strong ($k$ large); this localized exact smoothness is aimed at crack-tip analyses and is reported separately \cite{mixfem}.

\paragraph{(2) Via exact reproduction.}
If $U_i=\varphi|_{\omega_i}$ for a single smooth $\varphi$ on every patch, all interface jets of $U_i$ are traces of $\varphi$, so \eqref{eq:exact_jump_cancellation} holds up to the smoothness of $\varphi$. This is verified in Example~A Cases~1 and~3: the exact field is known, and $\mathcal{V}_i$ is designed so that $U_i\equiv u$ patchwise. For a polynomial manufactured solution, standard polynomial CWLS suffices; for a nonpolynomial field, enrichment functions are added to $\mathcal{V}_i$ to reproduce $u$.

When the exact solution is not contained in $\mathcal{V}_i$, the generic situation is $U_i=u+e_i$ with patch-dependent $e_i\neq 0$, and interface coherence is asymptotic rather than exact (Section~\ref{sec:asymp_jump_taylor}). Example~A Case~2 and the remaining numerical benchmarks in Section~\ref{sec:numerical_verification} use this regime.

\begin{table}[htbp]
\centering
\caption{How interface coherence is enforced in MiGFEM at finite $h$.}
\label{tab:exact_interface_coherence}
\footnotesize
\setlength{\tabcolsep}{3.5pt}
\begin{tabular}{@{}p{0.12\linewidth}p{0.22\linewidth}p{0.30\linewidth}p{0.10\linewidth}p{0.20\linewidth}@{}}
\toprule
\textbf{Route} & \textbf{Setting} & \textbf{Coherence on $\Gamma$} & \textbf{$k$} & \textbf{Numerical use} \\
\midrule
Exact & Tensor-product Lagrange, structured mesh\newline (square system; not generic CWLS) & Shared-face value match & $0$ & strict $C^1$ MiGFEM, separate report \cite{ref49} \\
\midrule
Exact & Shared nodal patch on subdomain & Single $U_i$ on subdomain. & large & Smooth crack tip analysis, separate report \cite{mixfem}\\
\midrule
Exact & Known polynomial $u$; $\mathcal{V}_i$ built to contain $u$ & $U_i\equiv u$ on every patch & smoothness of $u$ & Example~A Case~1 \\
\midrule
Exact & Known nonpolynomial $u$; enriched $\mathcal{V}_i$ so nonpolynomial $u\in\mathcal{V}_i$ & $U_i\equiv u$ on every patch & smoothness of $u$ & Example~A Case~3 \\
\midrule
Asymptotic & General unknown $u$; unstructured mesh & $U_i=u+e_i$ at finite $h$ & --- & Example~A Case~2; The rest of the numerical examples in Sec.~\ref{sec:numerical_verification}  \\
\bottomrule
\end{tabular}
\end{table}

\section{Leading Leibniz Derivative}
\label{sec:leibniz_leading_derivative}

For the MiGFEM approximation $u^h=\sum_i N_i U_i$, the Leibniz expansion on element $E$ reads
\begin{equation}
\label{eq:uh_derivative_product}
D^\alpha u^h
=
\sum_{i\in\mathcal I(E)}\sum_{0\le\beta\le\alpha}
\binom{\alpha}{\beta}
D^\beta N_i\,D^{\alpha-\beta}U_i,
\qquad |\alpha|\le p,
\end{equation}
where $\mathcal I(E)$ is the active PoU index set on $E$, as in Section~\ref{sec:smoothness_mechanism}. Isolating the $\beta=0$ term defines the \emph{leading Leibniz derivative}
\begin{equation}
\label{eq:leading_leibniz_def}
\widetilde D^\alpha u^h
:=
\sum_{i\in\mathcal I(E)}N_i D^\alpha U_i,
\end{equation}
and the \emph{Leibniz remainder}
\begin{equation}
\label{eq:Ralpha_full}
R^\alpha
:=
D^\alpha u^h - \widetilde D^\alpha u^h
=
\sum_{0<\beta\le\alpha}\binom{\alpha}{\beta}
\sum_{i\in\mathcal I(E)}D^\beta N_i\,D^{\alpha-\beta}U_i.
\end{equation}

Substituting $U_i=u+e_i$ into $R^\alpha$ and separating the $u$- and $e_i$-contributions, the $u$-part factors as $\sum_i D^\beta N_i$ on each element and vanishes by the PoZ identity (see Appendix~\ref{app:leibniz_remainder} for the explicit split). Hence only the patch errors contribute:
\begin{equation}
\label{eq:Ralpha_ei_only}
R^\alpha
=
\sum_{0<\beta\le\alpha}\binom{\alpha}{\beta}
\sum_{i\in\mathcal I(E)}D^\beta N_i\,D^{\alpha-\beta}e_i.
\end{equation}

The leading derivative $\widetilde D^\alpha u^h$ is globally continuous: each $D^\alpha U_i$ is smooth on its patch, while the $N_i$-blending inherits the $C^0$ regularity of the PoU. The remainder is driven entirely by $e_i$ and satisfies $\|R^\alpha\|_{L^\infty(E)} = \mathcal O(h^{p+1-|\alpha|})$ (Appendix~\ref{app:leibniz_remainder}), the same asymptotic order as the reconstruction error \eqref{eq:ei_mls_bound}. Therefore, using $\widetilde D^\alpha u^h$ in place of the full derivative does not reduce the formal approximation order. Moreover, when $u$ belongs to the local reconstruction space on every patch ($e_i\equiv 0$), we have $R^\alpha\equiv 0$ and $\widetilde D^\alpha u^h = D^\alpha u^h$ exactly.

In the discretization framework that follows, weak Galerkin retains the full Leibniz derivative \eqref{eq:uh_derivative_product} with integration by parts, while strong-form collocation employs the pointwise-continuous operator $\widetilde D^\alpha u^h$ on the identical nodal unknowns.

\section{Discretizations of Strong and Weak Forms}
\label{sec:unified_discretizations}

This section formulates weak Galerkin (WG) and strong-form collocation on the \emph{same} mesh-intrinsic trial space introduced in Section~\ref{sec:migfem_approximation}. Both discretizations share identical nodal unknowns $\{u_k\}$ and local approximants $\{U_i\}$, differing only in how the governing differential operator is evaluated and enforced. The weak formulation employs the full Leibniz derivative (\ref{eq:uh_derivative_product}) within a standard Galerkin variational framework, whereas the strong formulation replaces the full derivative by the globally continuous leading Leibniz derivative $\widetilde D^\alpha u^h$ (\ref{eq:leading_leibniz_def}) and enforces the residual at discrete points or over local subdomains.

\subsection{Model problems}

\smallskip
\noindent\textit{Second-order: linear elasticity.}
Let $\Omega\subset\mathbb R^d$ ($d=2$) be a bounded domain with boundary $\partial\Omega = \overline{\partial\Omega_D\cup\partial\Omega_N}$, $\partial\Omega_D\cap\partial\Omega_N=\emptyset$. The governing equations of linear elastostatics are
\begin{align}
-\nabla\cdot\boldsymbol\sigma(\mathbf u) &= \mathbf f &&\text{in }\Omega, \label{eq:elasticity_strong}\\[2pt]
\boldsymbol\sigma &= \mathbb C:\boldsymbol\varepsilon(\mathbf u), \qquad
\boldsymbol\varepsilon(\mathbf u)=\tfrac12\bigl(\nabla\mathbf u + \nabla\mathbf u^{\mathsf T}\bigr),\\[2pt]
\mathbf u &= \mathbf g &&\text{on }\partial\Omega_D, \\[2pt]
\boldsymbol\sigma\cdot\mathbf n &= \mathbf t &&\text{on }\partial\Omega_N,
\end{align}
where $\mathbb C$ is the fourth-order elasticity tensor, $\mathbf n$ the outward unit normal, $\mathbf f$ the body force, $\mathbf g$ the prescribed displacement, and $\mathbf t$ the prescribed traction.

\smallskip
\noindent\textit{Fourth-order: clamped biharmonic equation.}
Let $\Omega\subset\mathbb R^2$. The clamped thin-plate problem reads
\begin{equation}
\Delta^2 u = f \quad\text{in }\Omega,\qquad
u = g_1,\quad \partial_n u = g_2 \quad\text{on }\partial\Omega,
\label{eq:biharmonic_strong}
\end{equation}
where $\partial_n$ denotes the outward normal derivative. This problem tests the ability of the mesh-intrinsic approximation to support fourth-order operators on $C^0$ meshes via the leading Leibniz derivative.

\subsection{Weak Galerkin formulation}

In the weak formulation, the spatial derivatives of the MiGFEM approximation are evaluated through the full Leibniz expansion (\ref{eq:uh_derivative_product}). For linear elasticity, substituting $u^h$ into the principle of virtual work and integrating by parts yields the standard variational problem: find $\mathbf u^h\in\mathcal V^h$ (the MiGFEM trial space) such that
\begin{equation}
\label{eq:weak_abstract}
a(\mathbf w,\mathbf u^h) = \ell(\mathbf w),
\qquad \forall\,\mathbf w\in\mathcal V^h,
\end{equation}
with the bilinear and linear forms
\begin{align}
a(\mathbf w,\mathbf u^h) &=
\int_\Omega \boldsymbol\varepsilon(\mathbf w):\mathbb C:\boldsymbol\varepsilon(\mathbf u^h)\,d\Omega
-
\int_{\partial\Omega_D}\bigl[\,\mathbf w\cdot\boldsymbol\sigma(\mathbf u^h)\cdot\mathbf n
+ \boldsymbol\sigma(\mathbf w)\cdot\mathbf n\cdot\mathbf u^h
- \beta\,\mathbf w\cdot\mathbf u^h\,\bigr]\,d\Gamma,\label{eq:weak_a}\\[4pt]
\ell(\mathbf w) &=
\int_\Omega \mathbf w\cdot\mathbf f\,d\Omega
+
\int_{\partial\Omega_N}\mathbf w\cdot\mathbf t\,d\Gamma
-
\int_{\partial\Omega_D}\bigl[\,\boldsymbol\sigma(\mathbf w)\cdot\mathbf n\cdot\mathbf g
- \beta\,\mathbf w\cdot\mathbf g\,\bigr]\,d\Gamma,\label{eq:weak_ell}
\end{align}
where $\beta$ is a Nitsche penalty parameter chosen sufficiently large to ensure coercivity of $a(\cdot,\cdot)$. The boundary integrals on $\partial\Omega_D$ enforce the essential (Dirichlet) condition weakly through the symmetric Nitsche method, which is necessary because MiGFEM trial functions $\phi_k$ (\ref{eq:phi_k_def}), being high-order local approximants, do not vanish along element edges on $\partial\Omega_D$ except at the node itself. Consequently, the standard finite element practice of omitting boundary terms at Dirichlet edges does not apply.

The weak formulation inherits directly from the full Leibniz derivative and the $C^0$ PoU. Its accuracy is bounded both by the patchwise approximation error and by the consistency error introduced through integration by parts on Dirichlet boundaries. The latter is a distinctive feature general to high-order methods: it does not arise in standard FEM, where basis functions vanish identically on Dirichlet edges.

\subsection{Strong-form collocation by leading Leibniz derivative}
\label{sec:strong_form_collocation}

\subsubsection{Fourth-order problem: clamped biharmonic equation}

The discrete biharmonic operator is constructed from the leading fourth-order derivatives:
\begin{equation}
\widetilde\Delta_h^2 u^h := \widetilde\partial_{xxxx} u^h + 2\widetilde\partial_{xxyy} u^h + \widetilde\partial_{yyyy} u^h,
\label{eq:leading_biharmonic_operator}
\end{equation}
where $\widetilde\partial_{x^m y^n}u^h = \sum_i N_i\,\partial_x^m\partial_y^n U_i$. Each local derivative $\partial_x^m\partial_y^n U_i$ is computed from the local approximation on patch $\omega_i$ (Appendix~\ref{app:cwls_kronecker}). The collocation system is assembled as follows.

Let $\mathcal N$ denote the set of all mesh nodes, partitioned into interior nodes $\mathcal N_{\rm int}$ and boundary nodes $\mathcal N_{\partial}$. Dirichlet conditions $u = g_1$ are imposed strongly: for each $\boldsymbol x_b\in\mathcal N_{\partial}$, the nodal value is set to $u^h(\boldsymbol x_b)=g_1(\boldsymbol x_b)$ and eliminated from the unknown vector, leaving $N_{\rm int}=|\mathcal N_{\rm int}|$ degrees of freedom. For each interior node $\boldsymbol x_i\in\mathcal N_{\rm int}$, the biharmonic equation is collocated as
\begin{equation}
\sum_{j\in\mathcal S_i} w^{\Delta^2}_{ij}\,u^h(\boldsymbol x_j)
= f(\boldsymbol x_i) \;-\! \sum_{b\in\mathcal S_i\cap\partial\Omega} w^{\Delta^2}_{ib}\,g_1(\boldsymbol x_b),
\label{eq:biharmonic_interior}
\end{equation}
where $\mathcal S_i$ is the nodal patch of $\boldsymbol x_i$ and $w^{\Delta^2}_{ij}$ are the local approximation coefficients for $\widetilde\Delta_h^2$. The normal-derivative condition $\partial_n u = g_2$ is enforced by appending one collocation row per boundary node:
\begin{equation}
\sum_{j\in\mathcal S_b} w^{\partial_n}_{bj}\,u^h(\boldsymbol x_j)
= g_2(\boldsymbol x_b) \;-\! \sum_{c\in\mathcal S_b\cap\partial\Omega} w^{\partial_n}_{bc}\,g_1(\boldsymbol x_c),
\label{eq:biharmonic_bc}
\end{equation}
with $w^{\partial_n}_{bj}$ the local approximation coefficients for the normal derivative operator $\partial_n$ evaluated in the outward normal direction $\mathbf n_b$ at $\boldsymbol x_b$.

The assembled system comprises $N_{\rm int}+|\mathcal N_{\partial}|$ equations in $N_{\rm int}$ unknowns, forming an overdetermined algebraic system. To accommodate the disparate operator scales---the biharmonic scaling $\widetilde\Delta_h^2 = O(h^{-4})$ versus the normal derivative scaling $\partial_n = O(h^{-1})$---each row is normalized by its Euclidean norm prior to solution. The overdetermined system is solved using a sparse direct solve of the normal equations.

\subsubsection{Second-order problems: linear elasticity}

For second-order elliptic systems, the strong residual of the discrete equilibrium equation is evaluated using the leading first-order derivative:
\begin{equation}
\label{eq:strong_residual_elasticity}
\widetilde{\mathcal L}\,\mathbf u^h(\boldsymbol\xi_m)
\;:=\;
-\nabla_{\!h}\!\cdot\widetilde{\boldsymbol\sigma}(\mathbf u^h(\boldsymbol\xi_m))
=
\mathbf f(\boldsymbol\xi_m),
\qquad \boldsymbol\xi_m\in\Omega,
\end{equation}
where $\widetilde{\boldsymbol\sigma} = \mathbb C:\widetilde{\boldsymbol\varepsilon}(\mathbf u^h)$ and $\widetilde{\boldsymbol\varepsilon}(\mathbf u^h) = \tfrac12(\widetilde\nabla\mathbf u^h + \widetilde\nabla\mathbf u^{h\mathsf T})$, with $\widetilde\nabla\mathbf u^h$ assembled componentwise as $\widetilde D^1 u^h = \sum_i N_i D^1 U_i$. Dirichlet conditions are imposed by nodal elimination as in the biharmonic case. Neumann (traction) conditions are collocated at boundary edge midpoints $\boldsymbol\xi_m\in\partial\Omega_N$:
\begin{equation}
\widetilde{\boldsymbol\sigma}(\mathbf u^h(\boldsymbol\xi_m))\cdot\mathbf n(\boldsymbol\xi_m)
=
\mathbf t(\boldsymbol\xi_m),
\qquad \boldsymbol\xi_m\in\partial\Omega_N.
\label{eq:neumann_collocation}
\end{equation}

Three strong-form residual-enforcement strategies are employed in the numerical examples. They differ exclusively in how the residual $\widetilde{\mathcal L}\mathbf u^h - \mathbf f$ is tested; the underlying trial space, nodal unknowns, and patch reconstructions are identical across all variants.

\subsection{Collocation and residual-vanishing variants}
\label{sec:collocation_variants}
  
\begin{itemize}
  \item \textbf{NC (nodal collocation).} The strong residual (\ref{eq:strong_residual_elasticity}) is enforced pointwise at every interior mesh node $\boldsymbol x_i\in\Omega$. The PDE rows together with the appended Neumann rows (\ref{eq:neumann_collocation}) and, for fourth-order problems, the normal-derivative rows (\ref{eq:biharmonic_bc}) form the discrete system. For second-order problems with only interior and Neumann rows, the system is square; for fourth-order problems with appended boundary rows, it is overdetermined.

  \item \textbf{CC (cell collocation).} The residual is enforced at element-based sampling points---typically the centroid or a set of interior quadrature points of each element $E\in\mathcal T_h$. The resulting system has more rows than unknowns and is solved in a least-squares sense. Denoting by $\mathcal E$ the set of elements and by $\{\boldsymbol\xi_E\}$ the chosen sampling points, the CC system reads
  \[
  \widetilde{\mathcal L}\mathbf u^h(\boldsymbol\xi_E) = \mathbf f(\boldsymbol\xi_E),\qquad \forall\,E\in\mathcal T_h.
  \]

  \item \textbf{SD (subdomain residual).} The residual is required to vanish in a weighted-average sense over node-centered control volumes $\Omega_c$. For each mesh node $\boldsymbol x_c$, define the control volume
  \[
  \Omega_c := \bigcup_{E\ni\boldsymbol x_c} \omega_{E,c},
  \]
  where $\omega_{E,c}$ is the sub-cell of element $E$ associated with node $\boldsymbol x_c$ (e.g.\ the region obtained by connecting the element centroid to the edge midpoints). The SD condition is
  \[
  \int_{\Omega_c}\bigl(\widetilde{\mathcal L}\mathbf u^h - \mathbf f\bigr)\,d\Omega = 0,
  \qquad \forall\,\boldsymbol x_c\in\mathcal N_{\rm int},
  \]
  with volume integrals evaluated by standard element quadrature (Dunavant rules on triangles). SD can be interpreted as a low-order finite-volume-type residual projection.
\end{itemize}

The three variants span a spectrum of testing strategies: NC provides the strongest pointwise enforcement, CC offers flexibility in overdetermined sampling, and SD supplies a localized integral averaging that mitigates the impact of pointwise singularities or rough data.

The numerical experiments in Section~\ref{sec:numerical_verification} assess these formulations on a sequence of benchmarks---inter-element jump diagnostics, polynomial patch tests, manufactured elasticity problems, biharmonic collocation on distorted meshes, and singular crack-tip enrichment---to verify the theoretical predictions of Sections~\ref{sec:smoothness_mechanism}--\ref{sec:leibniz_leading_derivative} and to characterize the accuracy and efficiency of each variant.

\section{Numerical Verification}
\label{sec:numerical_verification}

\subsection{Numerical setup}
\label{sec:num_setup}

The discretizations of Section~\ref{sec:unified_discretizations} are assessed on a sequence of benchmark problems. Unless otherwise noted, the computational domain is discretized by uniform triangular meshes with characteristic mesh size $h$, and local patch reconstructions employ polynomial degree $p\in\{4,6\}$. Volume integrals in the weak Galerkin (WG) and subdomain residual (SD) variants are evaluated using Dunavant quadrature rules on triangles.

Inter-element jump diagnostics are quantified by the metric
\begin{equation}
\label{eq:jump_metric}
J_{m} := \max_{\text{internal edges }\Gamma}\ \max_{\text{points on }\Gamma}\ \bigl|\text{jump of }\partial_{n}^{m}u^{h}\text{ across }\Gamma\bigr|,\quad m=0,1,\ldots,4.
\end{equation}
Elasticity benchmarks report relative errors in displacement, energy, and $H^1$ seminorm; biharmonic benchmarks use the discrete $\ell^2$ and $h^2$ seminorms defined in Section~\ref{sec:example_group_c_biharmonic}.

\subsection{Example group A: Interface jump behavior}
\label{sec:example_group_a_jumps}

\noindent\textbf{Objective.}
Assess the inter-element smoothness of the mesh-intrinsic approximation on $C^0$ meshes. The test verifies three complementary properties: exact cancellation of derivative jumps when the exact field is contained in the local reconstruction space (Cases~1 and~3), and asymptotic decay of derivative jumps when it is not (Case~2). Case~3 in particular supports the claim in Section~\ref{sec:smoothness_mechanism} that patchwise inclusion of the exact solution induces strict interface coherence and exact PoZ cancellation, not only for polynomials but also for nonpolynomial enriched local spaces.

\noindent\textbf{Setup.}
The nodal values are prescribed from a known exact field. Each local approximant $U_i$ is reconstructed by CWLS, and the global field $u^h$ is obtained by PoU blending. The inter-element jump metrics $J_m$ are measured on internal edges according to Eq.~\eqref{eq:jump_metric}. No global PDE solve is performed. 

\textbf{Case 1} considers the polynomial field $u(x,y)=(1+x+y)^5$ with degree-five local reconstruction on uniform meshes over $[-1,1]^2$. 

\textbf{Case 2} uses the smooth nonpolynomial vector field in Eq.~\eqref{eq:tag84} with a \emph{polynomial-only} local space of degrees $p=3,4,5$ under uniform mesh refinement.  

\textbf{Case 3} uses the same vector field as Case~2, but augments the patch trial space with the harmonics that generate each component---$\{\sin x,\cos x,\sin y,\cos y\}$ together with the products $\cos y\sin x$ and $\cos x\sin y$---on top of a linear ($p=1$) polynomial block (patch depth $s=3$). This embeds a priori solution knowledge in the sense of classical PUM/GFEM enrichment while retaining the mesh-intrinsic, extra-DoF-free construction.

\begin{equation}
\label{eq:tag84}
\mathbf{u}(x,y)=\bigl(\cos(y)\sin(x),\ \cos(x)\sin(y)\bigr)^{\mathsf T},\qquad (x,y)\in[-1,1]^2.
\end{equation}

\subsubsection{Case 1: Polynomial field---exact derivative jump cancellation}

\noindent\textbf{Results.}
Table~\ref{tab:case1_strict_polynomial} lists the maximum inter-element jump norms $J_m$ ($m=0,\ldots,4$). Because $u\in\mathcal P_p$ is reproduced patchwise ($U_i\equiv u$), the exact jumps $\llbracket\partial_n^m u\rrbracket_\Gamma$ vanish on every internal edge $\Gamma$. The reported $J_m$ nevertheless reflect two distinct diagnostic mechanisms.

The value jump $J_0$ remains at $\mathcal O(10^{-14})$ on all meshes simply because $u^h$ is globally conforming by definition; it is analytic rather than numerical.

For $m\ge 1$, $J_m$ is the maximum normal-derivative jump of the blended field $u^h$.
When $U_i\equiv u$ patchwise, \eqref{eq:exact_jump_cancellation} implies $\llbracket D^\alpha u^h\rrbracket_\Gamma=0$ in exact arithmetic.
The table confirms this: $J_1\sim\mathcal O(10^{-11})$ is consistent with roundoff, while
$J_m$ for $m\ge 2$ grows with $m$ and with mesh refinement as a result of roundoff amplification through repeated numerical differentiation---this is a floating-point artifact, not the $\mathcal O(h^{p+1-m})$ inter-element inconsistency of Case~2.

\begin{table}[H]
\centering
\caption{Case 1: Jump metrics for the polynomial field $u(x,y)=(1+x+y)^5$ with $p=5$.}
\label{tab:case1_strict_polynomial}
\begin{tabular}{ccccccc}
\toprule
Mesh & $h$ & $J_0$ & $J_1$ & $J_2$ & $J_3$ & $J_4$ \\
\midrule
$11\times 11$ & 0.2000 & $5.68\times 10^{-14}$ & $1.08\times 10^{-11}$ & $1.14\times 10^{-10}$ & $7.54\times 10^{-10}$ & $2.88\times 10^{-9}$ \\
$21\times 21$ & 0.1000 & $8.53\times 10^{-14}$ & $4.76\times 10^{-11}$ & $9.61\times 10^{-10}$ & $4.81\times 10^{-9}$ & $3.88\times 10^{-8}$ \\
$41\times 41$ & 0.0500 & $8.53\times 10^{-14}$ & $9.09\times 10^{-11}$ & $4.25\times 10^{-9}$ & $4.78\times 10^{-8}$ & $1.31\times 10^{-6}$ \\
\bottomrule
\end{tabular}
\end{table}

\subsubsection{Case 2: Smooth non-polynomial field---asymptotic derivative jump decay}

\begin{figure}[htbp]
\centering
\begin{minipage}{0.33\textwidth}
  \centering
  \includegraphics[width=\linewidth]{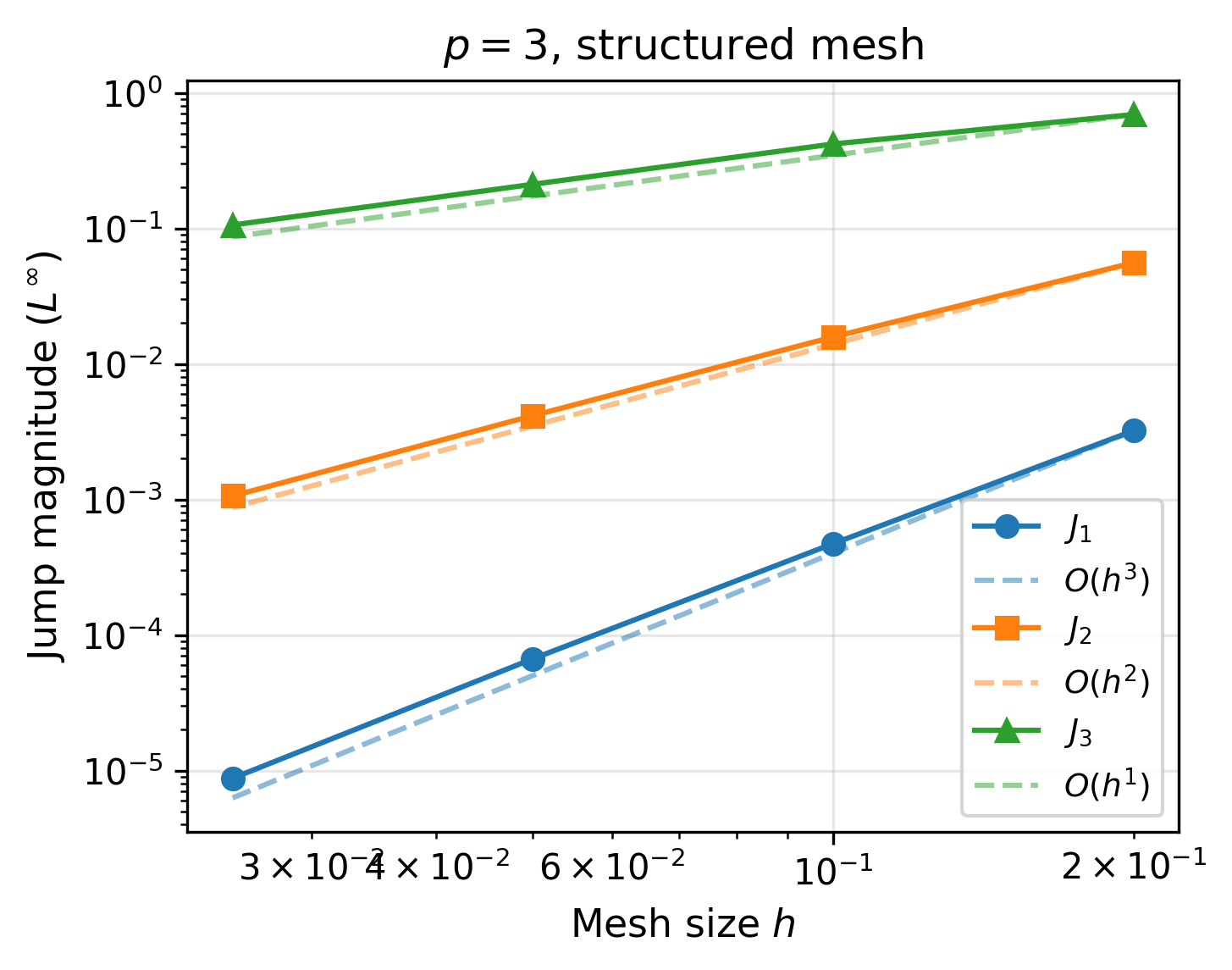}
  \centerline{(a) $p=3$}
\end{minipage}\hfill
\begin{minipage}{0.33\textwidth}
  \centering
  \includegraphics[width=\linewidth]{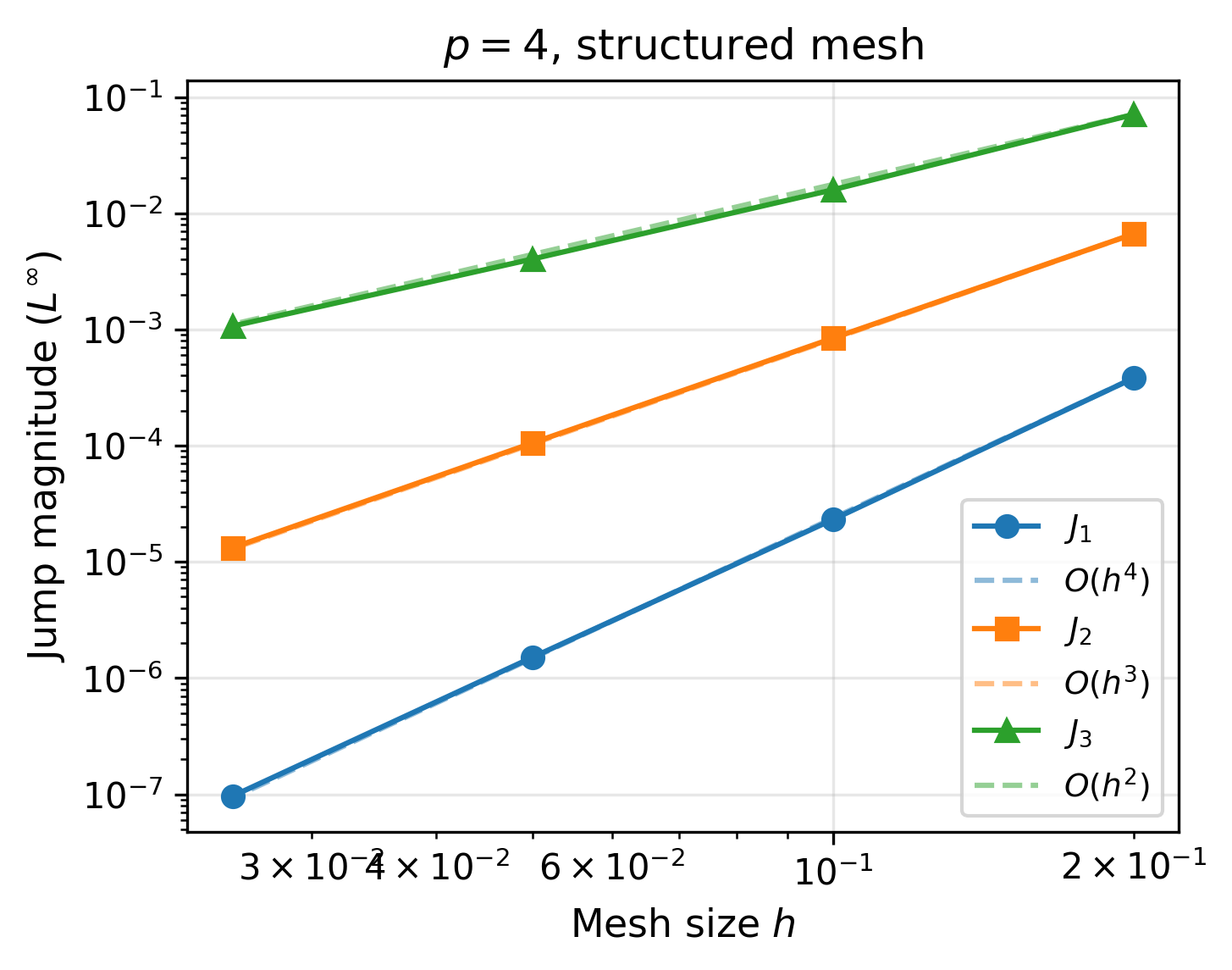}
  \centerline{(b) $p=4$}
\end{minipage}\hfill
\begin{minipage}{0.33\textwidth}
  \centering
  \includegraphics[width=\linewidth]{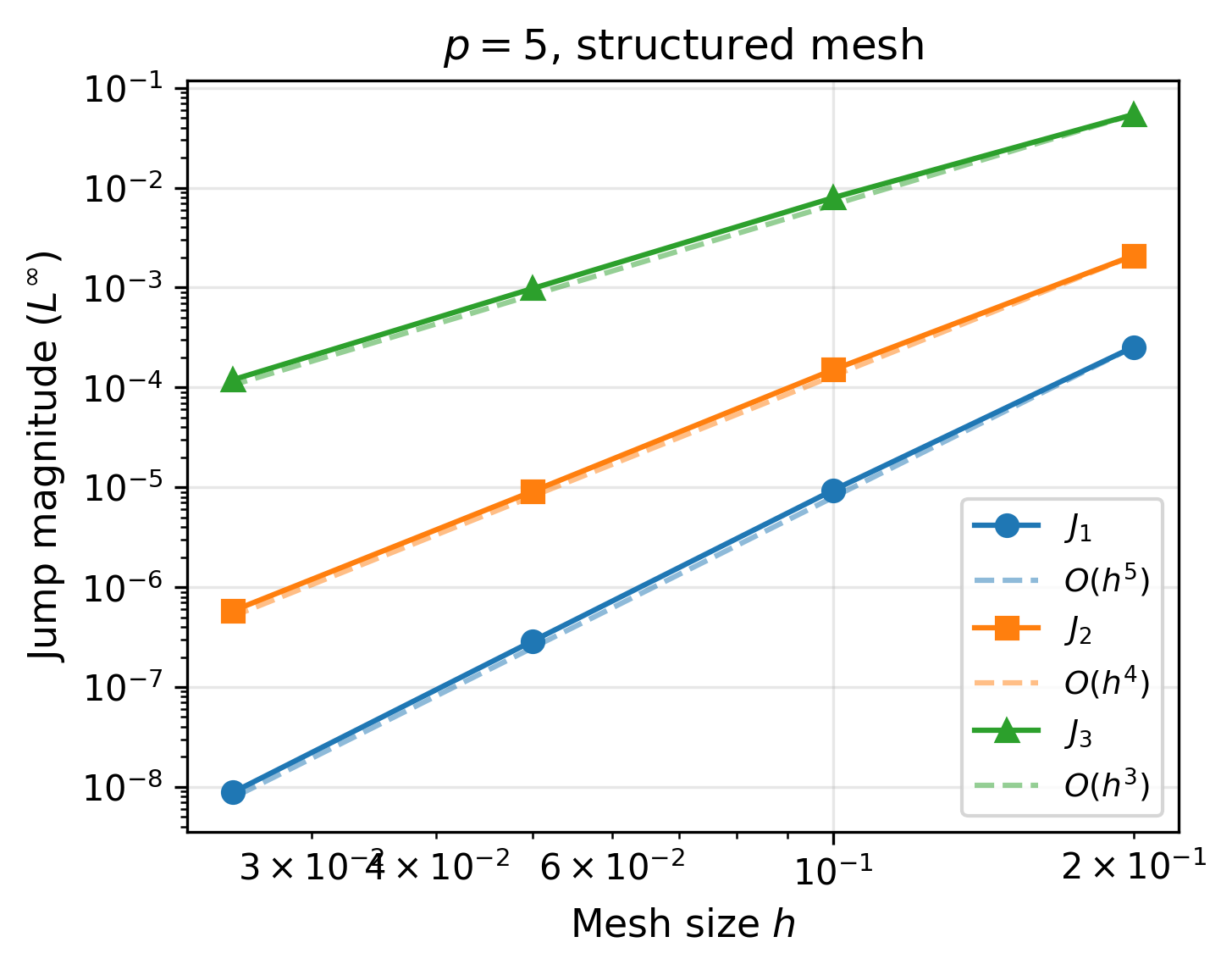}
  \centerline{(c) $p=5$}
\end{minipage}
\caption{Case 2: Inter-element jump metrics $J_m$ on uniform triangular meshes ($p=3,4,5$). Dashed lines indicate the reference slopes $\mathcal{O}(h^{p+1-m})$.}
\label{fig:case2_convergence_all}
\end{figure}

\noindent\textbf{Results.}
Figure~\ref{fig:case2_convergence_all} shows the decay of $J_m$ ($m=1,2,3$) under mesh refinement. The jump norms decrease monotonically with $h$, and the fitted slopes agree with the predicted order $\mathcal{O}(h^{p+1-m})$ from the asymptotic smoothness analysis. In particular, each increase in the derivative order reduces the observed convergence rate by approximately one power of $h$. The results confirm that the inter-element regularity of MiGFEM improves systematically with mesh refinement, even though the global trial space remains $C^0$.

\subsubsection{Case 3: Enriched local space---exact jump cancellation for a nonpolynomial field}
\label{sec:case3_enriched_jump}

\noindent\textbf{Results.}
Table~\ref{tab:case3_enriched_trig} lists $J_m$ for $m=0,\ldots,4$ on three uniform meshes. The value and low-order derivative jumps remain at roundoff level, in contrast to Case~2 where the same field with a polynomial-only space exhibits only $\mathcal{O}(h^{p+1-m})$ decay (Figure~\ref{fig:case2_convergence_all}). Table~\ref{tab:case3_poly_vs_enriched} summarizes the contrast at $h=0.1$: for $J_1$, the polynomial-only reconstruction ($p=4$) is $\mathcal{O}(10^{-5})$, whereas the enriched reconstruction is $\mathcal{O}(10^{-12})$. These results confirm that exact cancellation is a consequence of patchwise inclusion of the exact solution and strict interface coherence, with polynomial reproduction (Case~1) as the special case $\mathcal{P}_p\subset\mathcal{V}_i$.

\begin{table}[H]
\centering
\caption{Case 3: Jump metrics for the trigonometric vector field \eqref{eq:tag84} with enriched local reconstruction (linear monomials + harmonic enrichment, $s=3$).}
\label{tab:case3_enriched_trig}
\begin{tabular}{ccccccc}
\toprule
Mesh & $h$ & $J_0$ & $J_1$ & $J_2$ & $J_3$ & $J_4$ \\
\midrule
$11\times 11$ & 0.2000 & $4.44\times 10^{-16}$ & $2.60\times 10^{-12}$ & $3.63\times 10^{-11}$ & $2.93\times 10^{-10}$ & $8.75\times 10^{-10}$ \\
$21\times 21$ & 0.1000 & $4.44\times 10^{-16}$ & $8.26\times 10^{-12}$ & $2.93\times 10^{-10}$ & $2.65\times 10^{-9}$ & $2.04\times 10^{-8}$ \\
$41\times 41$ & 0.0500 & $4.44\times 10^{-16}$ & $5.42\times 10^{-11}$ & $5.90\times 10^{-10}$ & $3.64\times 10^{-8}$ & $4.73\times 10^{-7}$ \\
\bottomrule
\end{tabular}
\end{table}

\begin{table}[H]
\centering
\caption{Case 2 vs.\ Case 3 at $21\times 21$ ($h=0.1$): maximum inter-element jump $J_1$ for field \eqref{eq:tag84}. Case 2: polynomial-only CWLS ($p=4$); Case 3: linear monomials plus harmonic enrichment.}
\label{tab:case3_poly_vs_enriched}
\begin{tabular}{lcc}
\toprule
Local space & $J_1$ (max over components) & Interface behavior \\
\midrule
Polynomial $\mathcal{P}_4$ only (Case 2) & $2.32\times 10^{-5}$ & $\mathcal{O}(h^{p})$ asymptotic coherence \\
Polynomial $\mathcal{P}_1$ + harmonics (Case 3) & $8.26\times 10^{-12}$ & Exact coherence (machine precision) \\
\bottomrule
\end{tabular}
\end{table}

\subsubsection{Discussion: jump metrics $J_m$ versus the leading Leibniz derivative}
\label{sec:example_a_discussion_jm_leading}

The quantities $J_m$ in \eqref{eq:jump_metric} measure inter-element jumps of the \emph{full} Leibniz derivative $D^\alpha u^h$; Example~A is therefore a numerical check of the asymptotic smoothness statements in Section~\ref{sec:smoothness_mechanism}. We do not tabulate leading Leibniz derivative jumps here because $\widetilde D^\alpha u^h=\sum_i N_i D^\alpha U_i$ is globally continuous and jump-free by construction, independently of the patch reconstruction error. The tables instead validate the jump cancellation mechanism that motivates the use of $\widetilde D^\alpha u^h$ in strong-form collocation.

\subsection{Example group B: Elastic patch test (full solve)---exact derivative jump cancellation}
\label{sec:example_group_b_reproduction}

\noindent\textbf{Objective.}
This patch test is a manufactured-solution benchmark for linear elasticity. The prescribed displacement field is a degree-$p$ polynomial and is therefore exactly contained in the local approximation space $\mathcal P_p$. Under this condition, the PoZ identity together with interface coherence enforces exact cancellation of inter-element derivative jumps for all $|\alpha|\le p$ (Section~\ref{sec:smoothness_mechanism}). The purpose of the test is to verify polynomial reproducibility of the displacement field, the associated strain/stress derivatives, and the ability of the full weak- and strong-form solvers to recover the exact polynomial state on a globally $C^0$ mesh.

\subsubsection{High-order smoothness patch test}

\noindent\textbf{Setup.}
Plane-stress elasticity is considered on $\Omega=(0,5)^2$ with the manufactured polynomial displacement field
\begin{equation}
\label{eq:patch_u_poly}
\mathbf{u}(x,y)=
\begin{pmatrix}
(10+x+y)^p\\[2pt]
(10+x+y)^p
\end{pmatrix},
\end{equation}
from which the body force $\mathbf b$, strain tensor $\varepsilon$, and Cauchy stress $\sigma$ are derived analytically. The local reconstruction uses a $p=4$ polynomial basis with patch depth $s=3$ on an $11\times 11$ mesh. The problem is solved with the strong NC/CC/SD collocation variants, as well as with the weak Galerkin formulation using either Nitsche enforcement or direct elimination.

\noindent\textbf{Results.}
Table~\ref{tab:elasticity_patch_p4s3} lists relative errors after a full global solve. The strong-form variants NC, CC, and SD pass the patch test: displacement and derivative errors are at roundoff ($\mathcal{O}(10^{-10})$--$\mathcal{O}(10^{-8})$). With $\mathbf u\in\mathcal P_p$ patchwise (Section~\ref{sec:exact_regularity_mechanisms}), this shows that MiGFEM reproduces the exact polynomial elasticity state on a globally $C^0$ mesh---including exact interface coherence---and that MiGFEM with the leading Leibniz derivative supports a consistent primal strong discretization even without global $C^p$ conformity.

Weak Galerkin on the same space does not reach roundoff ($\mathcal{O}(10^{-7})$; Nitsche and direct elimination agree). The trial space is the same; the essential distinction is in the weak formulation. Integration by parts yields boundary terms that separate into natural (Neumann) and essential (Dirichlet) contributions. MiGFEM trial functions are high-order reconstructions tied to nodal values: they vanish at nodes but not, in general, along entire Dirichlet edges. The essential boundary integral therefore does not drop out automatically and must be kept on the left-hand side of the variational form together with the volume terms---unlike linear interpolants that are zero on essential boundaries. This structural feature contributes essentially to the gap between weak and strong forms; the boundary treatment error dominates when the bulk solution is already exact. Primal strong collocation does not pass through this integrated-by-parts essential boundary term and hence recovers the patch solution to roundoff when $\mathbf u\in\mathcal P_p$.

\begin{table}[h]
\centering
\caption{Elastic patch tests on a $11\times 11$ mesh ($p=4$, $s=3$). Relative errors in displacement and derivatives; not an interpolation-only test.}
\label{tab:elasticity_patch_p4s3}
\resizebox{\textwidth}{!}{%
\begin{tabular}{cccccc}
\toprule
Method & $\|e_u\|_{\infty}$ & $\|e_{D^1u}\|_{\infty}$ & $\|e_{D^2u}\|_{\infty}$ & $\|e_{D^3u}\|_{\infty}$ & $\|e_{D^4u}\|_{\infty}$ \\
\midrule
NC & $4.07\times10^{-10}$ & $6.22\times10^{-10}$ & $1.20\times10^{-9}$ & $9.60\times10^{-9}$ & $2.24\times10^{-8}$ \\
CC & $1.26\times10^{-9}$ & $8.88\times10^{-10}$ & $3.99\times10^{-9}$ & $2.61\times10^{-8}$ & $6.66\times10^{-8}$ \\
SD & $1.28\times10^{-9}$ & $6.38\times10^{-10}$ & $5.09\times10^{-9}$ & $5.08\times10^{-8}$ & $1.27\times10^{-7}$ \\
WG (Nitsche) & $2.29\times10^{-7}$ & $1.60\times10^{-7}$ & $2.19\times10^{-7}$ & $8.71\times10^{-7}$ & $1.68\times10^{-6}$ \\
WG (direct elimination) & $2.30\times10^{-7}$ & $1.66\times10^{-7}$ & $2.53\times10^{-7}$ & $7.03\times10^{-7}$ & $1.84\times10^{-6}$ \\
\bottomrule
\end{tabular}%
}
\end{table}

\subsubsection{Manufactured elasticity: MiGFEM strong vs.\ weak}
\label{sec:elasticity_migfem_performance}

\noindent\textbf{Objective.}
Assess the same mesh-intrinsic trial space in weak and strong formulations for a smooth manufactured elasticity problem. The test compares weak Galerkin (WG) with primal strong collocation on the same meshes and with the same local reconstructions, to verify whether the enriched approximation can support both variational and pointwise discretizations on unstructured $C^0$ meshes.

\noindent\textbf{Setup.}
Plane stress on $\Omega=(0,5)^2$ with
\begin{equation}
\label{eq:benchmark_u_trig}
\mathbf{u}(x,y) =
\begin{pmatrix}
\cos(y)\sin(x) \\
\cos(x)\sin(y)
\end{pmatrix},
\end{equation}
with the corresponding body force and traction. Uniform triangular meshes from $11\times 11$ to $81\times 81$ are used, with $p=4$ and $p=6$. NC and SD use Dirichlet elimination with collocated Neumann conditions (Section~\ref{sec:unified_discretizations}); WG uses Nitsche enforcement on the essential boundary.

\begin{figure}[H]
\centering
\begin{minipage}{0.95\textwidth}
  \centering
  \includegraphics[width=\textwidth]{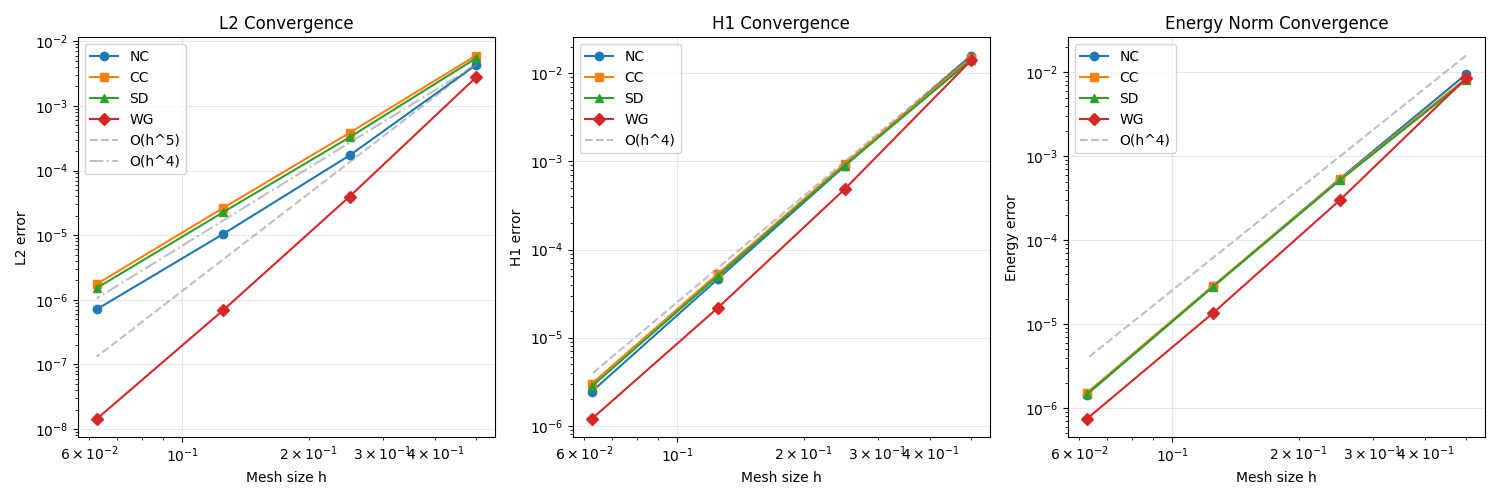}
  \captionline{(a) $p=4$}
\end{minipage}

\vspace{0.2cm}
\begin{minipage}{0.95\textwidth}
  \centering
  \includegraphics[width=\textwidth]{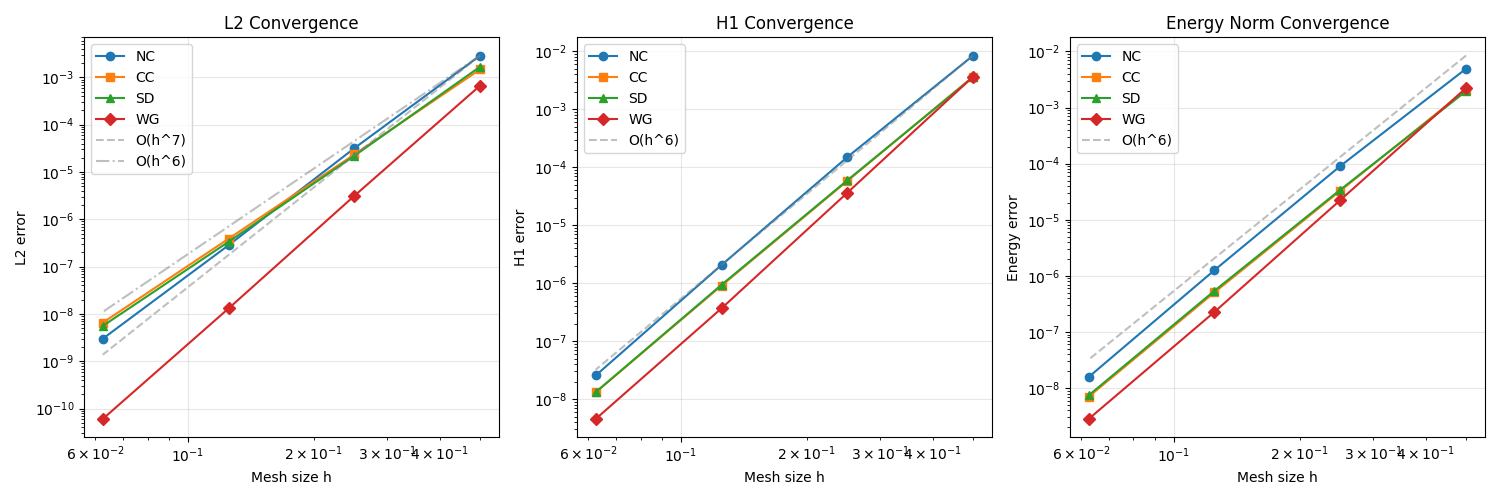}
  \captionline{(b) $p=6$}
\end{minipage}
\caption{MiGFEM NC, CC, SD, and WG on mixed-boundary elasticity (\eqref{eq:benchmark_u_trig}).}
\label{fig:elasticity_convergence}
\end{figure}

\noindent\textbf{Results.}
Figure~\ref{fig:elasticity_convergence} and Table~\ref{tab:elasticity_convergence_orders} summarize the $L^2$ and energy errors. Weak Galerkin, which retains the full Leibniz derivative in the volume form, attains the expected optimal rates in $L^2$ ($\approx p{+}1$) and energy ($\approx p$). NC, CC, and SD assemble the elasticity operator from the \emph{leading} Leibniz derivative on the \emph{same} local reconstructions and nodal unknowns; they converge monotonically in both norms and match WG closely in energy ($\approx p$ for $p=4,6$). That energy behavior is the primary indicator here: equilibrium and the strain-energy bilinear form depend on first-order strains built from $\widetilde D u^h$, and the observed $\mathcal O(h^p)$ energy rates show that replacing $D u^h$ by $\widetilde D u^h$ in the strong collocation residuals yields a consistent, stable strong-form discretization for this smooth manufactured problem, not a breakdown of the approximation space. The displacement $L^2$ rates of NC, CC, and SD are about one order below WG ($\approx p$ rather than $p{+}1$), while the polynomial patch test above recovers the exact state to roundoff when $\mathbf u\in\mathcal P_p$. These results separate the leading-derivative strong operator from the additional error introduced by primal strong-form collocation, which is examined next.

\begin{table}[H]
\centering
\caption{Fitted convergence orders (last three mesh levels) for NC, CC, SD, and WG.}
\label{tab:elasticity_convergence_orders}
\begingroup
\small
\setlength{\tabcolsep}{4pt}
\renewcommand{\arraystretch}{1.1}
\begin{tabular}{lcccc}
\toprule
Method & $L^2$ ($p=4$/$p=6$) & Energy ($p=4$/$p=6$) & Opt.\ $L^2$ ($p{+}1$) & Opt.\ energy ($p$) \\
\midrule
NC     & $4.2$ / $6.6$ & $4.2$ / $6.1$ & $5$ / $7$ & $4$ / $6$ \\
CC     & $3.9$ / $5.9$ & $4.2$ / $6.1$ & $5$ / $7$ & $4$ / $6$ \\
SD     & $4.0$ / $6.2$ & $4.2$ / $5.9$ & $5$ / $7$ & $4$ / $6$ \\
WG     & $5.8$ / $7.9$ & $4.4$ / $6.6$ & $5$ / $7$ & $4$ / $6$ \\
\bottomrule
\end{tabular}%
\endgroup
\end{table}

\subsection{Example group C: Biharmonic strong-form solves}
\label{sec:example_group_c_biharmonic}

\noindent\textbf{Objective.}
This benchmark demonstrates that the same $C^0$ mesh-intrinsic trial space solves the biharmonic equation by collocation, on both regular and strongly distorted meshes, using the leading Leibniz derivative to construct the fourth-order discrete operator.

\begin{figure}[htbp]
\centering
\begin{minipage}{0.45\textwidth}
  \centering
  \includegraphics[width=\linewidth]{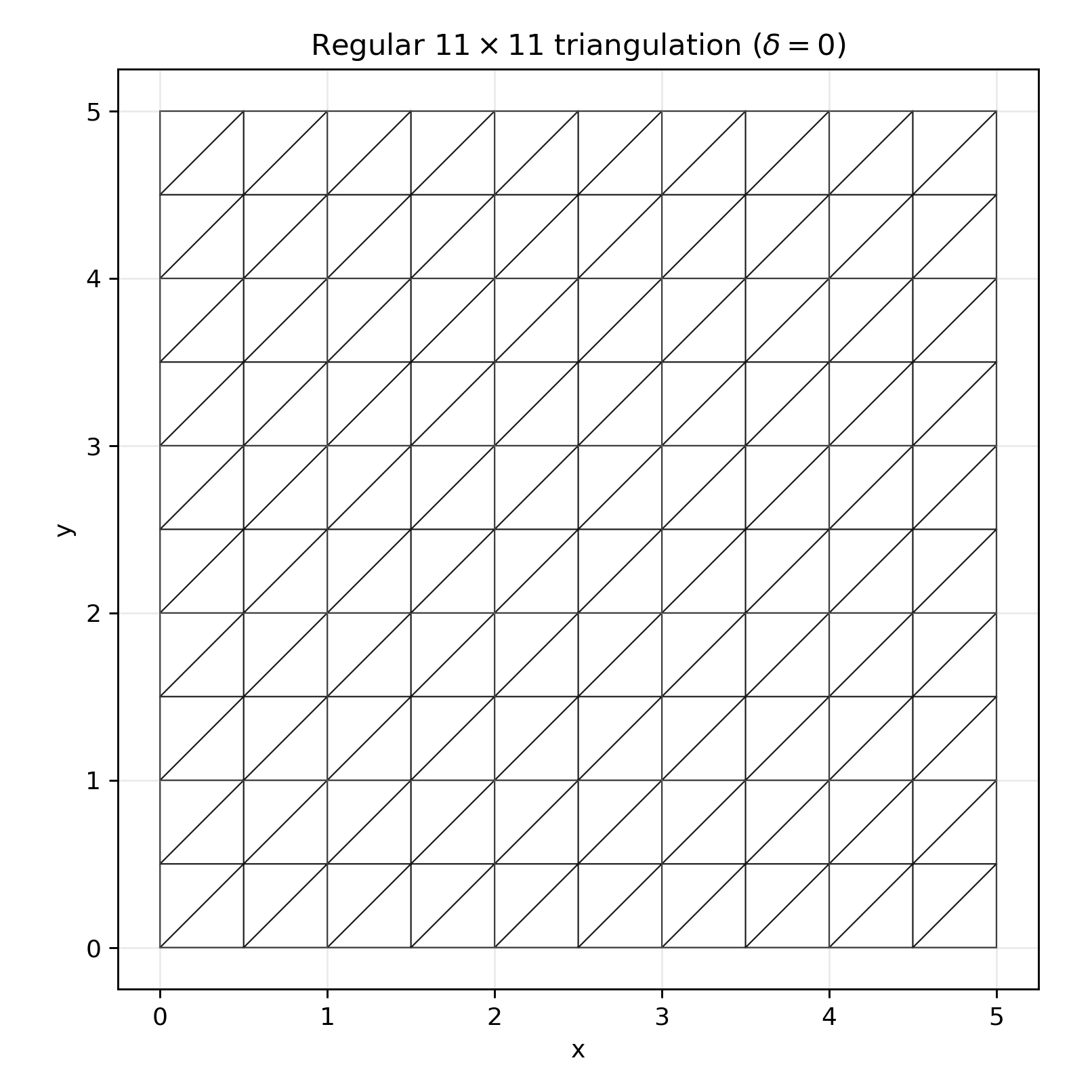}
  \captionline{(a) Regular mesh ($\delta=0.0h$)}
\end{minipage}\hfill
\begin{minipage}{0.45\textwidth}
  \centering
  \includegraphics[width=\linewidth]{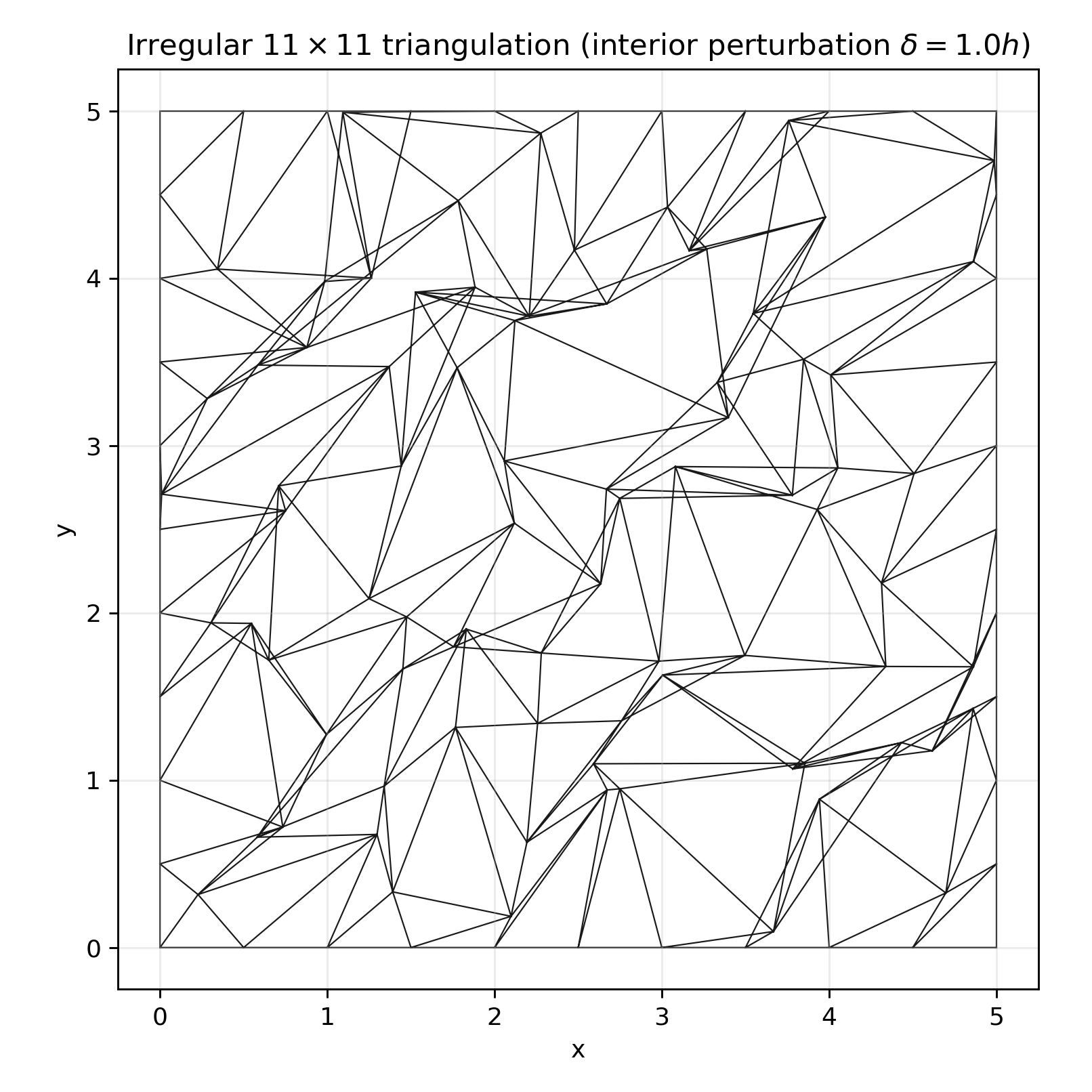}
  \captionline{(b) Strong perturbation ($\delta=1.0h$)}
\end{minipage}
\caption{Representative $11\times 11$ meshes for the biharmonic benchmark: regular ($\delta=0.0h$) and strongly distorted ($\delta=1.0h$).}
\label{fig:biharmonic_irregular_meshes_examples}
\end{figure}

\noindent\textbf{Setup.}
Consider the biharmonic problem on \(\Omega=(0,1)^2\):
\begin{equation}
\label{eq:biharmonic_nv}
\Delta^2 u = f \;\; \text{in }\Omega,\qquad
u = g_1,\;\; \partial_n u = g_2 \;\; \text{on }\partial\Omega,
\end{equation}
with manufactured solution \(u(x,y)=\sin(\pi x)\sin(\pi y)\) and \(f=4\pi^4\sin(\pi x)\sin(\pi y)\). The domain is discretized by uniform refinements from \(21\!\times\!21\) to \(101\!\times\!101\) nodes. Interior nodes are perturbed by \(\delta\in\{0,\,1.0h\}\) (regular and strongly distorted; see Figure~\ref{fig:biharmonic_irregular_meshes_examples}). The NC collocation variant is used: the PDE is enforced at all interior nodes; Dirichlet conditions are imposed strongly by eliminating boundary unknowns; normal-derivative conditions are enforced by appended collocation rows at each boundary node, with row-norm scaling applied to the entire system (see Section~\ref{sec:unified_discretizations} for the detailed assembly).

Errors are reported as discrete nodal norms \(\|e\|_{\ell^2}\) (root-mean-square over all mesh vertices) and \(|e|_{h^2}\) (same for second derivatives reconstructed via the local CWLS approximation). The normal-derivative boundary residual \(\varepsilon_{\partial_n}= \operatorname{rms}_{{\boldsymbol x}_b\in\partial\Omega} |\partial_n u^h({\boldsymbol x}_b)-g_2({\boldsymbol x}_b)|\) is also monitored to assess boundary-condition quality.

\begin{figure}[p]
\centering
\includegraphics[width=0.92\textwidth]{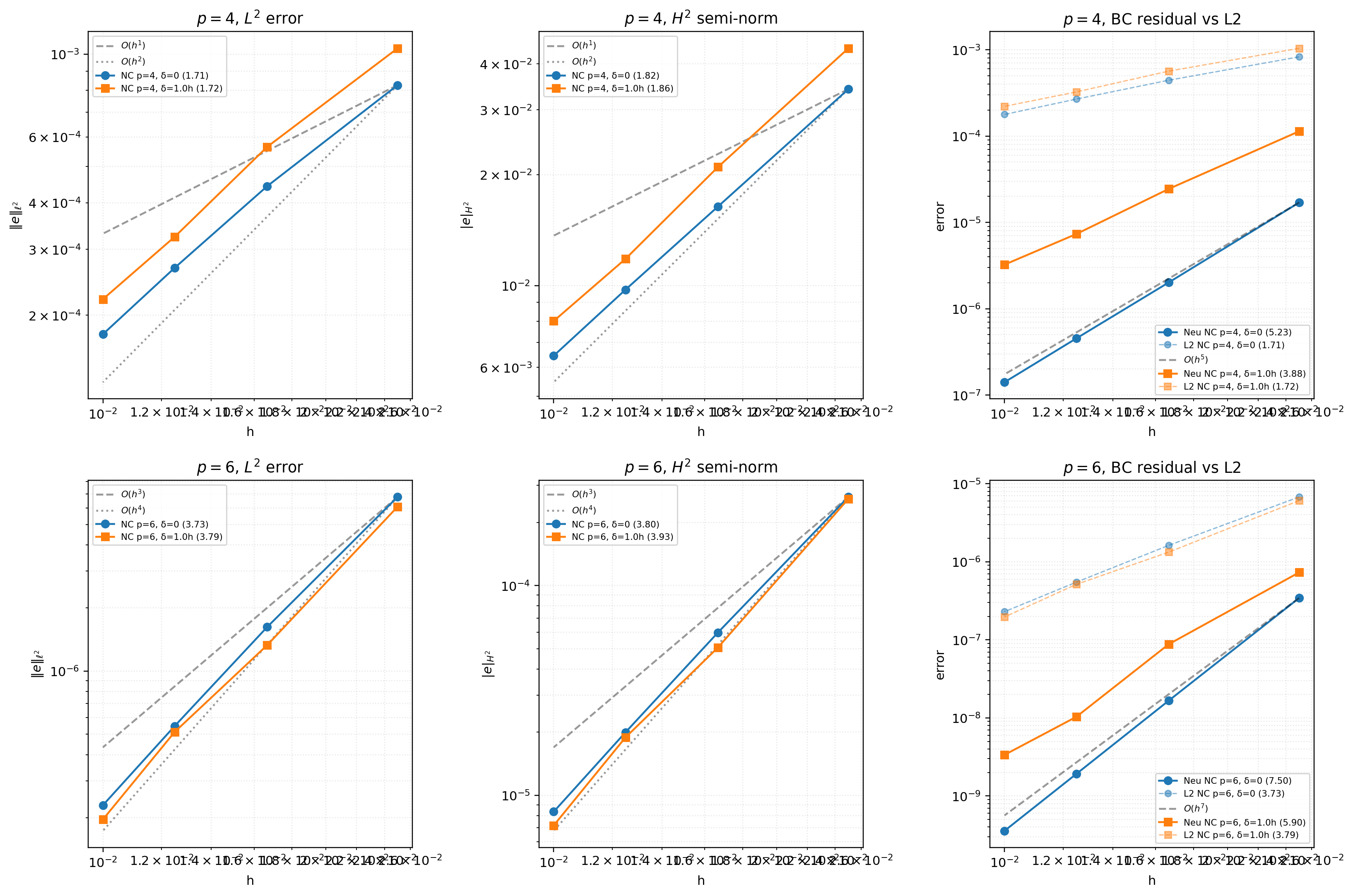}
\caption{Convergence of the biharmonic NC variant for $p=4$ and $p=6$ on regular ($\delta=0$) and distorted ($\delta=1.0h$) meshes. Left: $\|e\|_{\ell^2}$. Center: $|e|_{h^2}$. Right: normal-derivative boundary residual $\varepsilon_{\partial_n}$ versus $\|e\|_{\ell^2}$. Dashed reference lines: $\mathcal{O}(h^{p-3})$, $\mathcal{O}(h^{p-2})$, and $\mathcal{O}(h^{p+1})$ for the $\partial_n$ residual. Legend values are the average slope over the last three refinement levels.}
\label{fig:biharmonic_convergence}
\end{figure}

\noindent\textbf{Results.}
Figure~\ref{fig:biharmonic_convergence} shows the convergence of \(\|e\|_{\ell^2}\), \(|e|_{h^2}\), and the normal-derivative boundary residual \(\varepsilon_{\partial_n}\) for \(p=4\) and \(p=6\) on both regular (\(\delta=0\)) and strongly distorted (\(\delta=1.0h\)) meshes. The asymptotic convergence rates, estimated from the last three mesh levels, are summarized in Table~\ref{tab:biharmonic_slopes}.

\begin{table}[htbp]
\centering
\caption{Pairwise convergence slopes for the biharmonic NC variant. The last three mesh levels (\(h=1/60\to 1/80\to 1/100\)) define the asymptotic regime.}
\label{tab:biharmonic_slopes}
\begin{tabular}{c c c c c c c}
\toprule
\multirow{2}{*}{$p$} & \multirow{2}{*}{$\delta$} & \multicolumn{2}{c}{$\|e\|_{\ell^2}$ slopes} & \multicolumn{2}{c}{$|e|_{h^2}$ slopes} & $\varepsilon_{\partial_n}$ slope \\
\cmidrule(lr){3-4}\cmidrule(lr){5-6}
 & & $h_{60}\!\to\!h_{80}$ & $h_{80}\!\to\!h_{100}$ & $h_{60}\!\to\!h_{80}$ & $h_{80}\!\to\!h_{100}$ & (last three) \\
\midrule
4 & $0.0h$   & 1.75 & 1.83 & 1.81 & 1.85 & $\sim$5.2 \\
4 & $1.0h$   & 1.93 & 1.72 & 2.00 & 1.74 & $\sim$3.9 \\
6 & $0.0h$   & 3.78 & 3.89 & 3.81 & 3.91 & $\sim$7.5 \\
6 & $1.0h$   & 3.31 & 4.30 & 3.43 & 4.33 & $\sim$6.0 \\
\bottomrule
\end{tabular}
\end{table}

On regular meshes, the $\|e\|_{\ell^2}$ and $|e|_{h^2}$ errors decay at approximately $\mathcal{O}(h^{p-2})$: $\sim\!\mathcal{O}(h^{1.8})$ for $p=4$ and $\sim\!\mathcal{O}(h^{3.7})$ for $p=6$. The normal-derivative boundary residual \(\varepsilon_{\partial_n}\) converges substantially faster, at $\mathcal{O}(h^{p+1})$ ($\sim\!\mathcal{O}(h^{5.2})$ and $\sim\!\mathcal{O}(h^{7.5})$), confirming that boundary-condition error does not dominate the total error. Mesh distortion ($\delta=1.0h$) preserves the asymptotic convergence order for both $p=4$ and $p=6$, increasing only the prefactor by a modest factor ($<$25\% in $\|e\|_{\ell^2}$ at $101\!\times\!101$). The $\mathcal{O}(h^{p-2})$ scaling aligns with standard primal isogeometric biharmonic collocation theory (Appendix~\ref{app:biharmonic_order_gap}) and is attributed to the collocation mechanism rather than the trial space.

For reference, nodal injection of the manufactured solution into the MiGFEM trial space follows $\mathcal{O}(h^{p+1})$ in $\|e\|_{\ell^2}$ and $\mathcal{O}(h^{p-1})$ in $|e|_{h^2}$ (Appendix~\ref{app:biharmonic_order_gap}), underlining that the $\mathcal{O}(h^{p-2})$ full-solve rate is governed by the Poisson-like structure of the collocation system, not by approximation capability.

\begin{figure}[htbp]
\centering
\includegraphics[width=0.55\textwidth]{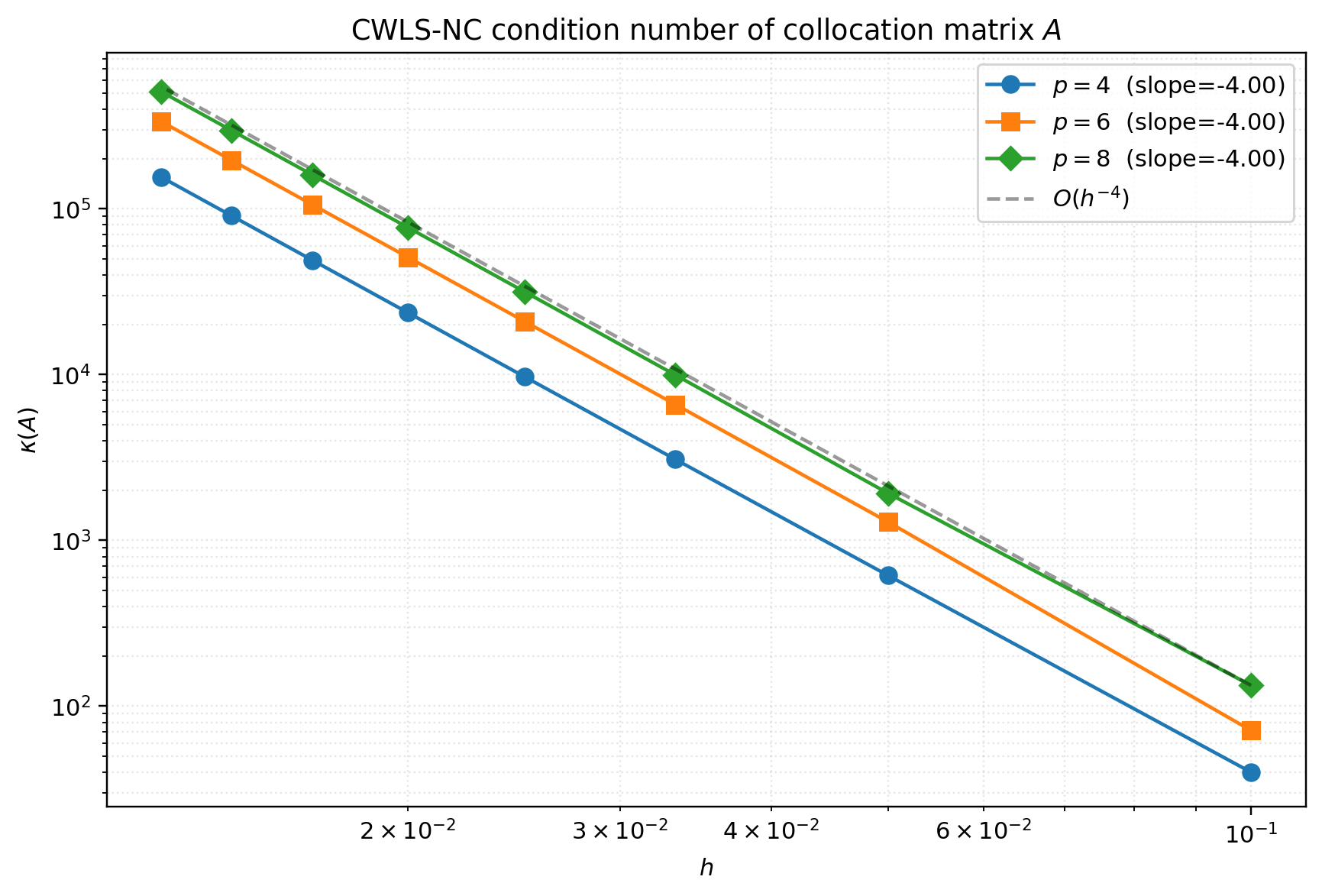}
\caption{Condition number of the row-normalized biharmonic collocation matrix $\mathbf A$ as a function of mesh size $h$. All three polynomial orders ($p=4,6,8$) exhibit $\kappa(\mathbf A)\propto h^{-4.0}$ asymptotically. Legend labels report the average log--log slope over the last three mesh levels.}
\label{fig:biharmonic_condition}
\end{figure}

\noindent\textbf{Stability.}
The stability of the collocation system is assessed by monitoring the condition number $\kappa(\mathbf A)$ of the row-normalized matrix for $p=4$, $p=6$, and $p=8$ on meshes from $11\times11$ to $81\times81$. As shown in Figure~\ref{fig:biharmonic_condition}, $\kappa(\mathbf A)$ grows at the asymptotic rate $\mathcal O(h^{-4.0})$ for all three polynomial orders, consistent with the fourth-order nature of the biharmonic operator. The condition number increases with polynomial order by moderate prefactor ratios: $p=6$ exceeds $p=4$ by a factor of approximately $2.2$, and $p=8$ exceeds $p=6$ by a factor of approximately $1.5$, across all mesh sizes---a modest increase that does not alter the asymptotic scaling. 

The observed $\mathcal O(h^{-4})$ is the irreducible operator-intrinsic scaling; the MiGFEM collocation formulation does not introduce additional sources of ill-conditioning. Together with the preserved convergence rates on strongly distorted meshes ($\delta=1.0h$) reported above, 
these results demonstrate that the MiGFEM strong-form collocation system---with the leading Leibniz derivative---remains stable under both polynomial order increase and severe mesh distortion on $C^0$ unstructured meshes for the tested range of polynomial degrees and mesh distortions.

\subsection{Example group D: Crack-tip singular patch test---strong vs. weak}
\label{sec:rw_tip}

\noindent\textbf{Objective.}
This benchmark is designed to highlight the main practical advantage of extending GFEM local enrichment from a weak-form formulation to a strong-form one: the same singular trial space can deliver both high accuracy and high efficiency, without paying the quadrature cost that usually dominates enriched weak-form solvers. In other words, once the local enrichment is moved to a collocation setting, the crack-tip basis is no longer tied to expensive element integration.

\noindent\textbf{Setup.}
The domain is the unit square $\Omega=[0,1]^2$ with a crack segment ending at $\mathbf x_{\mathrm{tip}}=(0.5,0.5)$. Plane strain isotropic elasticity is used with $E=1.0$ and $\nu=0.3$, and the manufactured solution is the classical Mode-I Westergaard solution with $K_{\mathrm{I}}=1.0$. A structured $11\times11$ triangular mesh is used. Both NC and WG employ the same degree-4 polynomial block plus four crack-tip branch functions, so the approximation space is identical in the two methods. The only difference is the discretization: NC enforces equilibrium by strong collocation, whereas WG uses weak-form Galerkin assembly with element quadrature. The crack tip is excluded from collocation.

\begin{figure}[htbp]
\centering
\includegraphics[width=0.95\linewidth]{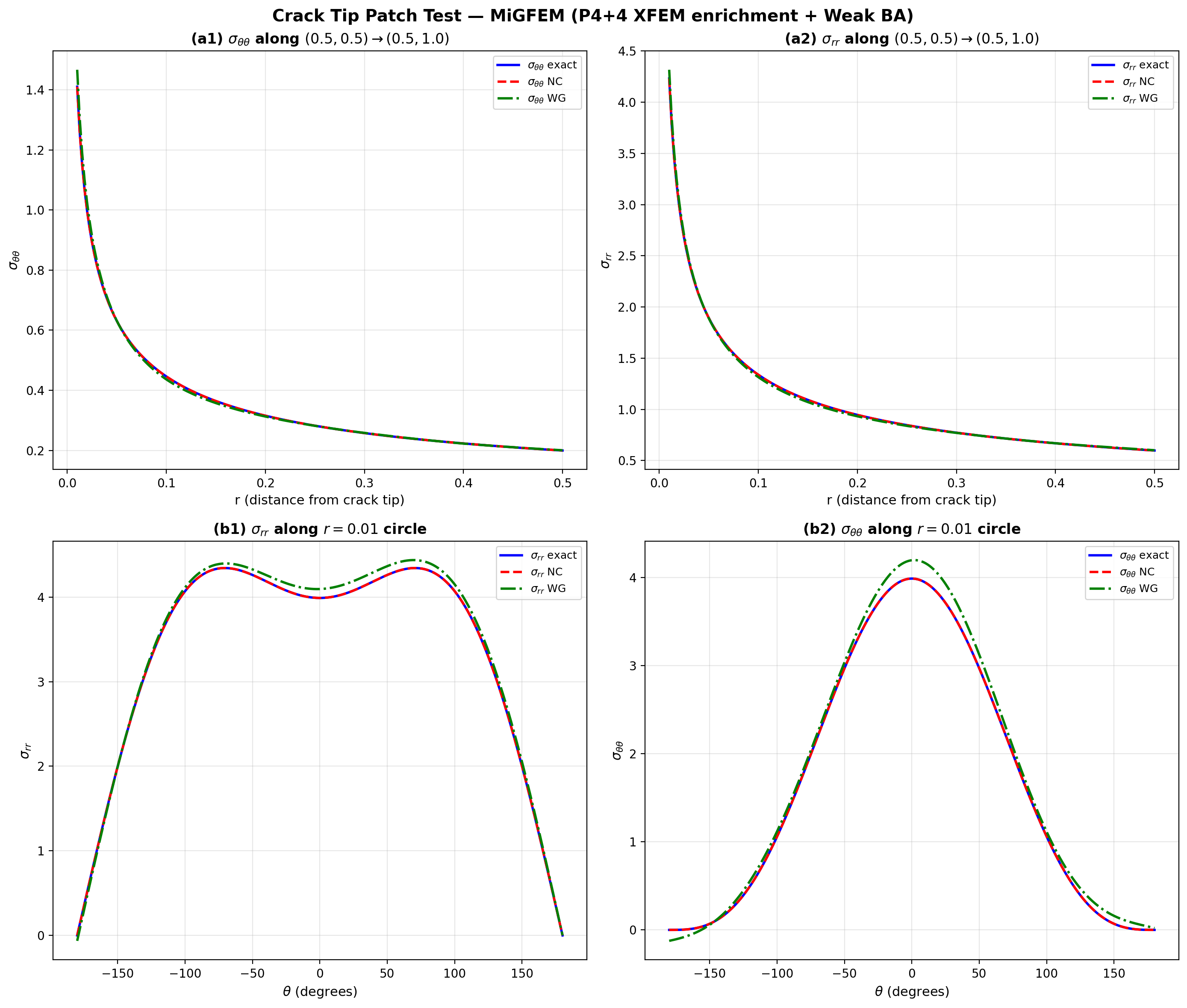}
\caption{Crack-tip patch test results for NC and WG.}
\label{fig:rw_tip_results}
\end{figure}

\noindent\textbf{Results.}
Figure~\ref{fig:rw_tip_results} compares the recovered stress profiles. The NC solution matches the exact singular field, while the WG curves show visible deviations near the tip. Table~\ref{tab:rw_tip_errors} shows that NC reaches machine precision, with $L^2$ displacement error $1.10\times10^{-13}$ and stress RMS errors below $10^{-12}$. In contrast, WG gives $1.16\times10^{-3}$ in $L^2$ and stress RMS errors of order $10^{-2}$. The computational impact is even more striking: Table~\ref{tab:rw_tip_efficiency} shows that the NC total cost is 0.45~s, whereas WG requires 274~s. The quadrature study in Table~\ref{tab:rw_tip_gauss} further confirms that increasing the Gauss order only slowly reduces the WG residual, while the assembly cost rises sharply.

\begin{figure}[htbp]
\centering
\includegraphics[width=0.8\linewidth]{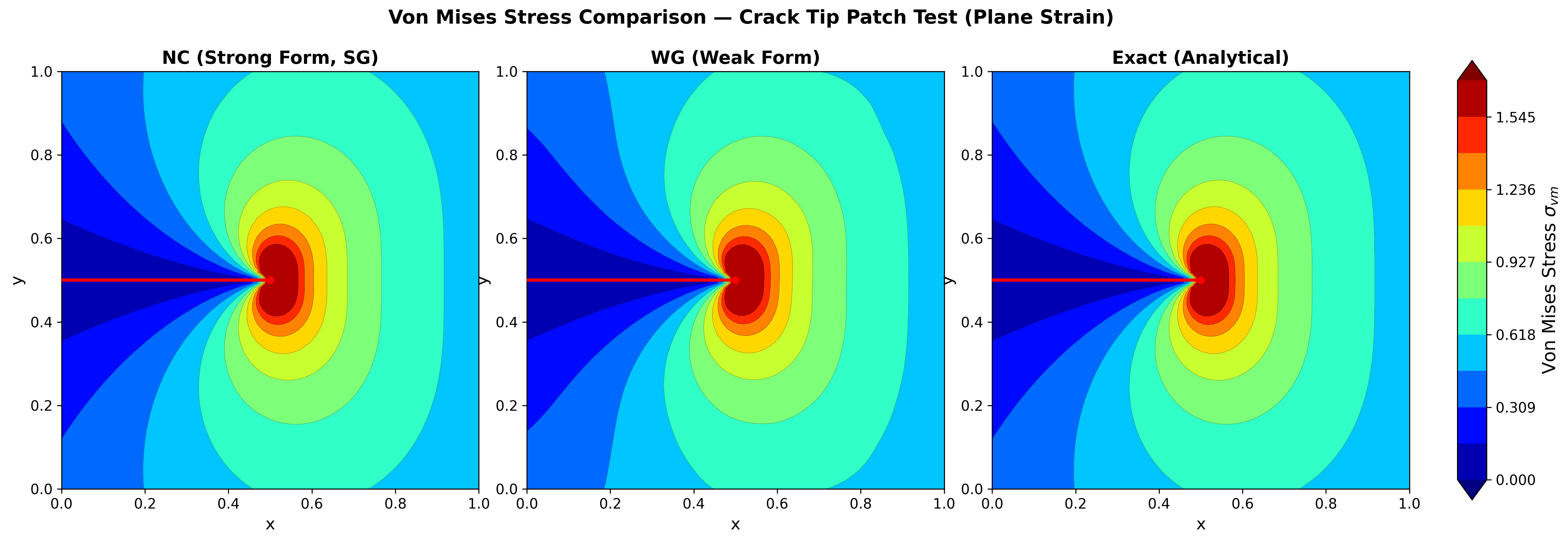}
\caption{Von Mises stress contour for the crack-tip patch test.}
\label{fig:rw_tip_vm}
\end{figure}

\begin{table}[htbp]
\centering
\small
\caption{Quantitative error summary for the crack-tip patch test.}
\label{tab:rw_tip_errors}
\begin{tabular}{@{} l c c @{}}
\toprule
\textbf{Error metric} & \textbf{NC} & \textbf{WG} \\
\midrule
$L^2$ displacement (interior nodes) & $1.10\times10^{-13}$ & $1.16\times10^{-3}$ \\
$\sigma_{\theta\theta}$ RMS along path & $7.01\times10^{-13}$ & $2.03\times10^{-2}$ \\
$\sigma_{rr}$ RMS along path & $3.36\times10^{-13}$ & $1.28\times10^{-2}$ \\
$\sigma_{rr}$ RMS on $r{=}0.01$ circle & $1.78\times10^{-13}$ & $2.02\times10^{-2}$ \\
$\max|u_x^{\mathrm{num}}-u_x^{\mathrm{ex}}|$ (interior) & $2.74\times10^{-13}$ & $3.32\times10^{-3}$ \\
\midrule
Global system residual & $2.83\times10^{-13}$ & $2.94\times10^{-1}$ \\
\bottomrule
\end{tabular}
\end{table}

\begin{table}[htbp]
\centering
\small
\caption{Computational cost for the crack-tip benchmark.}
\label{tab:rw_tip_efficiency}
\begin{tabular}{@{} l r r r @{}}
\toprule
\textbf{Operation} & \textbf{NC} & \textbf{WG} & \textbf{Ratio (WG/NC)} \\
\midrule
LS precompute (shared) & \multicolumn{2}{c}{0.37~s} & --- \\
BA precompute (NC only) & 0.39~s & --- & --- \\
Assembly & 0.004~s & 260~s & $\approx6\times10^4$ \\
Solve (sparse direct) & 0.06~s & 14~s & $\approx230$ \\
\midrule
Total (NC: BA+assemble+solve; WG: assemble+solve) & 0.45~s & 274~s & $\approx600$ \\
\midrule
Total DOFs & 252 & 252 & 1 \\
Unknowns (after elimination) & 168 & 168 & 1 \\
Quadrature points per DOF row & $0$ & $8\times200/168\approx 9.5$ & --- \\
System matrix nonzeros & $\approx 1.2\times10^4$ & $\approx 3.5\times10^5$ & $\approx30$ \\
\midrule
$L^2$ displacement error & $1.1\times10^{-13}$ & $1.2\times10^{-3}$ & --- \\
\bottomrule
\end{tabular}
\end{table}

\begin{table}[htbp]
\centering
\small
\caption{WG interior residual versus Gauss--Legendre quadrature order.}
\label{tab:rw_tip_gauss}
\begin{tabular}{@{} c c c @{}}
\toprule
\textbf{Gauss order} & \textbf{$\max|\mathbf K\mathbf u_{\mathrm{ex}}-\mathbf r|$ (interior)} & \textbf{Assembly time (s)} \\
\midrule
6  & $3.67\times10^{-2}$ & 21.5 \\
8  & $1.15\times10^{-2}$ & 139 \\
10 & $8.26\times10^{-3}$ & 693 \\
\bottomrule
\end{tabular}
\end{table}

\noindent\textbf{Discussion.}
Both methods use the same enriched trial space, so the Westergaard solution is represented exactly at the approximation level. The results underscore that weak-form assembly in singular enrichment problems is limited mainly by quadrature, whereas strong collocation avoids that cost and accuracy loss entirely. Extensive tests for holes, interfaces, and crack growth---where mesh-intrinsic enrichment can operate without mesh cutting and quadrature---are left for future work.

\section{Conclusions}
\label{sec:conclusions}

This paper revisits the mesh-intrinsic GFEM (MiGFEM) of \cite{ref47}---a formulation that constructs enriched local approximations on each nodal patch from the original nodal values alone and blends them through the standard $C^0$ partition of unity (PoU), introducing no extra global unknowns. By construction it avoids the linear dependence and additional-DoF complications of classical GFEM/XFEM.

The paper identifies a property not previously emphasized: \emph{asymptotic smoothness on standard $C^0$ unstructured meshes}. The mechanism rests on two ingredients. The first is the \emph{Partition of Zero} (PoZ): differentiating the PoU identity $\sum_i N_i\equiv 1$ gives $\sum_i D^\alpha N_i\equiv 0$ for $|\alpha|\ge 1$ elementwise. The second is \emph{interface coherence}: since all local reconstructions $U_i$ are built from the same shared nodal unknowns, neighboring patches are tied to a common data set. When the exact solution $u$ is reproduced patchwise ($U_i\equiv u$), the $U_i$ share the same interfacial jets and are derivative-coherent across $\Gamma$.

Through the Leibniz expansion of the MiGFEM derivative jump, PoZ factors out the common interfacial traces and the remaining sum of PoU derivative jumps vanishes. The exact cancellation is
\[
\llbracket D^\alpha u^h\rrbracket_\Gamma = 0,\qquad |\alpha|\le k+1,
\]
whenever the local approximants are $k$-order derivative-coherent on $\Gamma$; polynomial reproduction is the benchmark case. For a general $u\in C^{p+1}(\Omega)$ not contained in the local reconstruction space, writing $U_i=u+e_i$ isolates the exact field's contribution, which vanishes by PoZ. The residual jump is driven solely by the local errors $e_i$, yielding
\[
\|\llbracket D^\alpha u^h\rrbracket_\Gamma\|_{L^\infty(\Gamma)} \le \widehat C\,h^{p+1-|\alpha|},\qquad |\alpha|\le p.
\]
Local enrichment thus plays a dual role: it improves patchwise accuracy and simultaneously suppresses inter-element derivative jumps.

This smoothness mechanism leads naturally to a decomposition of the full Leibniz derivative,
\[
D^\alpha u^h = \underbrace{\sum_{i} N_i D^\alpha U_i}_{\widetilde D^\alpha u^h} \;+\; R^\alpha,\qquad \|R^\alpha\|_{L^\infty(E)} = \mathcal O(h^{p+1-|\alpha|}),
\]
where the \emph{leading Leibniz derivative} $\widetilde D^\alpha u^h$ is globally continuous and jump-free by construction. It serves as the derivative operator for primal strong-form collocation.

The framework unifies weak and strong discretizations on the same mesh-intrinsic trial space. The weak formulation retains the full Leibniz derivative with integration by parts; the strong formulation assembles the differential operator from $\widetilde D^\alpha u^h$ and enforces the residual by nodal collocation (NC), cell collocation (CC), or subdomain-averaged residuals (SD). All variants share identical nodal unknowns and local reconstructions.

Numerical results support the analysis:
\begin{itemize}
	\item Jump diagnostics confirm exact cancellation to roundoff under patchwise reproduction and $\mathcal O(h^{p+1-m})$ decay for general smooth fields, including nonpolynomial enriched spaces.
	\item Elasticity patch tests recover exact polynomial displacement, strain, and stress fields to roundoff with all three strong variants, while weak Galerkin incurs a small boundary-treatment error ($\sim10^{-7}$) inherent to the Nitsche formulation on the enriched $C^0$ space.
	\item Manufactured elasticity yields energy rates of $\mathcal O(h^p)$ for NC and SD, matching weak Galerkin and confirming that the leading Leibniz derivative provides a consistent strong-form operator for second-order elliptic systems.
	\item Biharmonic \(\Delta^2\) solves using NC converge stably on regular and strongly distorted meshes at approximately \(\mathcal O(h^{p-2})\), with the normal-derivative boundary residual converging at the faster \(\mathcal O(h^{p+1})\) rate, confirming effective boundary-condition enforcement. The order deficit relative to nodal injection is attributed to the collocation mechanism---not to the trial space itself (Appendix~\ref{app:biharmonic_order_gap}).
\end{itemize}

In summary, mesh-intrinsic enrichment enriches the local approximation space and induces asymptotic inter-element smoothness simultaneously, without raising global continuity above $C^0$. The leading Leibniz derivative extends PoU-based local enrichment from a predominantly weak-form tool to a unified weak--strong framework on the original nodal layout. This is particularly relevant for material interfaces, voids, and propagating defects, where classical enriched weak formulations incur heavy quadrature and cut-cell costs \cite{refHansboHansbo2002,mixfem}. Asymptotic smoothness explains why enriched $C^0$ approximations can carry stronger inter-element regularity than conventional piecewise polynomials, and why strong collocation becomes practical on $C^0$ unstructured meshes.

% End of main text

\section{Appendix}
\appendix

\section{CWLS local approximation and Kronecker property}
\label{app:cwls_kronecker}
\label{sec:impl_cwls}

This appendix presents the constrained weighted least-squares (CWLS) reconstruction used in the present MiGFEM implementation on unstructured meshes, and shows how the nodal Kronecker property is preserved at the global level. The purpose of the CWLS construction is twofold: it provides a stable patchwise local approximation on general nodal stencils, and it enforces interpolation at the patch-star node so that the final MiGFEM basis remains nodally interpolatory.

\subsection{Patch definition}
\label{sec:patch_stencil}

For each mesh node $i\in\mathcal N$, called the \emph{patch-star node}, we define a nodal patch $\omega_i$ by collecting the elements contained in a prescribed number of element layers around $x_i$. Let
\begin{equation}
\mathcal S_i \subset \mathcal N
\end{equation}
denote the set of nodes contained in this patch. The patch size is chosen so that the stencil is sufficiently large for the intended local approximation space; in particular, one requires
\begin{equation}
|\mathcal S_i| \ge \dim \mathcal P_p
\end{equation}
when a polynomial reconstruction of degree $p$ is used.

\begin{figure}[htbp]
  \centering
  \includegraphics[width=0.8\linewidth]{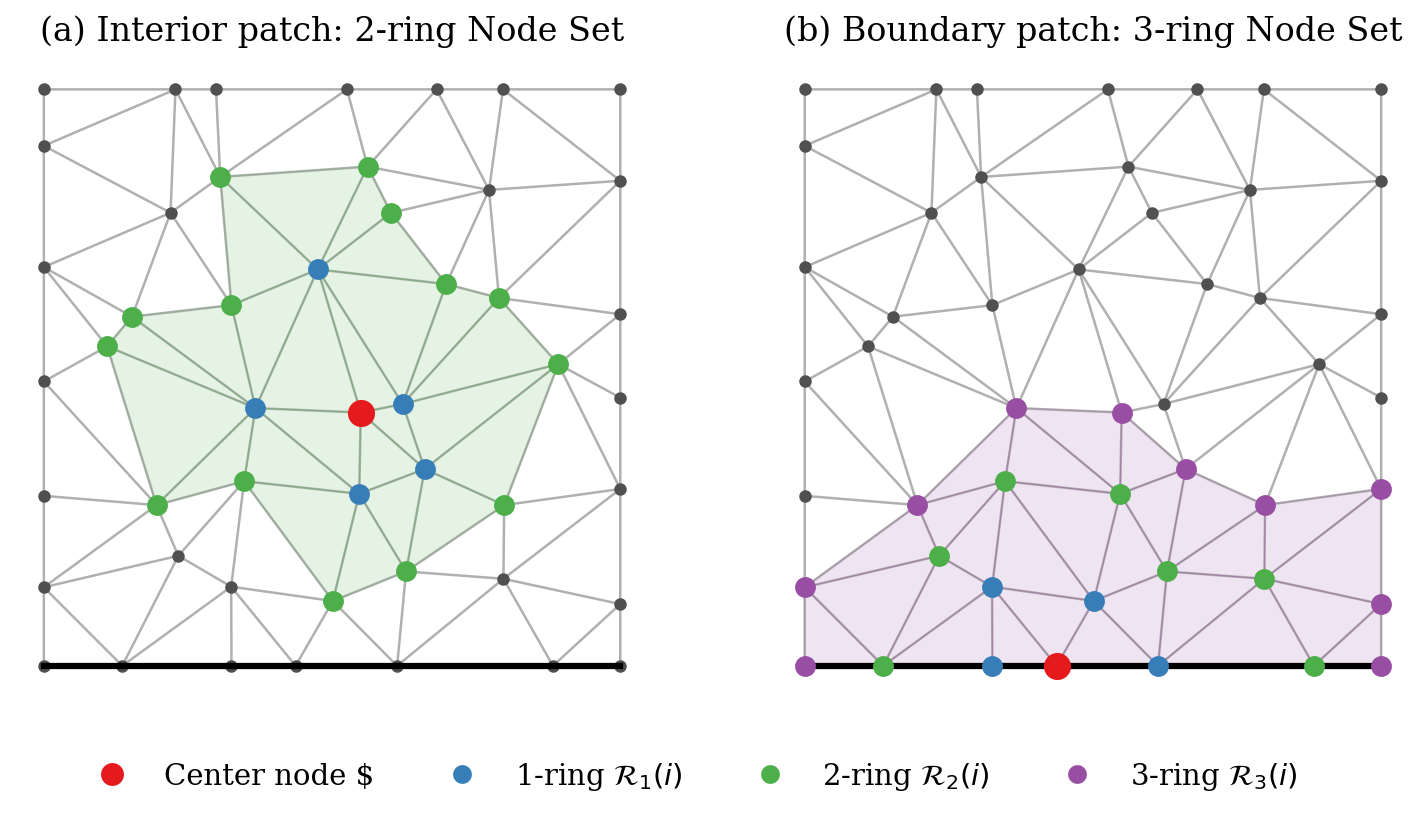}
  \caption{Patch node sets on an unstructured triangular mesh. (a) Interior patch, $s=2$. (b) Boundary patch, $s=3$.}
  \label{fig:boundary_stencil}
\end{figure}

In the triangular meshes considered in this work, the number of patch layers is selected empirically to ensure a robust and unisolvent reconstruction. For interior nodes, a typical choice is
\begin{equation}
s = p/2 + 1,
\end{equation}
while for boundary nodes the patch is extended inward so that the local stencil remains sufficiently rich. Examples of interior and boundary patch node sets are shown in Fig.~\ref{fig:boundary_stencil}.

The essential requirement is that each patch provide enough sampling points and a sufficiently regular geometry for the local reconstruction to be well posed. Once the patch is fixed, all local approximation coefficients are determined entirely from the nodal values on $\mathcal S_i$.

\subsection{CWLS local approximation}
\label{sec:cwls_formulation}

On each patch $\omega_i$, the local approximant is represented in a polynomial basis centered at the patch-star node. Let
\begin{equation}
\boldsymbol{\xi}=\frac{\boldsymbol{x}-\boldsymbol{x}_i}{\rho_i},
\end{equation}
where $\rho_i$ is a characteristic patch length used for scaling, and let $\boldsymbol{p}(\boldsymbol{\xi})$ be a basis of the polynomial space $\mathcal P_p$. The local approximant is written as
\begin{equation}
\label{eq:local_poly_param}
U_i(\boldsymbol{x})=\boldsymbol{p}(\boldsymbol{\xi})^\top \boldsymbol{a}_i,
\end{equation}
where $\boldsymbol{a}_i$ is the vector of unknown local coefficients.

Let $\boldsymbol{u}_{\mathcal S_i}$ collect the nodal values on the patch, and let $\mathbf P_i$ be the sampling matrix obtained by evaluating the polynomial basis at the patch nodes. Let $\mathbf W_i$ be a diagonal weight matrix used to emphasize samples closer to the patch-star node. Then the unconstrained weighted least-squares fit would minimize
\[
\frac12 \|\mathbf W_i(\mathbf P_i\boldsymbol{a}-\boldsymbol{u}_{\mathcal S_i})\|_2^2.
\]
However, in order to preserve nodal interpolation in the global MiGFEM approximation, the local approximant is additionally constrained to pass through the patch-star node:
\begin{equation}
U_i(\boldsymbol{x}_i)=u_i.
\end{equation}
Since $\boldsymbol{\xi}=\boldsymbol{0}$ at $\boldsymbol{x}=\boldsymbol{x}_i$, this constraint can be written as
\begin{equation}
\boldsymbol{c}_i^\top \boldsymbol{a}=u_i,
\qquad
\boldsymbol{c}_i:=\boldsymbol{p}(\boldsymbol{0}).
\end{equation}

The CWLS reconstruction is therefore defined by the constrained minimization problem
\begin{equation}
\label{eq:cwls_min_problem}
\boldsymbol{a}_i
=
\arg\min_{\boldsymbol{a}}
\frac12 \|\mathbf W_i(\mathbf P_i\boldsymbol{a}-\boldsymbol{u}_{\mathcal S_i})\|_2^2
\quad
\text{subject to}
\quad
\boldsymbol{c}_i^\top \boldsymbol{a}=u_i.
\end{equation}

Introduce the normal-equation quantities
\begin{equation}
\mathbf A_i:=\mathbf P_i^\top \mathbf W_i^2 \mathbf P_i,
\qquad
\mathbf b_i:=\mathbf P_i^\top \mathbf W_i^2 \boldsymbol{u}_{\mathcal S_i}.
\end{equation}
The first-order optimality conditions lead to the KKT system
\begin{equation}
\label{eq:cwls_kkt_star}
\begin{bmatrix}
\mathbf A_i & \boldsymbol{c}_i \\
\boldsymbol{c}_i^\top & 0
\end{bmatrix}
\begin{bmatrix}
\boldsymbol{a}_i \\
\lambda_i
\end{bmatrix}
=
\begin{bmatrix}
\mathbf b_i \\
u_i
\end{bmatrix},
\end{equation}
where $\lambda_i$ is the Lagrange multiplier associated with the interpolation constraint.

Provided that the patch stencil is well posed, namely that $\mathbf P_i$ has full column rank and hence $\mathbf A_i$ is nonsingular, the coefficient vector $\boldsymbol{a}_i$ can be written explicitly through the Schur-complement formula:
\begin{equation}
\label{eq:ai_schur}
\boldsymbol{a}_i
=
\mathbf A_i^{-1}\mathbf b_i
-
\frac{\boldsymbol{c}_i^\top \mathbf A_i^{-1}\mathbf b_i-u_i}
{\boldsymbol{c}_i^\top \mathbf A_i^{-1}\boldsymbol{c}_i}
\,\mathbf A_i^{-1}\boldsymbol{c}_i.
\end{equation}
This expression shows that the local approximation depends linearly on the nodal values in the patch. The second term arises from the patch-star constraint.

\subsection{Local reconstruction shape functions}

Because the CWLS system is linear in the patch nodal values, the local approximant can be rewritten in nodal-value form:
\begin{equation}
\label{eq:cwls_local_shape_expansion}
U_i(\boldsymbol{x})
=
\sum_{k\in\mathcal S_i}\Psi_k^{(i)}(\boldsymbol{x})\,u_k,
\end{equation}
where $\Psi_k^{(i)}$ are the local reconstruction shape functions associated with patch $\omega_i$.

To derive them explicitly, define
\begin{equation}
\mathbf Z_i:=\mathbf A_i^{-1}\mathbf P_i^\top \mathbf W_i^2,
\end{equation}
and let $\mathbf e_k$ denote the unit sampling vector on the patch node set $\mathcal S_i$, i.e. the vector whose $k$-th patch entry is one and all others are zero. Then each local reconstruction mode can be written as
\begin{equation}
\boldsymbol{m}_k^{(i)}
=
\mathbf Z_i \mathbf e_k
-
\frac{\boldsymbol{c}_i^\top \mathbf Z_i \mathbf e_k-\delta_{ik}}
{\boldsymbol{c}_i^\top \mathbf A_i^{-1}\boldsymbol{c}_i}
\,\mathbf A_i^{-1}\boldsymbol{c}_i,
\end{equation}
and hence
\begin{equation}
\label{eq:mk_closed}
\Psi_k^{(i)}(\boldsymbol{x})
=
\boldsymbol{p}(\boldsymbol{\xi})^\top \boldsymbol{m}_k^{(i)}.
\end{equation}

Substituting these local basis functions into \eqref{eq:cwls_local_shape_expansion} yields the reconstructed patch field directly in terms of the original nodal unknowns. In particular, the interpolation constraint at the patch-star node is built into the local basis:
\begin{equation}
\Psi_k^{(i)}(\boldsymbol{x}_i)=\delta_{ik},
\qquad
U_i(\boldsymbol{x}_i)=u_i.
\end{equation}
Thus each local approximant interpolates the nodal value at its own patch center while remaining least-squares optimal, in the weighted sense, over the rest of the patch stencil.

\subsection{Derivative evaluation}

Since the local approximant is represented by a single polynomial field on the patch, its derivatives are obtained directly by differentiating the polynomial basis. From \eqref{eq:cwls_local_shape_expansion},
\begin{equation}
\label{eq:cwls_derivative}
D^\alpha U_i(\boldsymbol{x})
=
\sum_{k\in\mathcal S_i} D^\alpha \Psi_k^{(i)}(\boldsymbol{x})\,u_k.
\end{equation}
Using the scaled coordinate $\boldsymbol{\xi}$, one has
\begin{equation}
D^\alpha \Psi_k^{(i)}(\boldsymbol{x})
=
\rho_i^{-|\alpha|}
\bigl(D^\alpha_{\boldsymbol{\xi}}\boldsymbol{p}(\boldsymbol{\xi})\bigr)^\top
\boldsymbol{m}_k^{(i)}.
\end{equation}
These derivatives are used directly in the MiGFEM assembly, both for weak-form terms and for the leading Leibniz derivative employed in strong-form collocation.

\subsection{Kronecker property of the global MiGFEM approximation}

The local constraint $U_i(\boldsymbol{x}_i)=u_i$ is precisely what allows the global MiGFEM field to retain the Kronecker delta property at mesh nodes. Recall that the global approximation is constructed by partition-of-unity blending:
\begin{equation}
u^h(\boldsymbol{x})=\sum_{i\in\mathcal N} N_i(\boldsymbol{x})\,U_i(\boldsymbol{x}),
\end{equation}
where the standard finite element PoU basis satisfies
\begin{equation}
N_i(\boldsymbol{x}_j)=\delta_{ij}.
\end{equation}
Evaluating $u^h$ at a mesh node $\boldsymbol{x}_j$ gives
\begin{equation}
u^h(\boldsymbol{x}_j)
=
\sum_{i\in\mathcal N} N_i(\boldsymbol{x}_j)\,U_i(\boldsymbol{x}_j)
=
U_j(\boldsymbol{x}_j)
=
u_j.
\end{equation}
Therefore the global MiGFEM approximation interpolates the nodal values exactly:
\begin{equation}
\label{eq:kronecker_appendix}
u^h(\boldsymbol{x}_j)=u_j.
\end{equation}

This property is important both conceptually and computationally. Conceptually, it shows that the mesh-intrinsic enrichment does not alter the nodal meaning of the unknowns: they remain physical nodal values, exactly as in standard FEM. Computationally, it simplifies the treatment of essential boundary conditions and preserves direct compatibility with nodal finite element data structures.

In summary, the CWLS reconstruction provides a general patchwise local approximation on unstructured meshes, while the interpolation constraint at the patch-star node guarantees that the resulting MiGFEM basis remains nodally interpolatory. This combination is one of the key ingredients that allows MiGFEM to enrich the approximation without introducing any additional degrees of freedom.

\section{Derivation details for inter-element jump bounds}
\label{app:jump_decay_proof}

This appendix collects the algebraic steps behind \eqref{eq:exact_jump_cancellation} and the asymptotic estimate \eqref{eq:asymp_cp_rate} in Section~\ref{sec:asymp_jump_taylor}. The argument is local to an internal interface $\Gamma$ and uses the MiGFEM representation $u^h=\sum_i N_i U_i$, the PoZ identity \eqref{eq:poz_jump}, the reconstruction error bound \eqref{eq:ei_mls_bound}, and the standard inverse estimate $\|D^\beta N_i\|_{L^\infty(E)}\le C_N h^{-|\beta|}$ for linear FE shape functions.

\subsection{Exact cancellation under interface coherence}

For $k$-th order derivative coherence on $\Gamma$, common traces $U_\Gamma$ and $D^\mu U_\Gamma$ exist for $|\mu|\le k$. Substituting into \eqref{eq:tag21} yields the factorization \eqref{eq:jump_factorization}. Each sum $\sum_{i\in\mathcal I(\Gamma)}\llbracket D^\beta N_i\rrbracket_\Gamma$ vanishes by \eqref{eq:poz_jump}, hence $\llbracket D^\alpha u^h\rrbracket_\Gamma=0$ for $|\alpha|\le k+1$.

\subsection{Split into exact field and reconstruction error}

With $U_i=u+e_i$ on each patch, expansion of \eqref{eq:tag21} gives the decomposition $\llbracket D^\alpha u^h\rrbracket_\Gamma = J_u(\alpha;\Gamma) + J_e(\alpha;\Gamma)$ with
\begin{equation}
\label{eq:jump_u_factorization}
J_u(\alpha;\Gamma)
=
u\big|_\Gamma
\sum_{i\in\mathcal I(\Gamma)}
\llbracket D^\alpha N_i\rrbracket_\Gamma
+
\sum_{0<\beta<\alpha}
\binom{\alpha}{\beta}
\left(D^{\alpha-\beta}u\right)\big|_\Gamma
\sum_{i\in\mathcal I(\Gamma)}
\llbracket D^\beta N_i\rrbracket_\Gamma ,
\end{equation}
and $J_e(\alpha;\Gamma)$ defined by the same formula with $u$ replaced by $e_i$. PoZ annihilates $J_u$, leaving \eqref{eq:jump_error_only}.

\subsection{Asymptotic estimate \eqref{eq:asymp_cp_rate}}

For each $i\in\mathcal I(\Gamma)$,
\[
\|\llbracket D^\beta N_i\rrbracket_\Gamma\|_{L^\infty(\Gamma)}
\le
\|D^\beta N_i\|_{L^\infty(E^+)}+\|D^\beta N_i\|_{L^\infty(E^-)}
\le C_N h^{-|\beta|}
\]
by the standard inverse estimate $\|D^\beta N_i\|_{L^\infty(E)}\le C_N h^{-|\beta|}$. For traces of $e_i$, \eqref{eq:ei_mls_bound} gives $\|D^\mu e_i\|_{L^\infty(\omega_i)}\le C_e h^{p+1-|\mu|}$. Each summand in $J_e$ therefore satisfies
\[
\|\llbracket D^\beta N_i\rrbracket_\Gamma\,(D^{\alpha-\beta}e_i)|_\Gamma\|_{L^\infty(\Gamma)}
\le C\, h^{p+1-|\alpha|},
\]
and the same bound holds for the $\beta=\alpha$ terms with $e_i|_\Gamma$. The set $\mathcal I(\Gamma)$ has cardinality bounded independently of $h$ on shape-regular meshes, and the number of multi-indices $\beta$ in the Leibniz sum depends only on $\alpha$. Summation yields \eqref{eq:asymp_cp_rate} with $\widehat C$ depending on $p$, $|\mathcal I(\Gamma)|$, and the mesh and reconstruction parameters in Section~\ref{sec:asymp_jump_taylor}, but not on $h$.

\subsection{Relation to MLS references}

The bound \eqref{eq:ei_mls_bound} is the standard $C^{p+1}$ reproduction estimate for weighted least-squares on sufficiently dense, shape-regular stencils; see \cite{ref41,ref42,Wendland2004}. The unisolvence and Wendland-type error bound apply uniformly on the mesh family under the shape-regular and stencil conditions of Section~\ref{sec:asymp_jump_taylor}. No additional regularity beyond $C^{p+1}(\Omega)$ is imposed on $u$.

\section{Leibniz cross-term remainder}
\label{app:leibniz_remainder}

This appendix summarizes the size of the Leibniz remainder $R^{\alpha}$ introduced in Section~\ref{sec:leibniz_leading_derivative}. From the definition
\[
R^\alpha
=
\sum_{0<\beta\le\alpha}\binom{\alpha}{\beta}
\sum_{i\in\mathcal I(E)}D^\beta N_i\,D^{\alpha-\beta}U_i,
\]
write $U_i=u+e_i$ and separate the $u$- and $e_i$-terms:
\[
R^\alpha
=
\underbrace{\sum_{0<\beta\le\alpha}\binom{\alpha}{\beta}
\bigl(D^{\alpha-\beta}u\bigr)
\sum_{i\in\mathcal I(E)}D^\beta N_i}_{\text{PoZ $\Rightarrow\; 0$}}
\;+\;
\sum_{0<\beta\le\alpha}\binom{\alpha}{\beta}
\sum_{i\in\mathcal I(E)}D^\beta N_i\,D^{\alpha-\beta}e_i.
\]
The $u$-term vanishes because $\sum_{i\in\mathcal I(E)}D^\beta N_i\equiv 0$ elementwise by the PoZ identity \eqref{eq:poz_def}. Hence the remainder is driven entirely by the local reconstruction errors:
\begin{equation}
R^\alpha
=
\sum_{0<\beta\le\alpha}\binom{\alpha}{\beta}
\sum_{i\in\mathcal I(E)}D^\beta N_i\,D^{\alpha-\beta}e_i.
\end{equation}

A first consequence is that whenever the chosen local reconstruction space represents the exact solution on each patch, one has $e_i\equiv 0$ and therefore
\begin{equation}
\label{eq:Ralpha_zero_poly}
R^\alpha\equiv 0,
\qquad |\alpha|\le p.
\end{equation}
In that case, the leading Leibniz derivative coincides exactly with the full Leibniz derivative:
\begin{equation}
\label{eq:ld_equals_leibniz_poly}
\widetilde D^\alpha u^h=D^\alpha u^h.
\end{equation}
Polynomial reproduction is the simplest representative case, but the same conclusion holds for any enrichment space that captures the local solution exactly.

For general smooth solutions $u\in C^{p+1}(\Omega)$, the local reconstruction estimate
\begin{equation}
\|D^{\mu}e_i\|_{L^\infty(\omega_i)}\le C_e h^{p+1-|\mu|},
\qquad |\mu|\le p,
\end{equation}
holds on shape-regular patches \cite{ref41,ref42,Wendland2004}. Combined with the inverse bound $\|D^\beta N_i\|_{L^\infty(E)}\le C_N h^{-|\beta|}$, each term in $R^\alpha$ is of order $h^{p+1-|\alpha|}$. Summing over the finitely many active indices and admissible multi-indices yields
\begin{equation}
\label{eq:Ralpha_asymp_order}
\|R^\alpha\|_{L^\infty(E)}=\mathcal O(h^{p+1-|\alpha|}),
\qquad |\alpha|\le p.
\end{equation}
Therefore, omitting $R^{\alpha}$ in strong-form collocation does not reduce the asymptotic order of the local reconstruction. Moreover, if the enrichment basis captures the dominant local structure of the exact solution, the constant hidden in \eqref{eq:Ralpha_asymp_order} may be much smaller than in the purely polynomial case.

\section{Nodal-injection diagnostic for the biharmonic trial space}
\label{app:biharmonic_order_gap}

\noindent\textbf{Objective.}
Isolate the approximation capability of the mesh-intrinsic trial space for the biharmonic problem by injecting exact nodal values directly, without a global PDE solve. This test determines whether the observed order gap in the full strong-form solve (Section~\ref{sec:example_group_c_biharmonic}) is caused by the local approximation or by the collocation mechanism.

\noindent\textbf{Setup.}
The manufactured solution is
\[
u(x,y)=\sin(\pi x)\sin(\pi y), \qquad (x,y)\in(0,1)^2,
\]
with the corresponding body force and boundary data chosen to make this field exact. The biharmonic operator is evaluated by the leading Leibniz derivative. A sequence of uniformly refined structured meshes ($11\times11$ through $81\times81$) is used with polynomial orders $p=4$ and $p=6$, on both regular ($\delta=0$) and strongly perturbed ($\delta=1.0h$) meshes. Exact nodal values are prescribed directly and no global PDE solve is carried out, so the only source of error is the patchwise reconstruction itself.

\subsection{Nodal injection: approximation-space diagnostic}

\begin{figure}[htbp]
  \centering
  \begin{minipage}{0.90\textwidth}
    \centering
    \includegraphics[width=\linewidth]{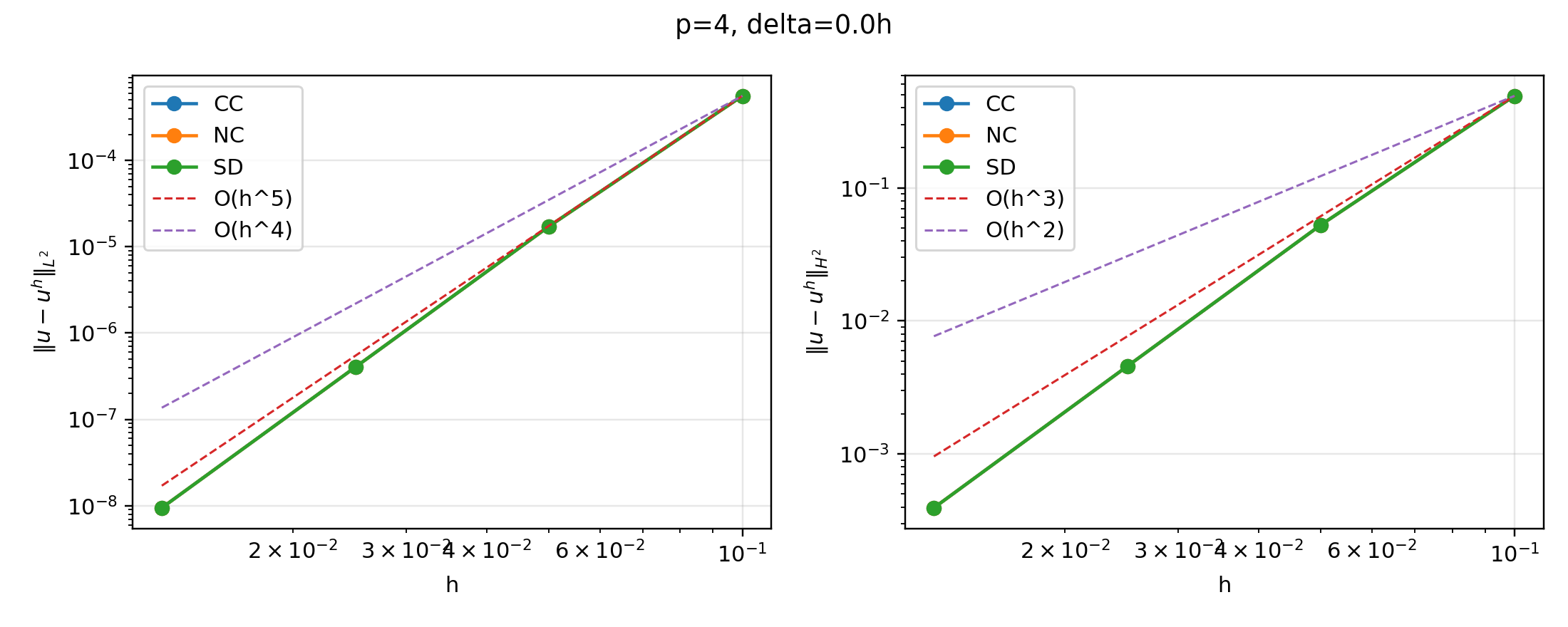}
    \captionline{(a) $p=4$, regular mesh ($\delta=0.0h$)}
  \end{minipage}
  
  \vspace{0.2cm}
  \begin{minipage}{0.90\textwidth}
    \centering
    \includegraphics[width=\linewidth]{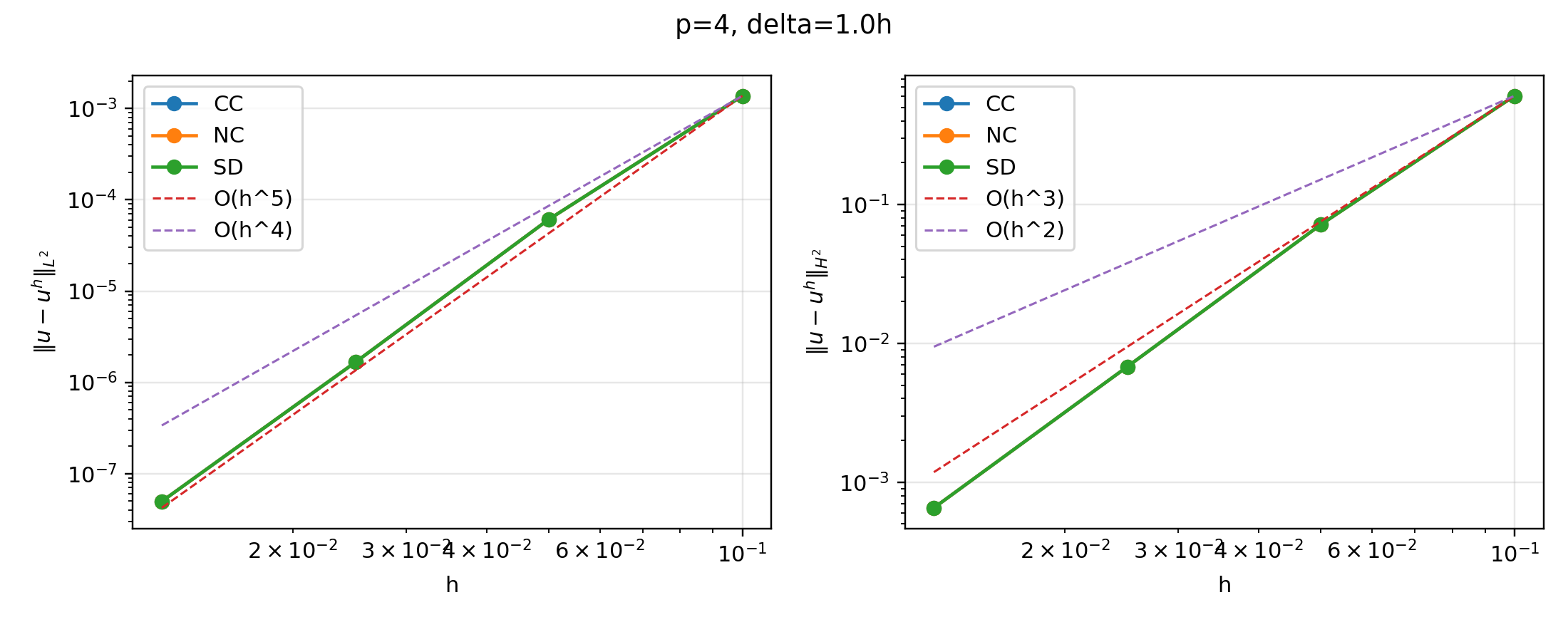}
    \captionline{(b) $p=4$, strong perturbation ($\delta=1.0h$)}
  \end{minipage}
  \caption{Nodal-injection diagnostic for the biharmonic benchmark at $p=4$.}
  \label{fig:biharmonic_nodal_injection}
  \end{figure}
  
  \begin{figure}[htbp]
  \centering
  \begin{minipage}{0.90\textwidth}
    \centering
    \includegraphics[width=\linewidth]{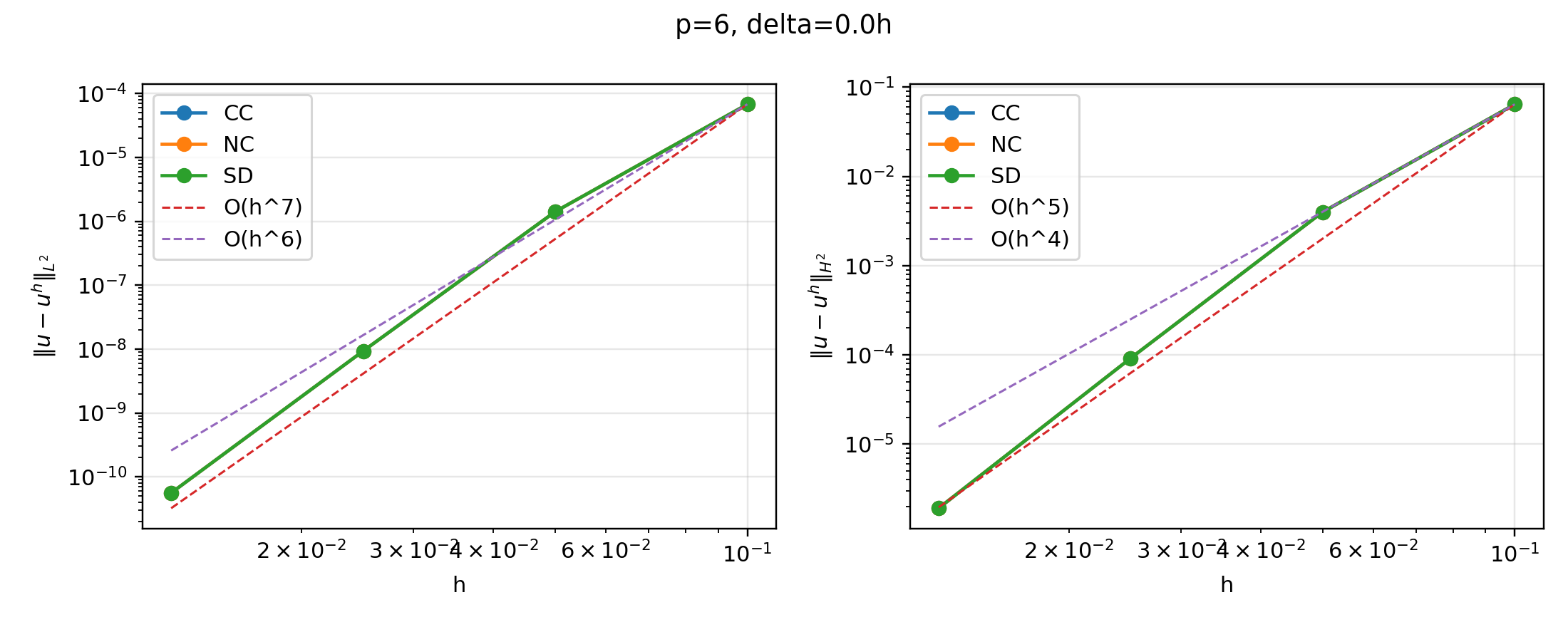}
    \captionline{(a) $p=6$, regular mesh ($\delta=0.0h$)}
  \end{minipage}
  
  \vspace{0.2cm}
  \begin{minipage}{0.90\textwidth}
    \centering
    \includegraphics[width=\linewidth]{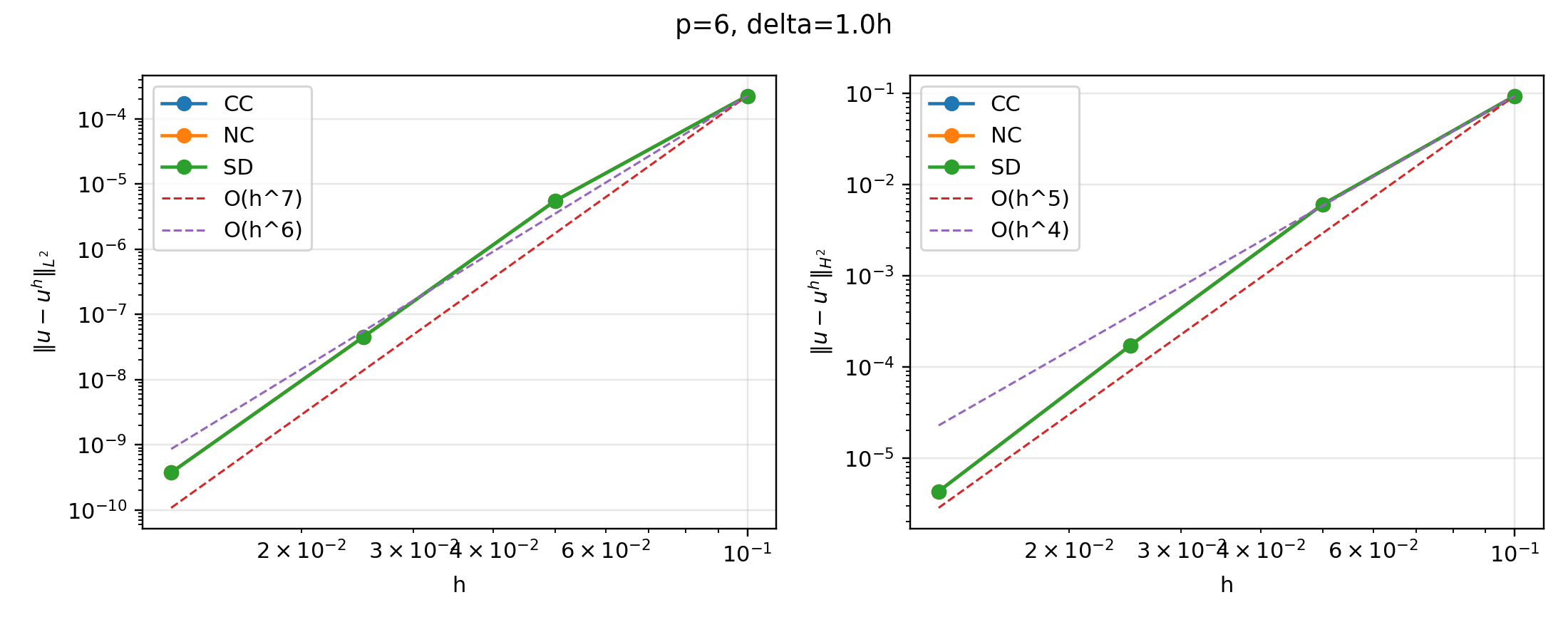}
    \captionline{(b) $p=6$, strong perturbation ($\delta=1.0h$)}
  \end{minipage}
  \caption{Nodal-injection diagnostic for the biharmonic benchmark at $p=6$.}
  \label{fig:biharmonic_nodal_injection_p6}
  \end{figure}

\noindent\textbf{Results.}
For nodal injection, the discrete field is obtained directly from exact nodal values---no PDEs solved, so the only source of error is the reconstruction itself. Rates $p+1$ and $p-1$ are observed for $L_h^2$ and $H_h^2$ errors respectively, and are consistent with the expected interpolation order $\mathcal{O}(h^{p+1})$ on both regular and strongly perturbed meshes. The results indicate that the mesh-intrinsic reconstruction and the leading Leibniz evaluation do not introduce any additional loss of accuracy at the approximation level. Kapl et al.\ \cite{iga-c-2} report analogous suboptimal convergence for isogeometric collocation on the biharmonic equation over multi-patch planar domains; Auricchio et al.\ \cite{iga-c-1} analyze isogeometric collocation more generally and note reduced rates relative to Galerkin IGA in several settings. These observations and the nodal-injection diagnostic above collectively indicate that the order gap in the full solve is attributable to the collocation mechanism rather than the trial space.

\section{Comparison with standard \texorpdfstring{$p$}{p}-FEM and PUFEM}
\label{subsec:mi_pfem_pu_mini_case}

\noindent\textbf{Objective.}
Compare MiGFEM WG with same-order $p$-FEM and PUFEM at comparable vector DOFs (accuracy, conditioning proxy, runtime).

\noindent\textbf{Setup.}
MiGFEM WG ($p=4$, $s=3$), $P_4$ $p$-FEM, degree-four PUFEM; same field and meshes. Boundary treatment differs across methods (see discussion).

\begin{figure}[htbp]
\centering
\includegraphics[width=0.8\linewidth]{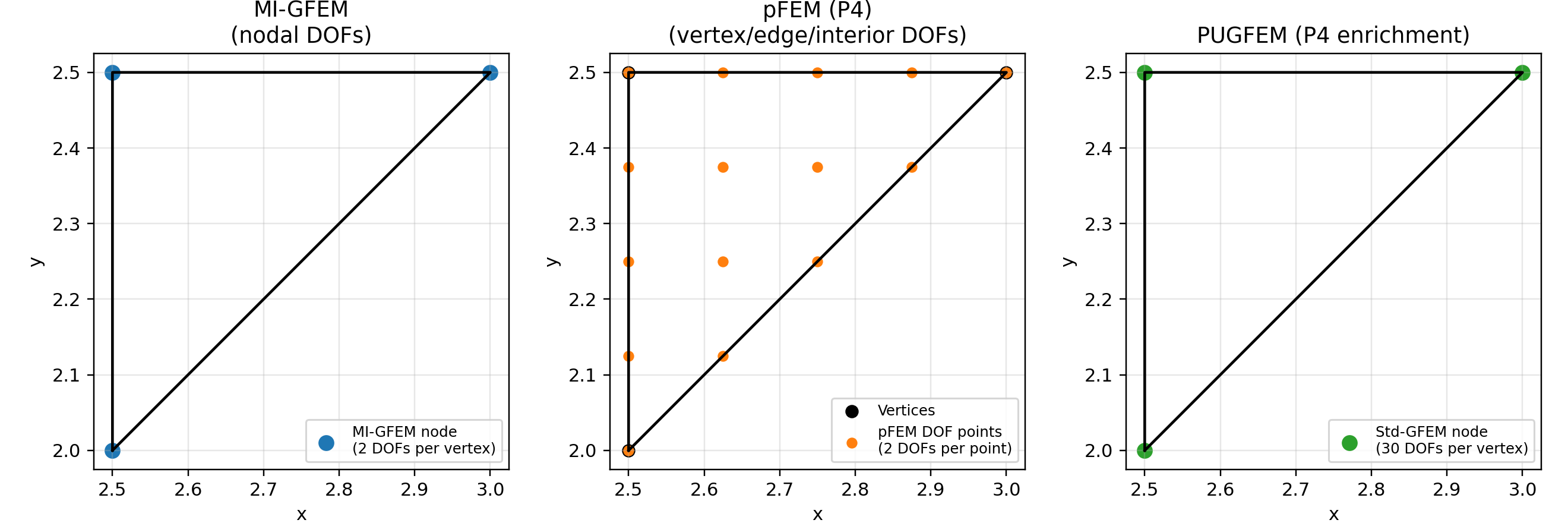}
\caption{Illustration of DOF layouts for MiGFEM WG (nodal DOFs), standard $P_4$ $p$-FEM (vertex/edge/interior DOFs), and PUFEM (enriched nodal DOFs) on a representative triangle.}
\label{fig:triangle_dofs_mi_pfem_std_p4}
\end{figure}

\begin{figure}[htbp]
\centering
\includegraphics[width=\linewidth]{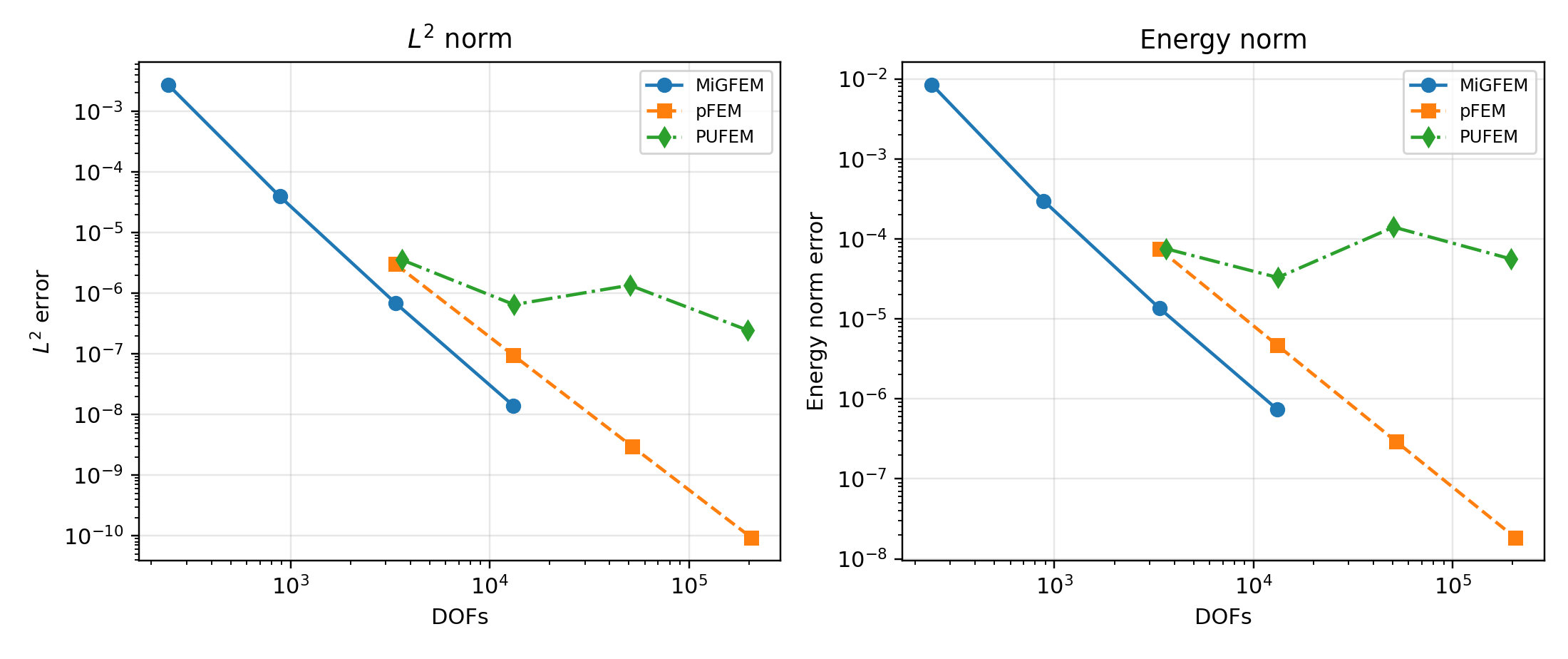}
\caption{DOF-based comparison of MiGFEM WG, standard $p$-FEM, and PUFEM on $(0,5)^2$ with $p=4$. Left: $L^2$ error versus total vector DOFs. Right: energy-norm error versus total vector DOFs.}
\label{fig:mi_pfem_pu_dof_comparison}
\end{figure}

\noindent\textbf{Results.}
Table~\ref{tab:mi_pfem_pu_elasticity} (representative mesh levels) and Figure~\ref{fig:mi_pfem_pu_dof_comparison} (error versus DOF trends) summarize the comparison.
On $11\times 11$, $L^2$ errors span roughly $\mathcal{O}(10^{-3})$ to $\mathcal{O}(10^{-6})$ with DOFs from $\mathcal{O}(10^2)$ to $\mathcal{O}(10^3)$; on $81\times 81$, all three methods reach $10^{-8}$--$10^{-11}$ in $L^2$.
At comparable DOF levels in this benchmark, MiGFEM shows favorable error levels relative to same-order $p$-FEM; this is reported as an empirical observation under the present discretization and boundary-treatment settings.
All three methods show formal high-order convergence on this smooth field. MiGFEM uses fewer DOFs per accuracy level but pays patch-reconstruction cost; $p$-FEM has the smallest $\kappa_\mathrm{diag}$; unstabilized PUFEM shows very large conditioning (stabilized GFEM variants such as SGFEM \cite{ref43} can mitigate this). Absolute errors are not strictly comparable because boundary enforcement differs.

\begin{table}[H]
  \centering
  \caption{Weak-form elasticity on $(0,5)^2$ with $p=4$: MiGFEM WG, standard $p$-FEM, and PUFEM on structured meshes.}
  \label{tab:mi_pfem_pu_elasticity}
  \begingroup
  \small
  \setlength{\tabcolsep}{3pt}
  \renewcommand{\arraystretch}{1.05}
  \begin{tabular}{cccccccc}
    \toprule
    Mesh & Method & $N_\mathrm{dof}$ &
    $\|\mathbf u-\mathbf u^h\|_{L^2}$ &
    $\|\mathbf u-\mathbf u^h\|_{E}$ &
    $\kappa_\mathrm{diag}(K)$ & $t_{\mathrm{total}}$ (s) & $t_{\mathrm{total}}/1000\,N_\mathrm{dof}$ (s) \\
    \midrule
    $11\times 11$ &
    MiGFEM &
    242 &
    $2.763\times 10^{-3}$ &
    $8.486\times 10^{-3}$ &
    $3.748\times 10^{3}$ &
    $5.96$ &
    $24.6$ \\
    &
    $p$-FEM &
    3362 &
    $3.011\times 10^{-6}$ &
    $7.460\times 10^{-5}$ &
    $8.440$ &
    $3.51$ &
    $1.04$ \\
    &
    PUFEM &
    3630 &
    $3.540\times 10^{-6}$ &
    $7.525\times 10^{-5}$ &
    $2.682\times 10^{14}$ &
    $8.66$ &
    $2.39$ \\
    \midrule
    $21\times 21$ &
    MiGFEM &
    882 &
    $3.953\times 10^{-5}$ &
    $2.987\times 10^{-4}$ &
    $3.748\times 10^{3}$ &
    $3.34$ &
    $3.79$ \\
    &
    $p$-FEM &
    13122 &
    $9.511\times 10^{-8}$ &
    $4.700\times 10^{-6}$ &
    $8.440$ &
    $9.54\times 10^{1}$ &
    $7.27$ \\
    &
    PUFEM &
    13230 &
    $6.492\times 10^{-7}$ &
    $3.285\times 10^{-5}$ &
    $6.866\times 10^{16}$ &
    $2.91\times 10^{1}$ &
    $2.20$ \\
    \midrule
    $41\times 41$ &
    MiGFEM &
    3362 &
    $6.881\times 10^{-7}$ &
    $1.362\times 10^{-5}$ &
    $3.748\times 10^{3}$ &
    $1.24\times 10^{1}$ &
    $3.69$ \\
    &
    $p$-FEM &
    51842 &
    $2.980\times 10^{-9}$ &
    $2.944\times 10^{-7}$ &
    $8.440$ &
    $2.30\times 10^{2}$ &
    $4.44$ \\
    &
    PUFEM &
    50430 &
    $1.356\times 10^{-6}$ &
    $1.403\times 10^{-4}$ &
    $1.758\times 10^{19}$ &
    $4.07\times 10^{2}$ &
    $8.06$ \\
    \midrule
    $81\times 81$ &
    MiGFEM &
    13122 &
    $1.406\times 10^{-8}$ &
    $7.408\times 10^{-7}$ &
    $3.748\times 10^{3}$ &
    $1.79\times 10^{2}$ &
    $13.7$ \\
    &
    $p$-FEM &
    206082 &
    $9.318\times 10^{-11}$ &
    $1.841\times 10^{-8}$ &
    $8.440$ &
    $2.96\times 10^{3}$ &
    $14.4$ \\
    &
    PUFEM &
    196830 &
    $2.465\times 10^{-7}$ &
    $5.589\times 10^{-5}$ &
    $4.500\times 10^{21}$ &
    $1.71\times 10^{3}$ &
    $8.68$ \\
    \bottomrule
  \end{tabular}
  \endgroup
\end{table}

\subsubsection{WG post-processing: full vs.\ leading Leibniz derivative}
\label{sec:wg_leibniz_postprocess}

\noindent\textbf{Objective.}
On the smooth nonpolynomial field \eqref{eq:benchmark_u_trig}, the Leibniz remainder $R^\alpha$ in $D^\alpha u^h=\widetilde D^\alpha u^h+R^\alpha$ is $\mathcal{O}(h^{p+1-|\alpha|})$ when patches fit a $C^{p+1}$ solution (Appendix~\ref{app:leibniz_remainder}), so it should not dominate the mesh-refinement rate. This test asks whether \emph{weak Galerkin post-processing}---forming stresses or strain energy from the discrete displacement using either the full or the leading Leibniz derivative---shows the \emph{same} asymptotic convergence order for the energy error. It is not a polynomial patch test and does not require $\widetilde D^\alpha u^h$ and $D^\alpha u^h$ to agree pointwise.

\noindent\textbf{Setup.}
Plane-stress WG on $\Omega=(0,5)^2$ with \eqref{eq:benchmark_u_trig}; uniform refinement; local degrees $p=1,2,4$. After each solve, the same $\mathbf{u}^h$ is used; constitutive derivatives enter via full Leibniz or leading Leibniz only.

\begin{figure}[htbp]
\centering
\begin{minipage}{0.32\textwidth}
  \centering
  \includegraphics[width=\linewidth]{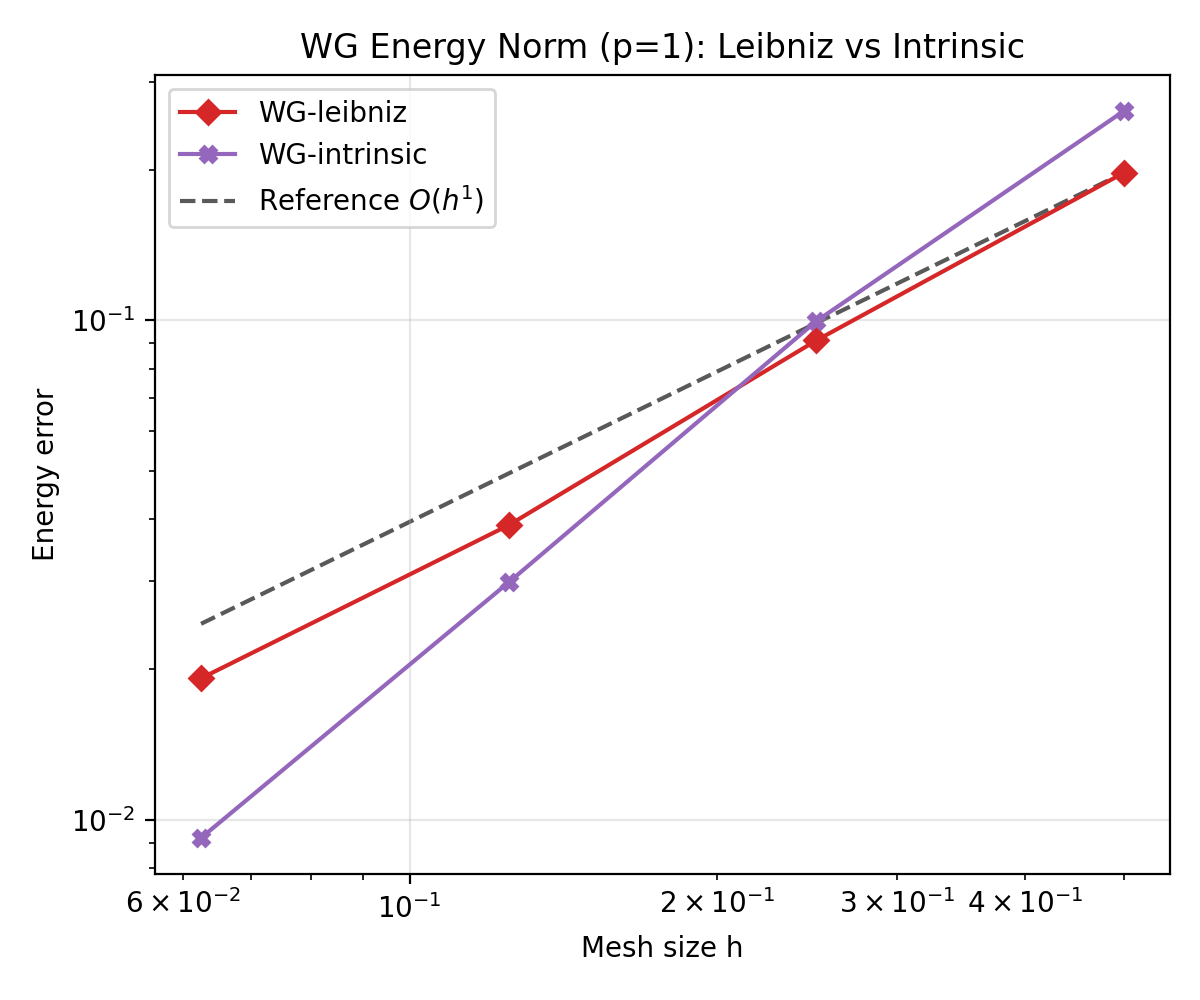}
  \captionline{(a) $p=1$}
\end{minipage}\hfill
\begin{minipage}{0.32\textwidth}
  \centering
  \includegraphics[width=\linewidth]{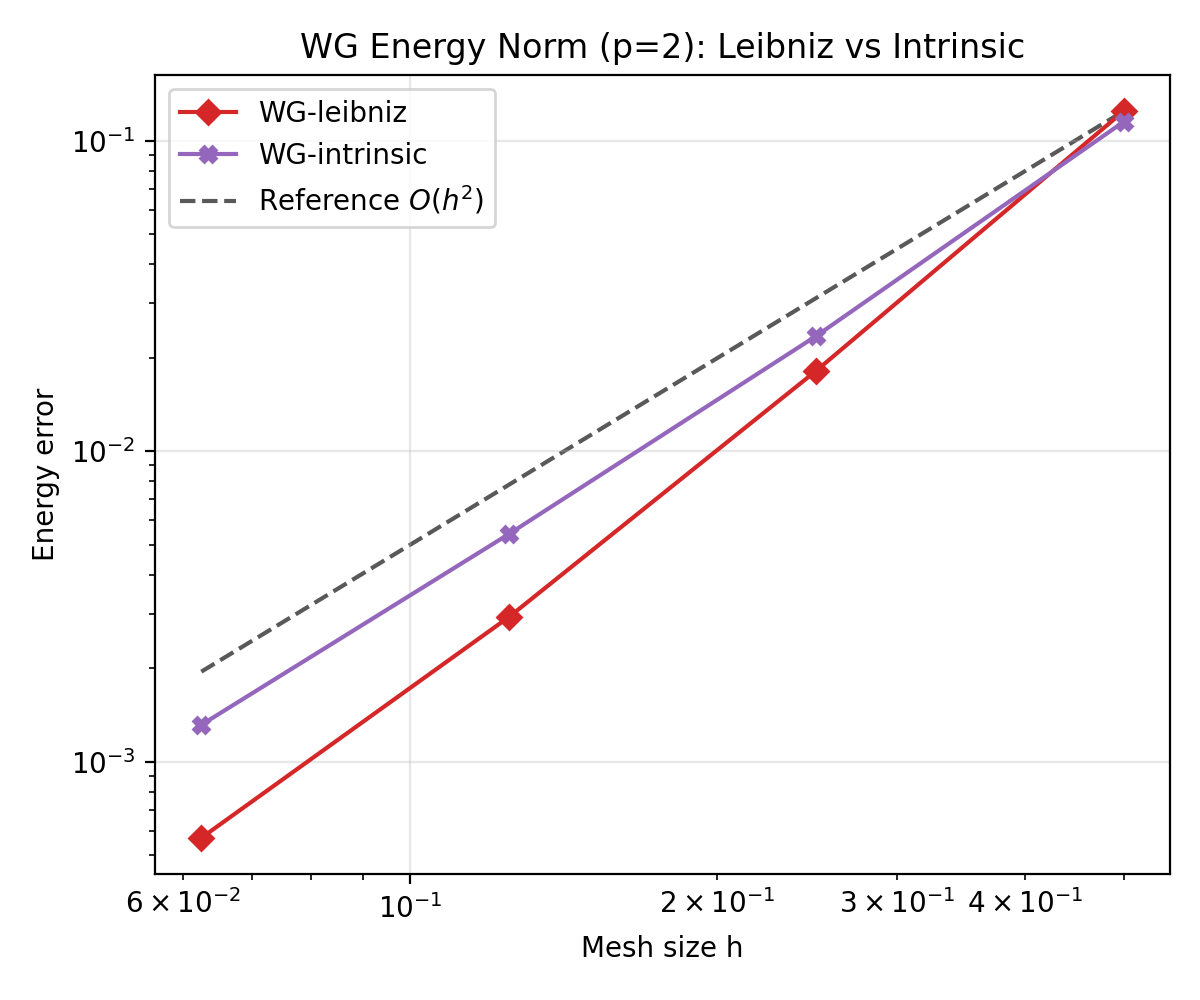}
  \captionline{(b) $p=2$}
\end{minipage}\hfill
\begin{minipage}{0.32\textwidth}
  \centering
  \includegraphics[width=\linewidth]{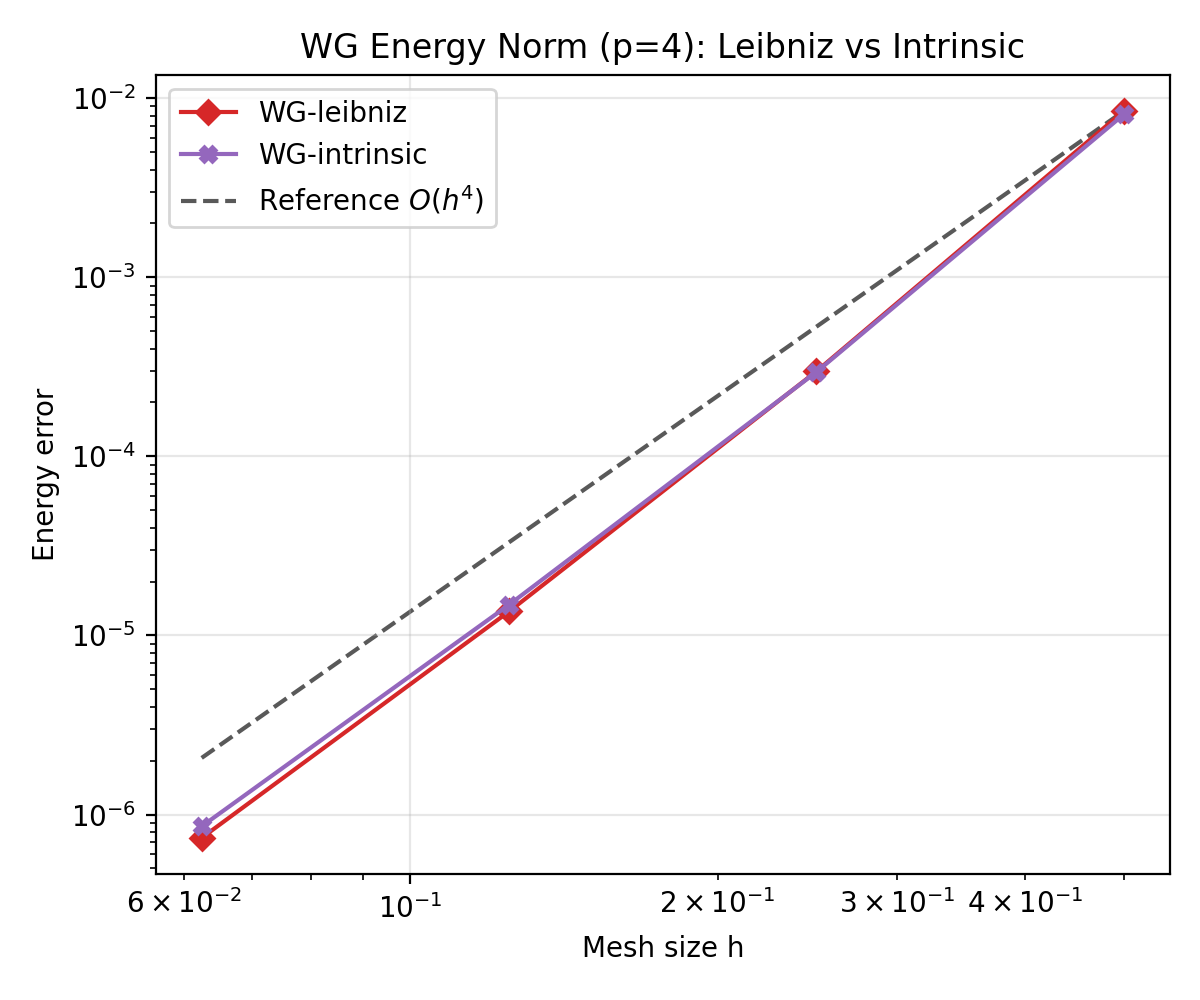}
  \captionline{(c) $p=4$}
\end{minipage}
\caption{Field \eqref{eq:benchmark_u_trig}: WG energy error vs.\ refinement, full vs.\ leading Leibniz post-processing ($p=1,2,4$).}
\label{fig:wg_energy_intrinsic_vs_leibniz_p124}
\end{figure}

\noindent\textbf{Results.}
Figure~\ref{fig:wg_energy_intrinsic_vs_leibniz_p124} plots energy error versus mesh size for both post-processing choices. Fitted slopes match in the asymptotic regime; coarse meshes at low $p$ can show a constant offset between the two curves. Dropping $R^\alpha$ in post-processing does not reduce the observed order relative to the full Leibniz evaluation.

\subsubsection{Stiffness conditioning}

\noindent\textbf{Objective.}
Compare growth of $\kappa_2(K_{II})$ under refinement for MiGFEM variants and standard FEM.

\paragraph{Setup.}
Same field \eqref{eq:benchmark_u_trig}; meshes $11\times 11$--$41\times 41$; MiGFEM NC, CC, SD, WG ($p=4$, $s=3$), $P_1$ and $P_4$ FEM. After Dirichlet elimination, $\kappa_2(K_{II})=\lambda_{\max}/\lambda_{\min}$ on the interior block.

\noindent\textbf{Results.}
\begin{table}[htbp]
\centering
\caption{Conditioning data for the elasticity benchmark on $(0,5)^2$ (interior block $K_{II}$).}
\label{tab:conditioning_kappa_raw}
\resizebox{\textwidth}{!}{%
\begin{tabular}{cccccccc}
\toprule
Mesh & $h$ & NC & CC & SD & WG & FEM $P_1$ & FEM $P_4$ \\
\midrule
$11\times11$ & 0.5000 & $4.04\times10^{1}$ & $5.63\times10^{1}$ & $5.74\times10^{1}$ & $6.04\times10^{1}$ & $6.22\times10^{1}$ & $3.15\times10^{3}$ \\
$21\times21$ & 0.2500 & $1.67\times10^{2}$ & $2.31\times10^{2}$ & $2.35\times10^{2}$ & $2.49\times10^{2}$ & $2.52\times10^{2}$ & $1.26\times10^{4}$ \\
$41\times41$ & 0.1250 & $6.73\times10^{2}$ & $9.30\times10^{2}$ & $9.46\times10^{2}$ & $1.01\times10^{3}$ & $1.01\times10^{3}$ & $5.04\times10^{4}$ \\
\bottomrule
\end{tabular}
}
\end{table}

\begin{figure}[htbp]
\centering
\includegraphics[width=0.6\linewidth]{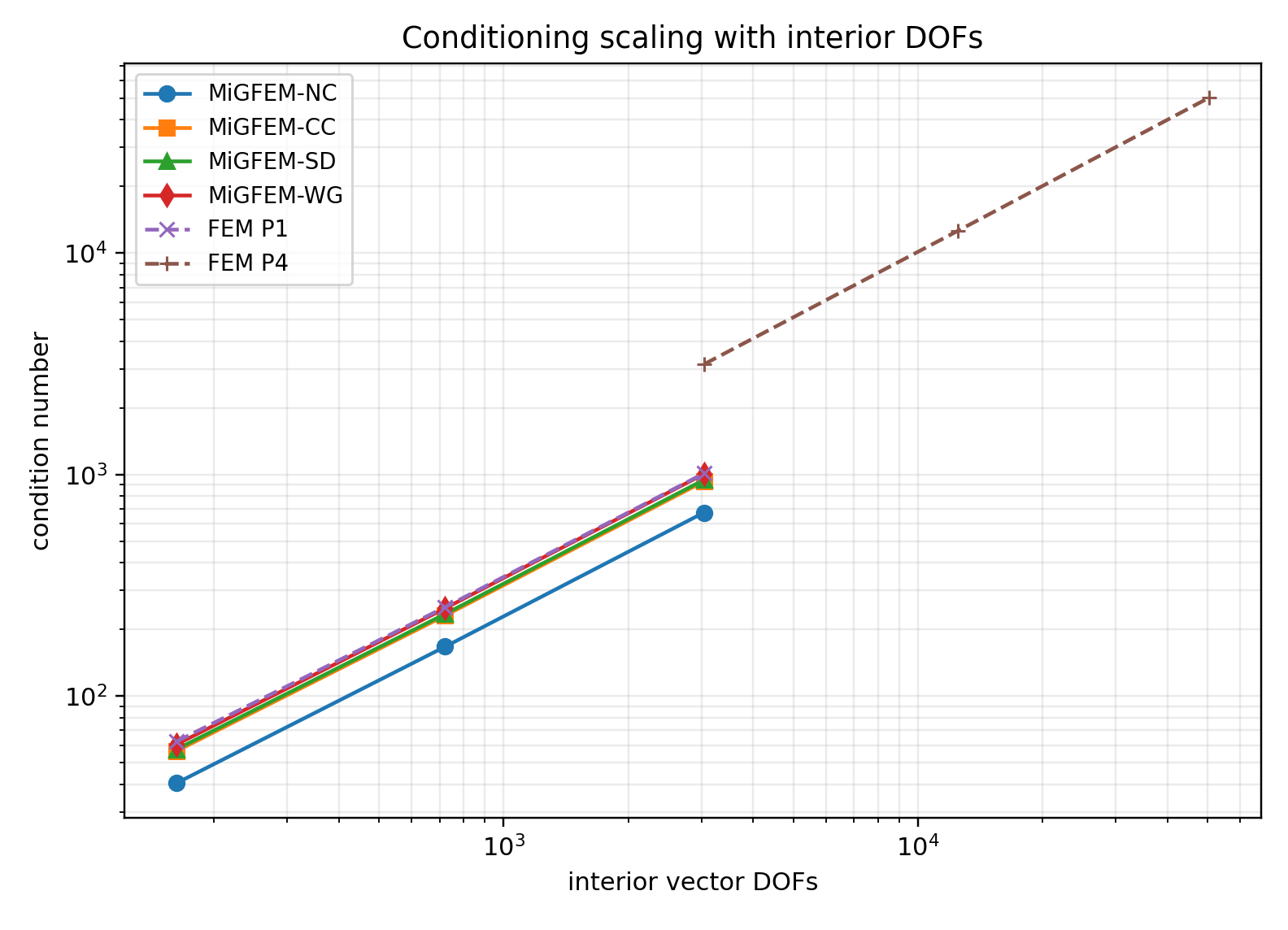}
\caption{Conditioning against interior vector DOFs for MiGFEM (NC/CC/SD/WG), FEM P1, and FEM P4.}
\label{fig:conditioning_vs_dof}
\end{figure}
Table~\ref{tab:conditioning_kappa_raw} and Figure~\ref{fig:conditioning_vs_dof} show monotone growth; MiGFEM and $P_1$ reach $\mathcal{O}(10^3)$ on the finest mesh, $P_4$ about $5\times10^4$. Fitted log--log slopes are near unity; differences are mostly prefactors. On this short mesh sequence, trends match expected $\mathcal{O}(h^{-2})$ scaling for second-order elliptic problems in 2D.

\subsubsection{Sensitivity to mesh perturbations}

\noindent\textbf{Objective.}
Assess robustness of MiGFEM and $p$-FEM under irregular meshes.

\noindent\textbf{Setup.}
Interior nodes perturbed by $0.3h$ or $0.5h$; boundary fixed; refinement to $81\times 81$; $p=4$.

\begin{figure}[htbp]
\centering
\includegraphics[width=\textwidth]{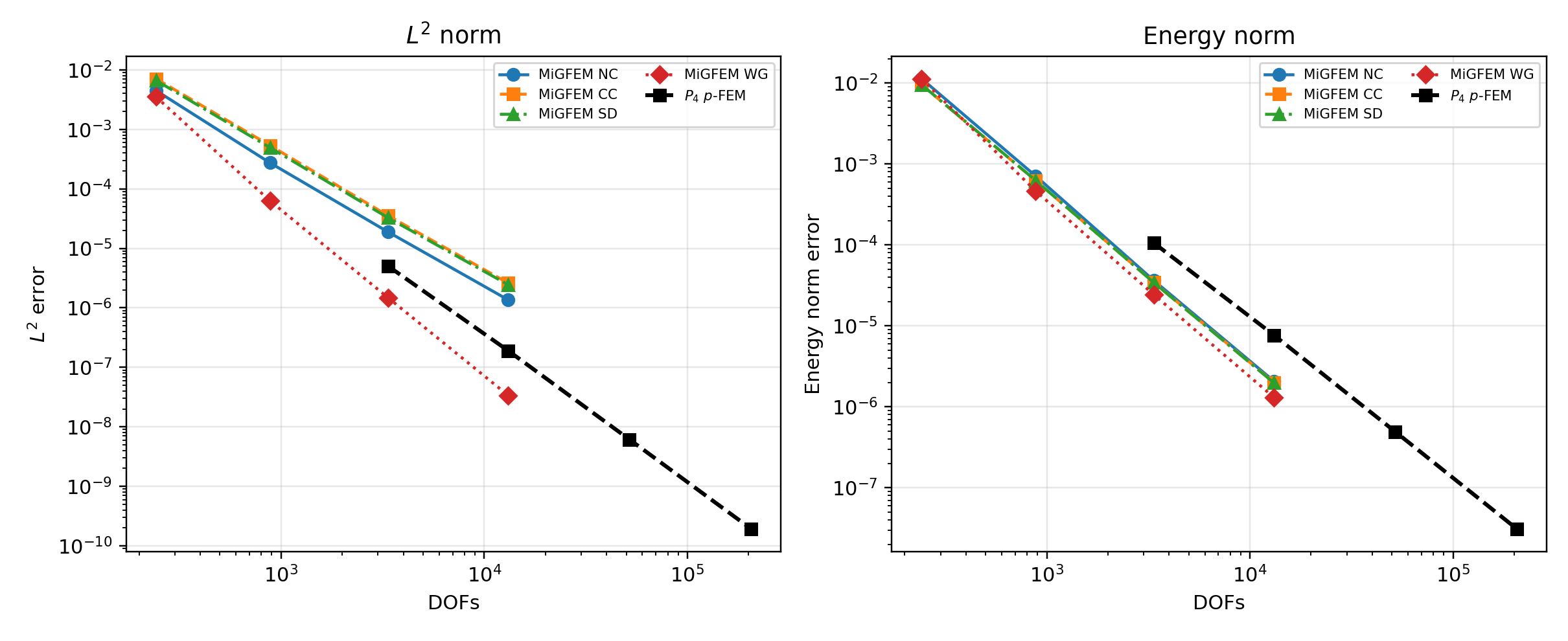}
\caption{Elasticity benchmark on perturbed meshes with $\delta=0.3h$: DOF-based $L^2$ and energy-norm errors for NC, CC, and SD, and WG, and $P_4$ $p$-FEM.}
\label{fig:weak_elasticity_irregular_mesh_p4_0p3}
\end{figure}

\begin{figure}[htbp]
\centering
\includegraphics[width=\textwidth]{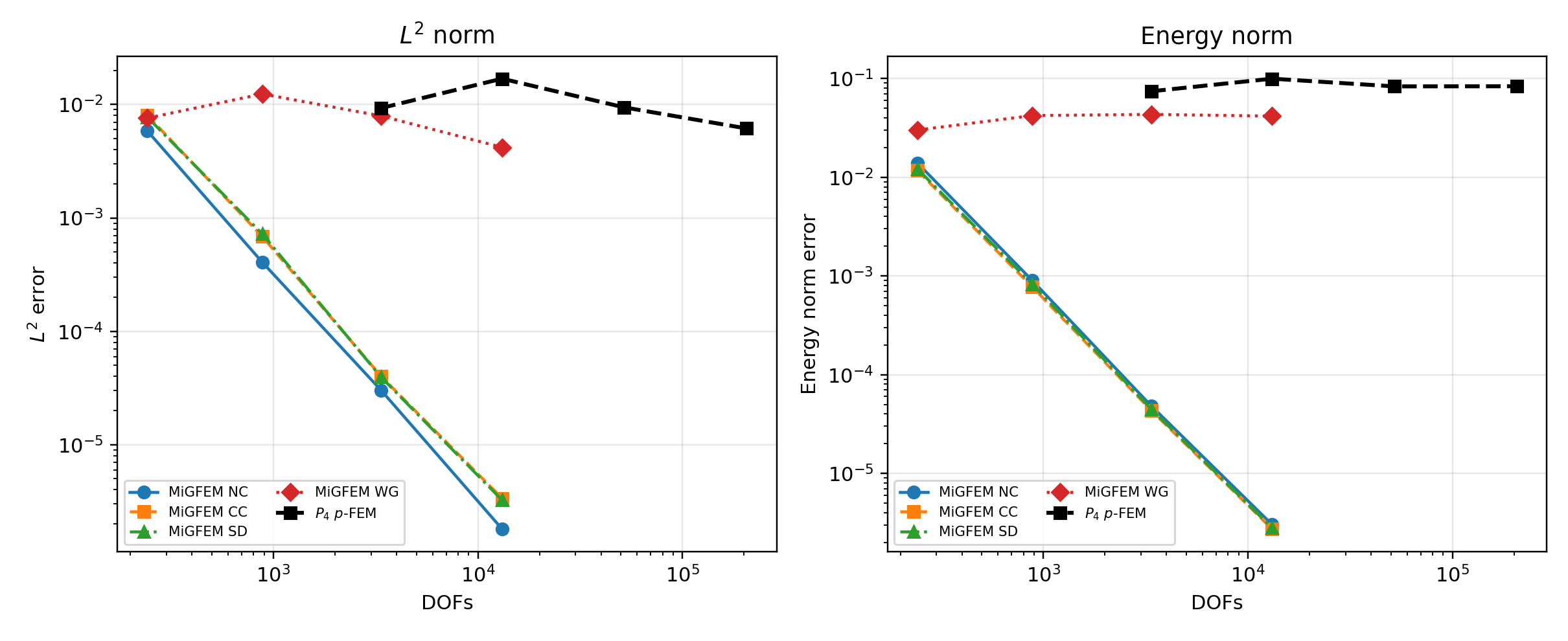}
\caption{Elasticity benchmark on perturbed meshes with increased distortion ($\delta=0.5h$): DOF-based $L^2$ and energy-norm errors for MiGFEM NC/CC/SD, MiGFEM WG, and $P_4$ $p$-FEM.}
\label{fig:weak_elasticity_irregular_mesh_p4_0p5}
\end{figure}

\noindent\textbf{Results.}
Figures~\ref{fig:weak_elasticity_irregular_mesh_p4_0p3} and \ref{fig:weak_elasticity_irregular_mesh_p4_0p5}. All MiGFEM variants show monotone $L^2$ and energy convergence; WG tracks $p$-FEM closely. CWLS reconstruction remains stable on the tested distortion levels; strong-form variants are not more sensitive than weak form in these norms.

\begin{figure}[htbp]
\centering
\begin{minipage}{0.49\textwidth}
  \centering
  \includegraphics[width=\linewidth]{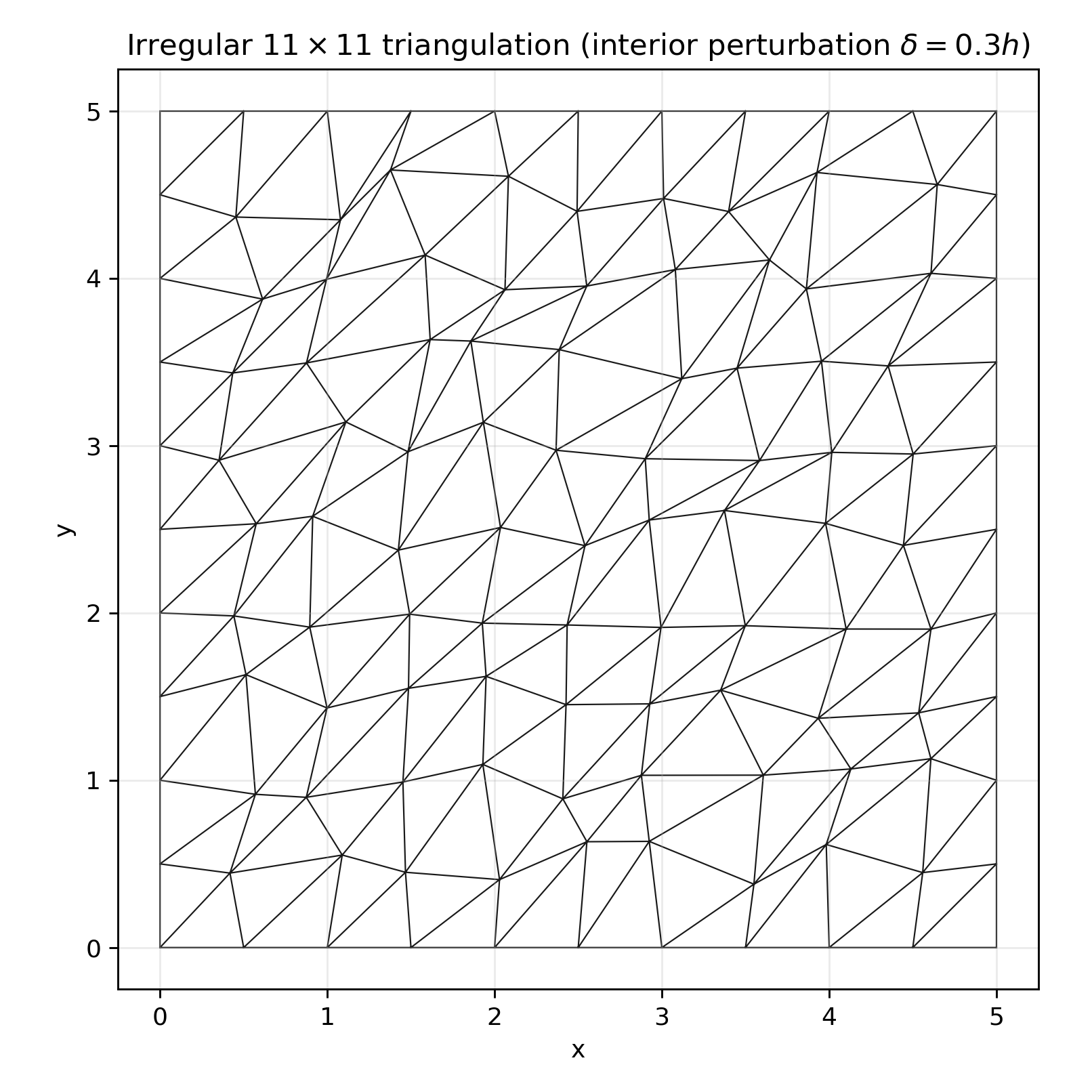}
  \captionline{(a) Perturbation $\delta=0.3h$}
\end{minipage}\hfill
\begin{minipage}{0.49\textwidth}
  \centering
  \includegraphics[width=\linewidth]{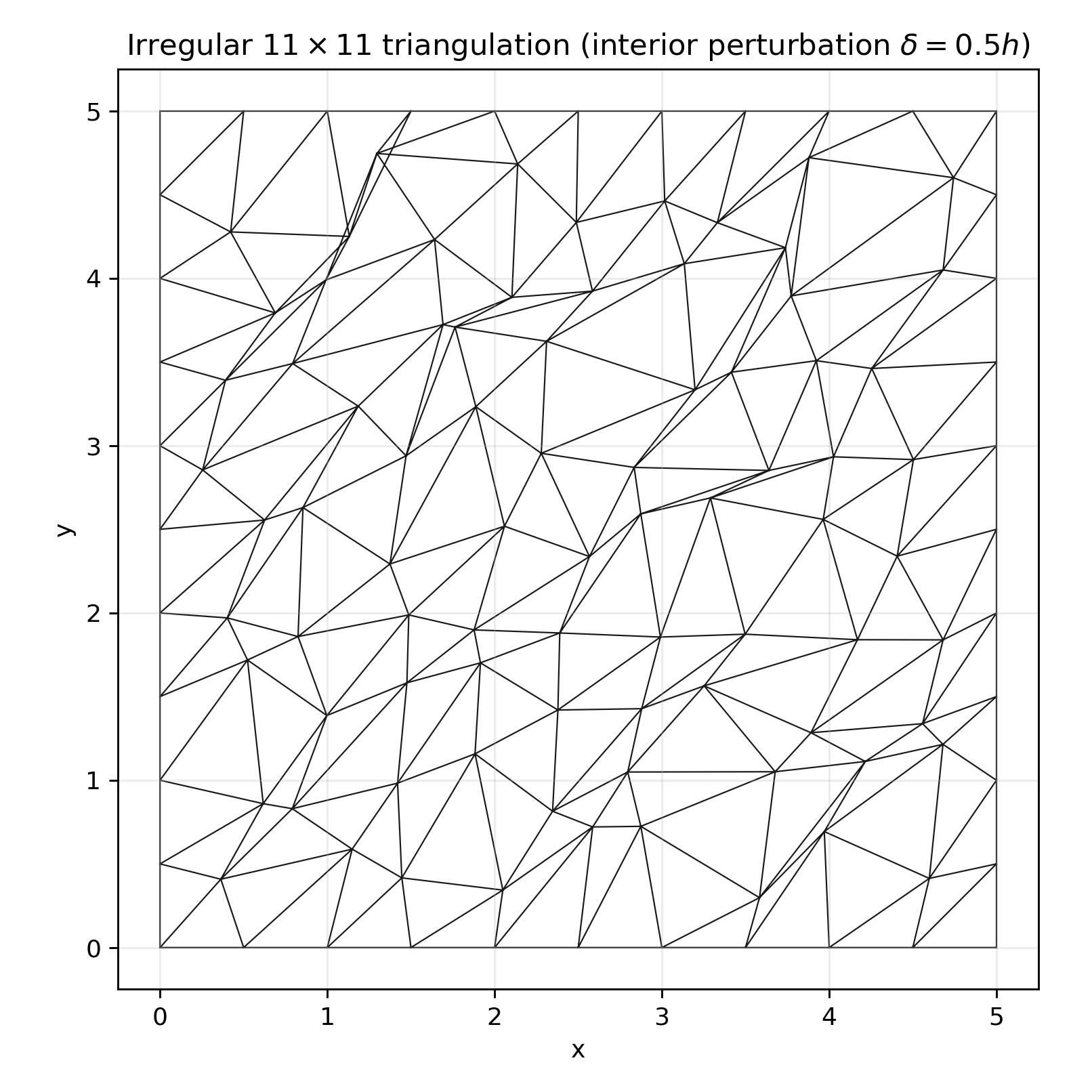}
  \captionline{(b) Perturbation $\delta=0.5h$}
\end{minipage}
\caption{Representative irregular $11\times 11$ triangulations on $(0,5)^2$ obtained by perturbing interior nodes by fractions $\delta \in \{0.3h, 0.5h\}$ of the local mesh size while keeping boundary nodes fixed.}
\label{fig:irregular_meshes_examples}
\end{figure}

\end{document}